\documentclass[11pt,bezier,amstex,a4wide]{article}  

\usepackage{color}              
\usepackage{graphicx,setspace}
\usepackage[]{amsmath}
\usepackage{amssymb}
\usepackage{hyperref}
\usepackage{enumerate, amsfonts, amsthm, bbm, psfrag}
\usepackage{fullpage}

\definecolor{MyDarkBlue}{rgb}{0,0.08,0.50}
\definecolor{BrickRed}{rgb}{0.65,0.08,0}

\hypersetup{
colorlinks=true,       
    linkcolor=MyDarkBlue,          
    citecolor=BrickRed,        
    filecolor=red,      
    urlcolor=cyan           
}

\newcommand{\funcDD}{\mathcal{V}}
\graphicspath{{Fig/}}

\newtheorem{theorem}{Theorem}[]

\newtheorem{proposition}[theorem]{Proposition}
\theoremstyle{definition}

\theoremstyle{definition}

\DeclareMathOperator{\sgn}{\text{\upshape sgn}}

\newcommand{\E}{\text{\itshape E}}
\newcommand{\Hup}{\text{\upshape H}}
\newcommand{\N}{\text{\itshape N}}
\newcommand{\V}{\text{\itshape V}}
\newcommand{\Ai}{\text{\upshape Ai}}
\newcommand{\Br}{\text{\upshape Br}}
\newcommand{\He}{\text{\upshape He}}

\newcommand{\Oup}{\text{\itshape O}}
\newcommand{\M}{\text{\itshape M}}
\newcommand{\U}{\text{\itshape U}}
\newcommand{\Hcal}{\mathcal{H}}
\newcommand{\Ocal}{\theta}

\newcommand{\Vcal}{\mathcal{V}}
\newcommand{\Mcal}{\mathcal{M}}

\newcommand{\e}{e}

\begin{document}

\title{\Large{\bf{Spectral gap of the Erlang A model\\ in the Halfin-Whitt regime}}}

\author{Johan S.H. van Leeuwaarden\footnotemark[1] \and
        Charles Knessl\footnotemark[2]
        }

\date{\today}

\maketitle

\footnotetext[1]{Eurandom, Technische Universiteit Eindhoven, Department of Mathematics and Computer Science,
P.O.\ Box 513, 5600 MB  Eindhoven, The Netherlands.
        Email: \texttt{j.s.h.v.leeuwaarden@tue.nl}
}

\footnotetext[2]{University of Illinois at Chicago, Department of Mathematics, Statistics and Computer Science, 815 South Morgan Street, Chicago, IL 60607-7045, USA.
        Email: \texttt{knessl@uic.edu}
}

\begin{abstract}
We consider a hybrid diffusion process that is a combination of two
Ornstein-Uhlenbeck processes with different
restraining forces. This process serves as the heavy-traffic approximation to the Markovian many-server queue with abandonments in the critical Halfin-Whitt regime. We obtain an expression for the Laplace transform of the time-dependent probability distribution, from which the spectral gap is explicitly characterized. The spectral gap gives the exponential rate of convergence to equilibrium. We further give various asymptotic results for the spectral gap, in the limits of small and large abandonment effects. It turns out that convergence to equilibrium becomes extremely slow for overloaded systems with small abandonment effects.
\end{abstract}

\section{Introduction}
Within the fields of stochastic processes and queueing theory, the Halfin-Whitt regime refers to a mathematical way of establishing economies-of-scale in many-server queueing systems like call centers (see \cite{gans}). The Halfin-Whitt regime in fact prescribes a scaling under which the many-server systems converge to limiting processes, which are for most systems diffusion processes. This paper deals with many-server systems in the Halfin-Whitt regime  with the additional feature that customers are impatient, and may abandon the system without being served. For such systems with abandonments, we are interested in the spectral gap, which is inversely related to the relaxation time or the speed at which a system reaches stationarity.  A large relaxation time in general indicates that replacing time-dependent characteristics by their stationary counterparts might lead to poor approximations. As it turns out, the rate at which customers renege (abandon the system) greatly influences the relaxation time.

In recent years, a large number of papers have dealt with the influence of reneging or abandonments on the system behavior (see, e.g.,~\cite{daitezcan,garnett,wardglynn,whitt04,whitt05,zeltyn04,zeltyn05}), and it is widely accepted that reneging is indeed one of the main factors driving the system performance. One of the key insights is that the system behavior strongly depends on whether it is stable or overloaded. By stable we mean that it can serve all customers, even if none of the customers would abandon the system. The spectral gap (or relaxation time) is also very different for stable or overloaded systems. In fact, we find that stable systems have a relatively short relaxation time, whereas the relaxation time of overloaded systems can become extremely large, particularly when the reneging rate is small.

The model we shall consider is the $M/M/s+M$ system, better known as the Erlang A model. This model is a standard Markovian many-server queueing system with Poisson arrivals, exponential service times, $s$ servers, and with the additional feature that customers that are waiting in the queue abandon the system after exponentially distributed reneging times. The queue length process in the Erlang A model, denoted by $(Q(t))_{t\geq 0}$, is a birth-death process. Whitt \cite{whitt04} (see also \cite{whitt06} and \cite{kangramanan}) derived a fluid approximation for the the steady-state behavior of the overloaded Erlang A model, and he further showed that a diffusion limit might provide refined approximations. Garnett, Mandelbaum and Reiman \cite{garnett}
proved a diffusion limit for the Erlang A model in the critical regime. In particular, they showed that under certain conditions a sequence of normalized queue length processes converges to a certain diffusion process $(X(t))_{t\geq 0}$. These conditions are in fact the ones that correspond to the Halfin-Whitt regime, in which the arrival rate $\lambda$ and the numbers of servers $s$ are scaled such that, while both $\lambda$ and $s$ increase toward infinity, the traffic intensity $\rho_0=\lambda/s$ approaches one, with
\begin{equation}\label{eqq1}
(1-\rho_0)\sqrt{s}\rightarrow \beta, \quad \beta\in(-\infty,\infty).
\end{equation}
The diffusion process $(X(t))_{t\geq 0}$ is a combination of two
Ornstein-Uhlenbeck (OU) processes with different
restraining forces, depending on whether the process is below or above zero.
The number of customers in the Erlang A model can be roughly expressed as $s+\sqrt{s}X(t)$ for $s$ sufficiently large. The diffusion process is generally easier to study than the birth-death process, and can thus be employed to obtain simple approximations for the system behavior.
The steady-state distribution of the diffusion can be easily obtained (see \eqref{steadyssty} below), but less is known about the time-dependent behavior. In this paper we shall present an explicit and asymptotic characterization of the spectral gap of  $(X(t))_{t\geq 0}$. The spectral gap of the diffusion process provides an understanding of the relaxation times for the Erlang A model in the Halfin-Whitt regime. Most importantly, we shall study in detail the impact on the spectral gap of the capacity parameter $\beta$ and the reneging rate $\eta$, which shall enhance our understanding of how the Erlang A model behaves for positive/negative $\beta$ and small/large values of $\eta$.

The diffusion process $(X(t))_{t\geq 0}$ also applies to the $G/M/s+M$ system, in the same asymptotic limit, which was proven by Whitt \cite{whitt05}.
Stochastic processes for more general systems with abandonments were obtained recently by Dai, He and Tezcan \cite{daitezcan} for the $G/Ph/s+M$ system. In this case, the limiting process is still a diffusion process, but it becomes multi-dimensional.
Zeltyn and Mandelbaum \cite{zeltyn05} derived approximations for the $M/M/n+G$ in the Halfin-Whitt regime.
In case of general service times, the limiting process is not even a diffusion process (see e.g.~\cite{reed,whitt05} for cases without reneging). Therefore, the one-dimensional diffusion process $(X(t))_{t\geq 0}$ strikes the proper balance between simplicity and tractability, while retaining the essential features of abandoning customers in many-server systems.

The Erlang A model is particularly interesting, as it incorporates three classical queueing systems as special cases. In the case of no reneging (with $\eta=0$) the Erlang A model reduces to the Erlang C model, or $M/M/s$ system. Halfin and Whitt \cite{halfinwhitt} established that the limiting process behaves as a Brownian motion above zero and an OU process below zero. In \cite{vlk} we have referred to this process as the {\it Halfin-Whitt diffusion}. For $\eta=1$ the Erlang A model becomes an infinite server queue or $M/M/\infty$ system, for which the stochastic-process limit is known to be an OU process \cite{IG}. For $\eta\to\infty$ the Erlang A model becomes the Erlang B model or $M/M/s/s$ system, in which case the stochastic-process limit is a reflected OU process (see \cite{robert,linetsky,wardglynn}). For the diffusion approximations, we show the $\eta\to\infty$ reduction in Section \ref{subsec:rou}. Our analysis of the spectral gap of $(X(t))_{t\geq 0}$ provides results for each of these three cases.

Mathematically, determining the transient distribution for the present diffusion process involves analyzing a Schr\"{o}dinger  type equation (see Section \ref{S1}) with a piecewise parabolic potential function, or, equivalently, a Fokker-Planck equation with a piecewise linear drift (see \eqref{diffusioneq} and \eqref{driftt} below). Such problems  arise in a variety of other applications, such as linear systems driven by white noise \cite{At67,AtC}, the Kramers' problem \cite{MP} and escape over potential barriers \cite{LRH}. Invariably, the solution involves the parabolic cylinder functions, and these we discuss in detail in Section \ref{S5}.

The key in determining the spectral gap is in fact determining the Laplace transform of the transient probability distribution over time. The spectral gap then follows from the dominant singularity of the Laplace transform. The main results are presented in Section \ref{sec:main}, the three special cases ($\eta=0,1,\infty$) are discussed in Section \ref{secpcases} and the proofs are given in Section \ref{sec:proofs}. Before the proofs we give some basic background on Schr\"{o}dinger equations (Section \ref{S1}) and parabolic cylinder functions (Section \ref{S5}), whose properties are heavily used later. In Section \ref{S8} we establish monotonicity properties of the spectral gap.

\section{Main results}\label{sec:main}

The diffusion process  $(X(t))_{t\geq 0}$ is a Markov process on the real line with continuous paths and density $p=p(x,t)=p(x,t;x_0;\beta,\eta)$ that satisfies the
forward Kolmogorov equation
\begin{equation}\label{diffusioneq}
\frac{\partial p}{\partial t}=-\frac{\partial}{\partial x}[a(x)p]+\frac{\partial^2p}{\partial x^2},
\end{equation}
where
\begin{equation}\label{driftt}
a(x)=\left\{
       \begin{array}{ll}
         -\beta-\eta x, & \hbox{$x\geq 0$}, \\
         -\beta-x, & \hbox{$x\leq 0$,}
       \end{array}
     \right.
\end{equation}
and (with $\delta(\cdot)$ the Dirac function and $p_x=\partial p/\partial x$)
\begin{align}
p(x,0)&=\delta(x-x_0),\label{444}\\
p(0^+,t)&=p(0^-,t), \quad p_x(0^+,t)=p_x(0^-,t),\label{555}
\end{align}
and $p(x,t)$ must decay as $x\rightarrow \pm\infty$.
The limiting distribution of the diffusion process is (see \cite{garnett})
\begin{equation}\label{steadyssty}
p(x,\infty;x_0;\beta,\eta)=C\left\{
       \begin{array}{ll}
          e^{-\frac{1}{2}\eta x^2}e^{-\beta x}, & \hbox{$x> 0$}, \\
          e^{-\frac{1}{2}     x^2}e^{-\beta x}, & \hbox{$x<0$,}
       \end{array}
     \right.
\end{equation}
where $C^{-1}=\int_0^\infty e^{-\frac{1}{2}\eta x^2}e^{-\beta x} dx+\int_{-\infty}^0e^{-\frac{1}{2}x^2}e^{-\beta x}dx$.

As shall be discussed in Section \ref{S1}, this problem has a purely discrete spectrum for all $\eta> 0$, and it is confined to the real axis. The spectral gap can thus be defined as the absolute value of the least negative eigenvalue of the operator in the right-hand side of \eqref{diffusioneq}. It governs the asymptotic rate of convergence to
the stationary distribution. An alternative description of the spectral gap is the absolute value of the singularity closest to the imaginary axis in the range ${\rm Re}(\theta)<0$ of the Laplace transform $\hat{p}$. Denote this dominant singularity by $\hat\theta$ and the spectral gap by $r(\beta,\eta)$. The relaxation time, which measures the time it takes for the system to approach its steady-state behavior, is defined as (see \cite{blancvandoorn,cohen})
\begin{equation}\label{deftau}
\tau=\inf\{T:p(x,t;x_0;\beta,\eta)-p(x,\infty;x_0;\beta,\eta)=O(e^{-t/T})\},
\end{equation}
and hence $\tau^{-1}=-{\rm Re}(\hat{\theta})=r(\beta,\eta)$. For this problem $\hat\theta$ is real, so that $-\hat\theta =r(\beta,\eta)$.
Our definition of the relaxation time in \eqref{deftau} assumes the initial condition $p(x,0)=\delta(x-x_0)$ in \eqref{444}, and then the approach to equilibrium is governed by $\lambda_1=r$. But we could certainly have initial conditions that would lead to a faster approach. For example, if $p(x,0)=p(x,\infty)$ then $p(x,t)=p(x,\infty)$ for all $t$ and equilibrium is attained instantaneously. We could also have initial distributions $p(x,0)$ that have zero projections on, say, the first $L$ eigenfunctions, and then the sums in \eqref{starr} and \eqref{starrstarr} below would be replaced by $\sum_{n=L+1}^\infty e^{-\lambda_n t}c_n\phi_n(x)$, where the $c_n$ may be computed in terms of $p(x,0)$. Then the approach to equilibrium would be governed by eigenvalue $\lambda_{L+1}$.


Here is the main result:

\begin{theorem}\label{thmspec}
The spectral gap of the diffusion process $(X(t))_{t\geq 0}$ is given by $r(\beta,\eta)=-\hat{\theta}$ where
$\hat\theta$ is the least negative solution to $\funcDD(\theta;\eta,\beta)=0
$
with
\begin{align}\label{mainD}
\funcDD(\theta;\eta,\beta)&=-\sqrt{\eta}D_{-\theta}(-\beta)D_{-\theta/\eta}'(\tfrac{\beta}{\sqrt{\eta}})-D_{-\theta}'(-\beta)D_{-\theta/\eta}(\tfrac{\beta}{\sqrt{\eta}})
\end{align}
$D_\nu(z)$ the parabolic cylinder function with index $\nu$ and argument $z$, and $D'_\nu(z)=\frac{d}{dz}D_\nu(z)$. If $\beta=0$, solving $\funcDD=0$ is equivalent to finding the roots of
\begin{equation}
\frac{\sqrt{\eta}}{\Gamma(\frac{\theta}{2\eta})\Gamma(\frac{1+\theta}{2})}+\frac{1}{\Gamma(\frac{\theta}{2})\Gamma(\frac12+\frac{\theta}{2\eta})}=0.
\end{equation}
\end{theorem}
We also note that the result in \eqref{mainD} corresponds to taking the limit of the discrete queueing model with the scaling in \eqref{eqq1} and then \eqref{deftau} examines what happens for large times.
We show in Appendix \ref{DISCRETE}  that \eqref{mainD} may also be obtained from the exact solution of the $M/M/s+M$ queue.
We have
\begin{proposition}
The spectral gap in the discrete $M/M/s+M$ model is the least negative solution to
\begin{equation}
\Delta(\theta)=F_s(\theta)H_{s-1}(\theta)-H_s(\theta)F_{s-1}(\theta)=0,
\end{equation}
where $F_n$ and $H_n$ are the contour integrals in \eqref{eqA8} and \eqref{eqA10}. For $s\to\infty$, with $\rho_0=1-\beta/\sqrt{s}$ and $\theta=O(1)$, the roots of $\Delta(\theta)$ may be approximated by those of $\funcDD$ in \eqref{mainD}.
This shows that the exchange of the limit in \eqref{eqq1} and of large time is permissible in this particular case.
\end{proposition}

\begin{figure}\centering
 \includegraphics[width= .7 \linewidth]{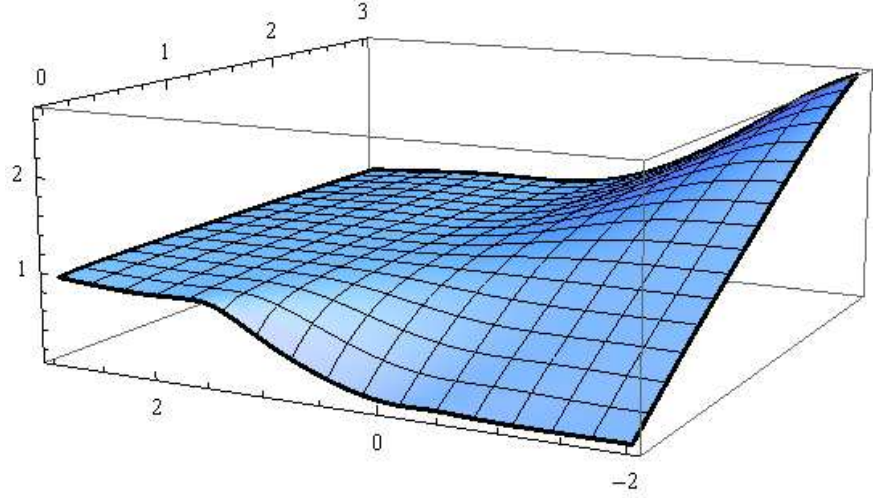}
  \caption{Spectral gap $r(\beta,\eta)$ for $\beta\in[-2,3]$ and $\eta\in(0,3]$.}
  \label{fig11}
\end{figure}

Theorem \ref{thmspec} is an implicit description of the spectral gap, and it can be used to calculate $r(\beta,\eta)$ numerically or asymptotically. For some values of $\beta$ and $\eta$ the spectral gap is shown in Figure \ref{fig11}.
We observe that the spectral gap decreases with $\beta$ if $\eta>1$ and increases with $\beta$ if $\eta<1$. If $\eta=1$
the spectral gap is $r(\beta,1)=1$ for all $\beta$, since then the problem reduces to a standard Ornstein-Uhlenbeck process (see Section \ref{ouf}).
This suggests that for systems with a large reneging rate, increasing the load (increasing $\rho_0$ and decreasing $\beta$) leads to shorter time scales for achieving equilibrium, while the opposite it true for small reneging rates.
We also see that $r$ increases as a function of $\eta$. Later, we establish the monotonicity of $r$ with $
\beta$, which is suggested by the numerical results in Figure \ref{fig11} (see Section \ref{S8}).

To further substantiate our findings, we accompany the observations from Figure \ref{fig11} by
results for the spectral gap in various asymptotic regimes. In order to do so we assume that $\eta\rightarrow 0$ (small abandonment rate). The asymptotics for $\eta\rightarrow\infty$ can be obtained from the following important symmetry result, which we establish in Section \ref{sec:proofs}.

\begin{proposition}\label{lemsym}For the density $p$ there holds the symmetry relation
\begin{equation}\label{sym9}
p(x,t;x_0;\beta,\eta)=\sqrt{\eta}\cdot p(-x\sqrt{\eta},t\eta;-x_0\sqrt{\eta};-\beta/\sqrt{\eta},1/\eta).
\end{equation}
For $\funcDD$ in \eqref{mainD},
\begin{equation}\label{sym10}
\funcDD(\theta;\eta,\beta)=-\sqrt{\eta}\cdot\funcDD(\theta/\eta;1/\eta;-\beta/\sqrt{\eta})
\end{equation}
and consequently,
\begin{equation}\label{sym11}
r(\beta,\eta)=\eta\cdot r(-\beta/\sqrt{\eta},1/\eta).
\end{equation}
\end{proposition}

Next we give five different asymptotic results for $r(\beta,\eta)$ as $\eta\rightarrow 0$. Then, from \eqref{sym11}
we can immediately get results for $\eta\rightarrow\infty$. We shall consider five ranges of $\beta$, with $\beta<0$, $\beta\approx 0$,
$0<\beta<\beta_*$, $\beta\approx \beta_*$ and $\beta>\beta_*$. Here, $\beta_*$ is the smallest positive solution to $D_{\beta^2/4}'(-\beta)=0$.
 We summarize in Table \ref{tabble} the five cases and where the asymptotic result may be found.

\begin{table}[h]
\begin{centering}
\begin{tabular}{c|c}\hline
    range of $\beta$  &  asymptotic result\\ \hline
	$\beta<0 $  & \eqref{eq:ertlem1} \\
    $\beta>\beta_*=1.85722\ldots$ & \eqref{13thirt}-\eqref{AAA}\\
    $0<\beta<\beta_*$ & \eqref{eq:ert}\\
    $\beta=\gamma\sqrt{\eta}=O(\sqrt{\eta})$ & \eqref{166} \\
    $\beta-\beta_*=\eta^{1/3}W=O(\eta^{1/3})$ & \eqref{1885gd}-\eqref{LLL}\\ \hline
\end{tabular}
\vskip.5cm
\caption{Five asymptotic regimes. \label{tabble}}
\end{centering}
\end{table}

\begin{proposition}\label{lem1}
For $\beta<0$ the spectral gap behaves asymptotically as $r(\beta,\eta)\sim\eta$ with the correction term
\begin{equation}\label{eq:ertlem1}
r(\beta,\eta)-\eta\sim-\tfrac{\beta\sqrt{\eta}}{\sqrt{2\pi}}e^{-\frac{\beta^2}{2\eta}}\left[1+\beta e^{\frac{\beta^2}{2}}\int_{-\infty}^\beta e^{-\frac{u^2}{2}}du\right],
\end{equation}
and hence $r-\eta$ is exponentially small as $\eta\rightarrow 0$.
\end{proposition}
Proposition \ref{lem1} describes the part at the far right end of Figure \ref{fig11}, where $r$ increases linearly with $\eta$. Since $\beta<0$, the diffusion process is mostly in the positive part of the state space, since the process has a positive drift for $0<x<-\beta/\eta$ and an equilibrium point at $x=-\beta/\eta=|\beta|/\eta\gg 1$. Hence, particularly when there is little reneging, one has to be far up in the state space before the process starts stabilizing. For the underlying queueing model, this scenario corresponds to large queues building up until enough customers renege so that the situation stabilizes.
For this scenario, the spectral gap $r=O(\eta)$ suggests large relaxation times. Note also from \eqref{steadyssty} that the steady-state distribution concentrates about $x=-\beta/\eta$. Table \ref{tableproposition1} compares exact and asymptotic results for $\beta=-1$.

\begin{table}[h!]
\begin{centering}
\begin{tabular}{c|cc}\hline
            & \multicolumn{2}{c}{$\beta=-1$}\\
    $\eta$  &   $r(\beta,\eta)-\eta$  & $\eqref{eq:ertlem1}$\\ \hline
   0.500&	   2.50092$\cdot 10^{-2}$&	   3.57325$\cdot 10^{-2}$\\
   0.400&	   1.91877$\cdot 10^{-2}$&	   2.48906$\cdot 10^{-2}$\\
   0.300&	   1.16366$\cdot 10^{-2}$&	   1.42105$\cdot 10^{-2}$\\
   0.200&	   4.29814$\cdot 10^{-3}$&	   5.04257$\cdot 10^{-3}$\\
   0.100&	   2.64792$\cdot 10^{-4}$&	   2.92685$\cdot 10^{-4}$\\
   0.050&	   1.32910$\cdot 10^{-6}$&	   1.39448$\cdot 10^{-6}$\\
   0.025&	   4.25017$\cdot 10^{-11}$&	   4.47665$\cdot 10^{-11}$\\ \hline
\end{tabular}
\vskip.5cm
\caption{Results for $\beta=-1$. \label{tableproposition1}}
\end{centering}
\end{table}


We next consider $\beta$ positive and sufficiently large, where we obtain a very different result for $r$.

\begin{proposition}\label{lem2}
For $\beta> \beta_*=1.85722\ldots$ and $\eta\rightarrow 0$
\begin{equation}\label{13thirt}
r(\beta,\eta)=r_0(\beta)+\mathcal{A}(\beta)\eta+O(\eta^2),
\end{equation}
where $r_0(\beta)$ is defined implicitly as the minimal positive solution to
\begin{equation}\label{defr0}
{D_{r_0}'(-\beta)}=\sqrt{\beta^2/4-r_0}{D_{r_0}(-\beta)}.
\end{equation}
The correction term is given by
\begin{equation}\label{AAA}
\mathcal{A}(\beta)=\frac{1}{2}\frac{\beta-\sqrt{\beta^2-4 r_0}}{\beta^2-4r_0}D_{r_0}(-\beta)\Big[\frac{\partial \tilde \funcDD}{\partial p}\Big|_{p=r_0(\beta)}\Big]^{-1},
\end{equation}
where $\tilde \funcDD (p,\beta)=D'_p(-\beta)-D_p(-\beta)\sqrt{\beta^2/4-p}$ {\rm(}so that $\tilde \funcDD (r_0(\beta),\beta)=0${\rm)}. We later show in Appendix \ref{appendixE} that $\mathcal{A}(\beta)>0$ so that $r(\beta,\eta)-r_0(\beta)$ is positive for sufficiently small $\eta$.
\end{proposition}

\begin{figure}[h]
\psfrag{a}{$p_1(\beta)$}
\psfrag{b}{$p_2(\beta)$}
\psfrag{c}{$p_3(\beta)$}
\psfrag{d}{$p_4(\beta)$}
\psfrag{e}{$p_5(\beta)$}
\psfrag{k}{$\beta$}
\psfrag{f}{$p=\beta^2/4$}
\psfrag{p}{$p$}
 \begin{center}
    \includegraphics[width = 0.5\linewidth]{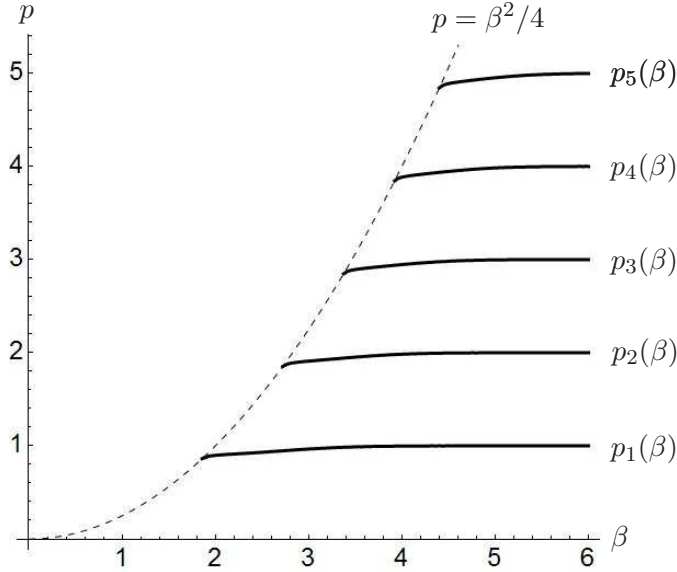}
 \end{center}
  \caption{Solutions to $\sqrt{\beta^2/4-p}D_p(-\beta)=D_p'(-\beta)$.}
 \label{fig2}
\end{figure}

The equation \eqref{defr0} corresponds to the discrete part of the spectrum of the Halfin-Whitt diffusion with no reneging (i.e., with $\eta=0$ in \eqref{driftt}); see the discussion after Proposition \ref{propCC}.
The various solution branches of \eqref{defr0} are demonstrated in Figure \ref{fig2}, where we plot the implicit function $\tilde \funcDD (p,\beta)=0$ for $\beta,p>0$.

Proposition \ref{lem2} applies to the flat part in Figure \ref{fig11}. Indeed, when $\beta$ is large enough, the spectral gap is hardly influenced by $\eta$. The diffusion process will spend most time below zero, near $x=-\beta$. A likely queueing scenario would be that queues hardly ever build up, which makes the impact of reneging customers negligible. As the spectral gap is $O(1)$, we expect relaxation times that are $O(1)$. While asymptotically $r(\beta,\eta)$ ranges from $r_0(\beta)$ to $1$, numerically this corresponds to the interval $(\beta^2_*/4,1)=(.86231,1)$, which is quite small, leading to the flatness of the surface in Figure  \ref{fig11} for $\beta>\beta_*$. Table \ref{tableproposition2} compares exact and asymptotic results for $\beta=2$.	
\begin{table}[h!]
\begin{centering}
\begin{tabular}{c|cc}\hline
            & \multicolumn{2}{c}{$\beta=2$}\\
    $\eta$  &   $r(\beta,\eta)$  & $r_0(\beta)$ \\\hline
   0.5000	&   0.98463  &  0.93229 \\
   0.2500	&   0.97072  &  0.93229\\
   0.1000	&   0.95576  &  0.93229\\
   0.0500	&   0.94741  &  0.93229\\
   0.0250	&   0.94150  &  0.93229\\
   0.0100	&   0.93671  &  0.93229\\
   0.0050	&   0.93470  &  0.93229\\
   0.0025	&   0.93356  &  0.93229\\
   0.0010	&   0.93282  &  0.93229\\ \hline
\end{tabular}
\vskip.5cm
\caption{Results for $\beta=2$. \label{tableproposition2}}
\end{centering}
\end{table}

We next consider $\beta>0$ but with $\beta<\beta_*$, in which case \eqref{defr0} has no positive solutions.
\begin{proposition}\label{lem3}
For $0<\beta<\beta_*=1.85722\ldots$ and $\eta\rightarrow 0$,
\begin{equation}\label{eq:ert}
r(\beta,\eta)=\frac{1}{4}\beta^2+\eta^{2/3}|a_0|\left(\frac{\beta}{2}\right)^{2/3}+\tfrac{1}{2}\eta\left(|a_0|+\beta
\frac{D_{\beta^2/4}(-\beta)}{D_{\beta^2/4}'(-\beta)}\right)+O(\eta^{4/3}),
\end{equation}
where $a_0=\max\{z:{\rm Ai}(z)=0\}=-2.33810\ldots$ is the least negative root of the Airy function.
\end{proposition}

From Proposition \ref{lem3}
we see that the asymptotic series now involves powers of $\eta^{1/3}$, which illustrates the lack of analyticity of $r(\beta, \eta)$ at $\eta=0$. Now $r(\beta,\eta)\sim \beta^2/4$ and for $\eta=0$ the spectral gap is in fact exactly $\frac{1}{4}\beta^2$ (see Section \ref{spec1}).
Table \ref{tableproposition3} compares exact and asymptotic results for $\beta=1$.

\begin{table}[h!]
\begin{centering}
\begin{tabular}{c|cc}\hline
            & \multicolumn{2}{c}{$\beta=1$}\\
    $\eta$  &   $r(\beta,\eta)$  & \eqref{eq:ert}\\\hline
   0.5000	&   0.87510  &  1.1778 \\
   0.2500	&   0.72686  &  0.83452\\
   0.1000	&   0.54242  &  0.56732\\
   0.0500	&   0.44074  &  0.44990\\
   0.0250	&   0.37193  &  0.37593\\
   0.0100	&   0.31673  &  0.31836\\
   0.0050	&   0.29217  &  0.29306\\
   0.0025	&   0.27664  &  0.27713\\
   0.0010	&   0.26450  &  0.26472\\ \hline
\end{tabular}
\vskip.5cm
\caption{Results for $\beta=1$\label{tableproposition3}.}
\end{centering}
\end{table}

When $\beta$ becomes small, both \eqref{eq:ertlem1} and \eqref{eq:ert}
become invalid, as the correction terms become larger than the leading term. Then a separate analysis leads to the following result.

\begin{figure}
\psfrag{a}{\footnotesize $\gamma$}
\psfrag{b}{\footnotesize $R$}
\begin{center}
 \includegraphics[width= .5 \linewidth]{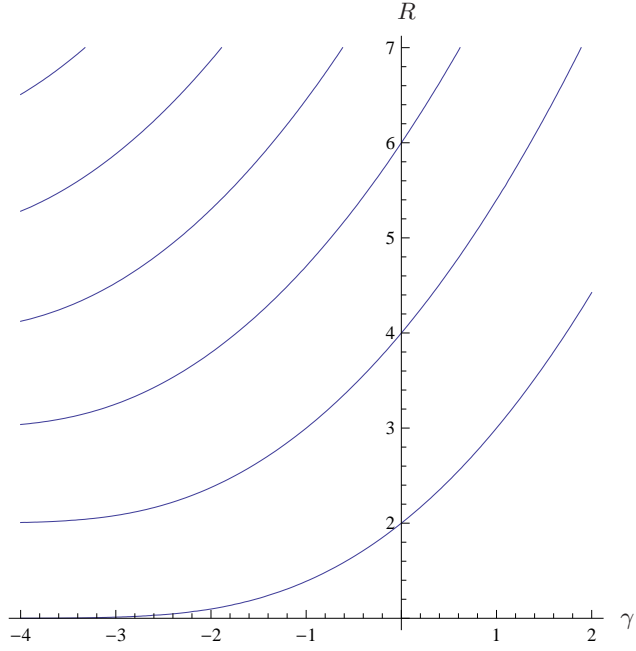}
 \end{center}
  \caption{A sketch of the solution branches of equation \eqref{eqrt} for $R>0$.}
  \label{fig3}
\end{figure}

\begin{proposition}\label{lem4}
Assume that $\beta$ is small, such that $\beta=\gamma\sqrt{\eta}=O(\sqrt{\eta})$ {\rm(}$\beta=0\leftrightarrow\gamma=0${\rm)}. Then
\begin{equation}\label{166}
r(\beta,\eta)\sim\eta R(\gamma),
\end{equation}
where $R$ is the minimal positive solution to
\begin{equation}\label{eqrt}
\gamma D_R(\gamma)=D_{1+R}(\gamma).
\end{equation}
\end{proposition}
Equation \eqref{eqrt} has infinitely many positive solutions, whose existence follows from
ODE theory, as discussed in \eqref{eqS17}-\eqref{eqS19}. Note that using the relations \eqref{eqP8} and \eqref{eqP9} below, \eqref{eqrt} is equivalent to $D_{R+1}(\gamma)+2D'_R(\gamma)=0$ or $D_R'(\gamma)+\frac{1}{2}\gamma D_R(\gamma)=0$.
Also, $R=0$ is a solution for any $\gamma$. In Figure \ref{fig3} we illustrate the solution branches of \eqref{eqrt} in the $(\gamma, R)$ plane, for $R>0$.

For $R\neq 0$, \eqref{eqrt} is also equivalent to $D_{R-1}(\gamma)=0$, and then the solutions
are precisely the positive eigenvalues of the reflected Ornstein-Uhlenbeck process (see the discussions in Section \ref{ouf} and surrounding \eqref{eqS17}-\eqref{eqS19} in Section \ref{S1}).
If $\gamma=0$ we can use \eqref{eqP7} to compute $D_{R-1}(0)$ and its roots are $R=2,4,6,\ldots$,
so that the spectral gap is $r(\beta,\eta)\sim 2\eta$ if $\beta=0$ or $\beta=o(\sqrt{\eta})$. For certain special values of $\gamma\neq 0$ we can also get some of the eigenvalues more explicitly. For example, we know from \eqref{eqP4} that $D_2(z)$ is proportional to the Hermite polynomial  ${\rm He}_2(z)=z^2-1$, so that  $D_2(1)=0$, and then $R(1)=3$ is the minimal root of \eqref{eqrt}. We also have $D_2(-1)=0$ so that when $\gamma=-1$, $R=3$ is a root of \eqref{eqrt}, but the minimal positive solution to \eqref{eqrt} is $R(-1)\approx 1.3882$ (which is illustrated in Figure \ref{fig3}).

In Appendix \ref{PROP6} we establish:
\begin{proposition}\label{propAA}
For $\gamma\to \pm \infty$, $R$ behaves as
\begin{equation}\label{p61}
R-1\sim-\frac{\gamma}{\sqrt{2\pi}} e^{-\gamma^2/2}, \quad \gamma\rightarrow-\infty
\end{equation}
and
\begin{equation}\label{p62}
R=\frac{\gamma^2}{4} +|a_0|\Big(\frac{\gamma}{2}\Big)^{2/3}[1+o(1)], \quad \gamma\rightarrow+\infty.
\end{equation}
\end{proposition}
Thus the approximation $r\sim \eta R$ on the $\gamma$-scale in Proposition
\ref{lem4} asymptotically matches to the results in Propositions \ref{lem1} and \ref{lem3}.
In fact, in view of \eqref{p62} the first two terms in \eqref{eq:ert} are a special case of \eqref{166}, but this is not true for the third ($O(\eta)$) term in \eqref{eq:ert}.
Note that $\eta \gamma^2/4=\beta^2/4$, $\eta \gamma^{2/3}=\eta^{2/3}\beta^{2/3}$, and if \eqref{eq:ertlem1} is expanded for small $\beta$, $r-\eta$ agrees with $\eta(R-1)$ for $\gamma\rightarrow-\infty$, in view of \eqref{p61}.

The results in \eqref{p61} and \eqref{p62} are also consistent with Figure \ref{fig11}, which suggests that $r(0,\eta)$ increases as a concave function of $\eta$.
The queueing counterpart is such that the load is one, and hence the reneging is necessary to alleviate the system. As $\eta$ becomes larger, more customers will leave the system, which reduces the queue lengths and, as seen from Proposition \ref{lem4}, shortens the relaxation times.

\begin{table}[h!]
\begin{centering}
\begin{tabular}{c|cc|cc|cc}\hline
            & \multicolumn{2}{c|}{$\gamma=1$}& \multicolumn{2}{c|}{$\gamma=0$}& \multicolumn{2}{c}{$\gamma=-1$}\\
    $\eta$  &   $r(\beta,\eta)$  & $\eta R(\gamma)$ &   $r(\beta,\eta)$  & $\eta R(\gamma)$&   $r(\beta,\eta)$  & $\eta R(\gamma)$\\ \hline
   0.50000	&   0.81266	 &  1.50000	&   0.65385	 &  1.00000	 &  0.54816	 &  0.69412\\
   0.25000	&   0.53164	 &  0.75000	&   0.38029	 &  0.50000	 &  0.29242	 &  0.34706\\
   0.10000	&   0.24948	 &  0.30000	&   0.16989	 &  0.20000	 &  0.12408	 &  0.13882\\
   0.05000	&   0.13266	 &  0.15000	&   0.08929	 &  0.10000	 &  0.06399	 &  0.06941\\
   0.02500	&   0.06896	 &  0.07500	&   0.04619	 &  0.05000	 &  0.03273	 &  0.03471\\
   0.01000	&   0.02848	 &  0.03000	&   0.01902	 &  0.02000	 &  0.01337	 &  0.01388\\
   0.00500	&   0.01446	 &  0.01500	&   0.00965	 &  0.01000	 &  0.00676	 &  0.00694\\
   0.00250	&   0.00731	 &  0.00750	&   0.00488	 &  0.00500	 &  0.00340	 &  0.00347\\
   0.00100	&   0.00295	 &  0.00300	&   0.00197	 &  0.00200	 &  0.00137	 &  0.00139\\ \hline
\end{tabular}
\vskip.5cm
\caption{Results for $\gamma=1,0,-1$. \label{tableproposition4}}
\end{centering}
\end{table}

Table \ref{tableproposition4} displays numerical results for $\gamma=1,0,-1$.

It remains to consider the case when $\beta\approx \beta_*$. Note that the correction term $\mathcal{A}(\beta)$ in \eqref{13thirt} develops a singularity as $\beta\downarrow\beta_*$, since $r_0(\beta)\rightarrow r_0(\beta_*)=\beta_*^2/4$. Also, the third ($O(\eta)$)  term in \eqref{eq:ert} becomes singular as
$\beta\uparrow\beta_*$, since by definition $\beta_*$ satisfies $D'_{\beta_*^2/4}(-\beta_*)=0$. Thus both \eqref{13thirt} and \eqref{eq:ert} cease to be valid near $\beta=\beta_*$ and we need a new expansion in this transition range.

\begin{proposition}\label{lem5}
For $\beta\approx \beta_*$ such that $\beta-\beta_*=\eta^{1/3}W$ with $W=O(1)$, and $\eta\rightarrow 0$,
\begin{align}\label{1885gd}
r(\beta,\eta)&=\frac{1}{4}\beta_*^2+\eta^{1/3}W\frac{\beta_*}{2}
+\eta^{2/3}\left(\frac{1}{4}W^2-\chi(W)\left(\frac{\beta_*}{2}\right)^{2/3}\right)+O(\eta),\nonumber\\
&=\frac{1}{4}\beta^2-\eta^{2/3}\Big(\frac{\beta_*}{2}\Big)^{2/3}\chi(W)+O(\eta),
\end{align}
where $\chi$ is the maximal solution to
\begin{equation}\label{199}
{\rm Ai}'(\chi)+\Big(\frac{2}{\beta_*}\Big)^{1/3}\cdot L\cdot W \cdot {\rm Ai}(\chi)=0
\end{equation}
with \begin{equation}\label{LLL}
L=\frac{1}{D_{\beta_*^2/4}(-\beta_*)}\left(\frac{d}{d\beta}[D_{\beta^2/4}'(-\beta)]\Big|_{\beta=\beta_*}\right)=2.73875\ldots.
\end{equation}
If $W=0$ {\rm($\beta=\beta_*$)}, then $\chi(0)=\max\{z:{\rm Ai}'(z)=0\}=-1.01870\ldots$, and as $W\rightarrow+\infty$, $\chi\rightarrow+\infty$.
\end{proposition}

\begin{figure}
\psfrag{a}{\footnotesize $w$}
\psfrag{b}{\footnotesize $\chi$}
\begin{center}
 \includegraphics[width= .5 \linewidth]{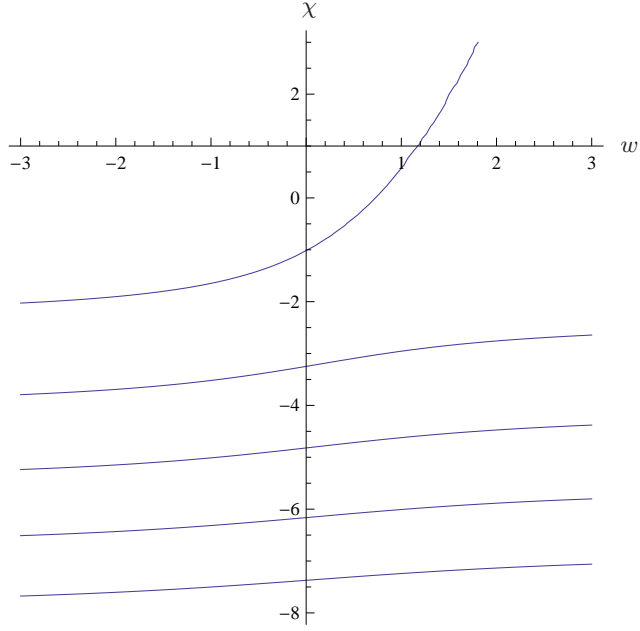}
 \end{center}
  \caption{A sketch of the solution branches of ${\rm Ai}'(\chi)+w {\rm Ai}(\chi)=0$.}
  \label{fig4}
\end{figure}

As discussed in \eqref{eqS20}-\eqref{eqS22} in Section \ref{S1}, the Sturm-Liouville ODE theory guarantees that there are infinitely many real solutions to \eqref{199}. The solution branches of \eqref{199} are illustrated in Figure \ref{fig4}.

We also note that if we order the roots of ${\rm Ai}(z)=0$ as $0>a_0>a_1>\cdots$ and the roots of ${\rm Ai}'(z)=0$ as $0>b_0>b_1>\cdots$, these roots interlace as $0>b_0>a_0>b_1>a_1>\cdots$, and this fact can be used to establish more directly that \eqref{199} has infinitely many solution branches, for any fixed $W$.

We have thus obtained the asymptotic connection between Propositions \ref{lem2} and \ref{lem3}.
Numerical results for the case $\beta=\beta_*$ are given in Table \ref{tableproposition5}.
\begin{table}[h!]
\begin{centering}
\begin{tabular}{c|ccc}\hline
            & \multicolumn{3}{c}{$\beta=\beta_*$}\\
    $\eta$  &   $r(\beta,\eta)$  & $\beta_*^2/4$ & \eqref{1885gd} \\ \hline
   0.5000	&   0.97803  &  0.86231 & 1.48841 \\
   0.2500	&   0.95673  &  0.86231 & 1.25673 \\
   0.1000	&   0.93129  &  0.86231 & 1.07644 \\
   0.0500	&   0.91493  &  0.86231 & 0.99721 \\
   0.0250	&   0.90139  &  0.86231 & 0.94729 \\
   0.0100	&   0.88770  &  0.86231 & 0.90845 \\
   0.0050	&   0.88016  &  0.86231 & 0.89138 \\
   0.0025	&   0.87462  &  0.86231 & 0.88062 \\
   0.0010	&   0.86966  &  0.86231 & 0.87225 \\ \hline
\end{tabular}
\vskip.5cm
\caption{Results for $\beta=\beta_*$. \label{tableproposition5}}
\end{centering}
\end{table}

This concludes our asymptotic analysis of the spectral gap in Theorem \ref{thmspec}. The proof of Theorem \ref{thmspec} follows immediately from an explicit expression for the Laplace transform $\hat{p}$ of the transient density, defined by
\begin{equation}
\hat{p}(x;\theta)=\int_{0}^\infty e^{-\theta t}p(x,t;x_0;\beta,\eta){\rm d}t, \quad {\rm Re}(\theta)>0.
\end{equation}
Defining the auxiliary function $\mathcal{M}$ by
\begin{align}\label{defMM}
\mathcal{M}(\theta;\eta,\beta)&=\sqrt{\eta}D_{-\theta}(\beta)D_{-\theta/\eta}'(\tfrac{\beta}{\sqrt{\eta}})-D_{-\theta}'(\beta)D_{-\theta/\eta}(\tfrac{\beta}{\sqrt{\eta}}),
\end{align}
we have the following result:
\begin{theorem}\label{greenpos}
Consider $x_0<0$, with $\funcDD$ in \eqref{mainD} and $\mathcal{M}$ in \eqref{defMM}, and assume that ${\rm Re}(\theta)>0$.
\begin{itemize}
\item[{\rm (i)}]For
 $x>0$,
\begin{align}\label{green1}
\hat{p}(x;\theta)&= e^{\frac{1}{2}\beta(x_0-x)}e^{\frac{1}{4}(x_0^2-\eta x^2)}\frac{D_{-\theta}(-x_0-\beta)D_{-\theta/\eta}(\frac{\eta x+\beta}{\sqrt{\eta}})}{\funcDD(\theta;\eta,\beta)}.
\end{align}
\item[{\rm (ii)}] For $x<x_0$,
\begin{align}\label{green2}
\hat{p}(x;\theta)=e^{\frac{1}{2}\beta(x_0-x)}&e^{\frac{1}{4}(x_0^2- x^2)}\frac{\Gamma(\theta)D_{-\theta}(-x-\beta)}{\sqrt{2\pi}}
 \nonumber\\
&\times\left(D_{-\theta}(x_0+\beta)+D_{-\theta}(-x_0-\beta)\frac{\mathcal{M}(\theta;\eta,\beta)}{\funcDD(\theta;\eta,\beta)}\right).
\end{align}
\item[{\rm (iii)}] For $x_0<x<0$,
\begin{align}\label{green3}
\hat{p}(x;\theta)=e^{\frac{1}{2}\beta(x_0-x)}&e^{\frac{1}{4}(x_0^2-x^2)}\frac{\Gamma(\theta)D_{-\theta}(-x_0-\beta)}{\sqrt{2\pi}}\nonumber\\
&\times\left(D_{-\theta}(x+\beta)+D_{-\theta}(-x-\beta)\frac{\mathcal{M}(\theta;\eta,\beta)}{\funcDD(\theta;\eta,\beta)}\right).
\end{align}
\end{itemize}
\end{theorem}
The proof of Theorem \ref{greenpos} is presented in Section \ref{prooflaplace}. Note that the results for $x_0>0$ follow
immediately from the symmetry relation \eqref{sym9}.

From \eqref{green1}-\eqref{green3}  we see that singularities of $\hat p$ may arise either due to those of $\Gamma(\theta)$ (which occur at $\theta=0,-1,-2,\ldots$) or from the zeros of $\funcDD(\theta;\eta,\beta)=0$. But in Appendix \ref{T8} we establish:
\begin{proposition}\label{propBB}
The singularities of \eqref{green1}-\eqref{green3} are precisely the solutions to $\funcDD(\theta;\eta,\beta)=0$.
\end{proposition}
Hence, the large-time behavior of the diffusion process is dominated by the least negative zero of $\funcDD$, which gives the result on the spectral gap in Theorem \ref{thmspec}.

From \eqref{green1}-\eqref{green3}, by evaluating the contour integral for the inversion of the Laplace transform $\hat p (x;\theta)$, we can obtain a spectral expansion of the form
\begin{equation}\label{starr}
p(x,t)=p(x,\infty)+ e^{\frac{1}{2}\beta(x_0-x)}e^{\frac{1}{4}(x_0^2-\eta x^2)}\sum_{n=1}^\infty e^{-\lambda_n t}\psi_n^{-}(x_0)\psi_n^{+}(x), \quad x>0,
\end{equation}
and
\begin{equation}\label{starrstarr}
p(x,t)=p(x,\infty)+ e^{\frac{1}{2}\beta(x_0-x)}e^{\frac{1}{4}(x_0^2- x^2)}\sum_{n=1}^\infty e^{-\lambda_n t}\psi_n^{-}(x_0)\psi_n^{-}(x), \quad x<0,
\end{equation}
where
\begin{equation}\label{}
\psi_n^{+}(x)=\sqrt{k_n}\frac{D_{\lambda_n/\eta}(\frac{\eta x+\beta}{\sqrt{\eta}})}{D_{\lambda_n/\eta}(\frac{\beta}{\sqrt{\eta}})},  \quad x>0,
\end{equation}
\begin{equation}\label{}
\psi_n^{-}(x)=\sqrt{k_n}\frac{D_{\lambda_n}(-x-\beta)}{D_{\lambda_n}(-\beta)},  \quad x<0,
\end{equation}
and
\begin{equation}\label{}
k_n=\frac{1}{\Delta_n^*}D_{\lambda_n}(-\beta)D_{\lambda_n/\eta}\Big(\frac{\beta}{\sqrt{\eta}}\Big), \quad \Delta_n^*=\frac{\partial \funcDD(\theta;\eta,\beta)}{\partial \theta}\Big|_{\theta=-\lambda_n}.
\end{equation}
The eigenfunctions $\psi_n^{\pm}$ then satisfy the orthogonality relation
\begin{equation}\label{OR}
\int_{-\infty}^0 \psi_n^-(x)\psi_m^-(x)dx+\int_{0}^\infty \psi_n^+(x)\psi_m^+(x)dx=\delta(n,m).
\end{equation}
When $\eta=1$ we have $\lambda_n=n$, $k_n=D_n^2(\beta)/(n!\sqrt{2\pi})$, and $D_n(\beta)=e^{-\beta^2/4}{\rm He}_n(\beta)$ so that \eqref{OR} reduces to \begin{equation}\label{}
\int_{-\infty}^\infty \frac{1}{n!\sqrt{2\pi}}{\rm He}_n(x+\beta){\rm He}_m(x+\beta)e^{-\frac12 (x+\beta)^2}dx=\delta(n,m).
\end{equation}
Note that the pole at $\theta=0$ of \eqref{green1}-\eqref{green3} corresponds to the steady state behavior $p(x,\infty)$, while the poles at $\theta=-\lambda_N$ and their residues lead to the decaying terms in \eqref{starr} and \eqref{starrstarr}. However, the spectral expansion does not yield any more insight than \eqref{green1}-\eqref{green3}.

\section{Three special cases}\label{secpcases}
We shall now consider the three special cases of the diffusion process that arise by setting $\eta$ equal to zero, one and infinity.

\subsection{The Halfin-Whitt diffusion}\label{spec1}
As $\eta\rightarrow 0$ we end up with a process that behaves like a Brownian motion with drift above zero and like an Ornstein-Uhlenbeck process below zero. In \cite{vlk} we have called this diffusion process the {\it Halfin-Whitt diffusion}, after
Halfin and Whitt \cite{halfinwhitt} who identified this process as a heavy-traffic limiting process for the $GI/M/s$ system.
The mean hitting time of the Halfin-Whitt diffusion was obtained in Maglaras and Zeevi \cite{maglaraszeevi}. Gamarnik and Goldberg \cite{goldberggamarnik} were the first to identify the spectral gap of the $M/M/s$ system, asymptotically in the Halfin-Whitt regime.
\begin{theorem}\label{theoremdave}{\rm (Gamarnik and Goldberg \cite{goldberggamarnik})}
Let $\beta_*=1.85722...$ represent the smallest positive solution to
$D_{\beta^2/4}'(-\beta)=0$. The spectral gap of the $M/M/s$ system, asymptotically in the Halfin-Whitt regime, is given by
\begin{equation}
r(\beta,0)=\left\{
               \begin{array}{ll}
                 \tfrac{1}{4}\beta^2, & \hbox{ $0<\beta\leq\beta_*$,} \\
                 r_0(\beta), & \hbox{ $\beta\geq\beta_*$.}
               \end{array}
             \right.
\end{equation}
\end{theorem}
To establish Theorem \ref{theoremdave}, Gamarnik and Goldberg used the framework of Karlin and McGregor \cite{km} for birth-death processes, and the result of Van Doorn \cite{vandoorn} on the spectral gap of the $M/M/s$ system. In \cite{goldberggamarnik} the starting point is the discrete $M/M/s$ model, and its spectral gap is then analyzed in the Halfin-Whitt regime \eqref{eqq1}. An alternative proof of Theorem \ref{theoremdave} was given by the authors in \cite{vlk} by deriving the expression for the Laplace transform $\hat{p}$ of the transient density in the diffusion limit, which shows that the limits of large time and \eqref{eqq1} may be, in this case, interchanged. Below we summarize the main result in  \cite{vlk}.
\begin{proposition}\label{propCC}
For $x_0<0$ and $x>0$ the Laplace transform of the transient density for the Halfin-Whitt diffusion with $\eta=0$ is
\begin{align}\label{277}
\hat{p}(x;\theta)&= e^{\frac{1}{4}x_0^2}e^{\frac{1}{2}\beta x_0}\frac{D_{-\theta}(-\beta-x_0)}{D_{-\theta}(-\beta)}
\frac{e^{-\frac{1}{2}x\beta-x\sqrt{\theta+\beta^2/4}}}{\sqrt{\theta+\beta^2/4}-D_{-\theta}'(-\beta)/D_{-\theta}(-\beta)}.
\end{align}
\end{proposition}
In Appendix \ref{3.1} we show that \eqref{277} indeed follows by letting $\eta\rightarrow 0^+$ in \eqref{green1}.
From \eqref{277} we see that there is a branch point at $\theta=- \beta^2/4$, and this will lead to a continuous spectrum in the range
${\rm Im}(\theta)=0$ and ${\rm Re}(\theta)\leq -\beta^2/4$. There is a pole at $\theta=0$ if $\beta>0$, while if $\beta<0$ the pole is absent. Note that  $D_0(-\beta)=\exp(-\beta^2/4)$. Other poles may occur at the roots of \eqref{defr0}, which we studied analytically and numerically in  \cite{vlk} (see also Figure \ref{fig2}).


\subsection{Free-space OU process}\label{ouf}
When $\eta=1$ it immediately follows from the process description that the diffusion process  $(X(t))_{t\geq 0}$ reduces to a free-space OU process,
for which it is known that (with $x_>=\max(x,x_0)$, $x_<=\min(x,x_0)$)
\begin{align}\label{lpou}
\hat{p}(x;\theta)&=
\frac{1}{\sqrt{2\pi}}e^{\frac{1}{4}(x_0^2-x^2)}e^{\frac{1}{2}\beta (x_0-x)}\Gamma(\theta)
D_{-\theta}(x_>+\beta)D_{-\theta}(-x_<-\beta).
\end{align}
Indeed, this result also follows from Theorem \ref{greenpos} using the Wronskian identity in \eqref{eqP9A}, which shows that
\begin{equation}
\funcDD(\theta;\eta,\beta)\rightarrow \frac{\sqrt{2\pi}}{\Gamma(\theta)}, \quad \eta\rightarrow 1.
\end{equation}
Expression \eqref{lpou} is obtained for example in \cite{T1962}, in the context of the harmonic oscillator (see \eqref{eqS1}-\eqref{eqS5} below).
Also, $\mathcal{M}\to 0$ (cf.~\eqref{defMM}) as $\eta\to 1$ and then \eqref{lpou} follows from \eqref{green1}-\eqref{green3}.
It is easy to invert the Laplace transform \eqref{lpou}, as its poles are at zero and at all negative integers. Hence, (see, e.g.,~\cite{keilson})
\begin{equation}\label{oudens2}
p(x,t;x_0;\beta,1)=\frac{e^{\frac{1}{4}(x_0^2-x^2)}e^{\frac{1}{2}\beta(x_0-x)}}{\sqrt{2\pi}}\sum_{n=0}^{\infty}\frac{D_n(x_0+\beta)D_n(x+\beta)}{n!}e^{-nt}, \quad x\in\mathbb{R}.
\end{equation}
Here $D_n(z)=e^{-z^2/4}{\rm He}_n(z)$ where ${\rm He}_n(\cdot)$ is the $n$th Hermite polynomial. Alternatively, there is the closed-form expression
\begin{equation}\label{oudens}
p(x,t;x_0;\beta,1)=\frac{1}{\sqrt{2\pi}}\sqrt{\frac{1}{1-e^{-2t}}}\exp\left[-\frac{(x+\beta-(x_0+\beta)e^{-t})^2}{2(1-e^{-2t})}\right], \quad x\in\mathbb{R}.
\end{equation}

\subsection{Reflected OU process}\label{subsec:rou}
As $\eta\rightarrow\infty$, the process will spend all its time below zero, and hence  $(X(t))_{t\geq 0}$ reduces to a reflected OU process (see Ward and Glynn \cite{wardglynn}, Linetsky \cite{linetsky} and Fricker et al.~\cite{robert}). In this limit we have $D_{-\theta/\eta}(\beta/\sqrt\eta)\rightarrow D_0(0)=1$ and $\sqrt{\eta}D_{-\theta/\eta}'(\beta/\sqrt\eta)\rightarrow -\beta/2$, and then, using \eqref{eqP8},
\begin{align}
\funcDD(\theta;\eta,\beta)&\rightarrow \theta D_{-1-\theta}(-\beta), \quad \eta\rightarrow \infty,\\
\mathcal{M}(\theta;\eta,\beta)&\rightarrow \theta D_{-1-\theta}(\beta), \quad \eta\rightarrow \infty,
\end{align}
which can be used to simplify \eqref{green2} and \eqref{green3}. Then the solution agrees with that in Xie and Knessl \cite[Eq. (4.3.2)]{XieKnessl1993}.

\section{Schr\"{o}dinger equations and spectral properties}\label{S1}

Here we give some basic background on
spectral properties that are relevant to
PDE's such as~(\ref{diffusioneq}). In particular we show
that the discreteness of the spectrum
for any $\eta>0$ follows from classic results
on the Schr\"{o}dinger equation.

We set $p=e^{-\lambda t}\phi(x)$ where $\lambda$ is a
spectral or eigenvalue parameter. Then~(\ref{diffusioneq}) and~(\ref{driftt}) lead to
\begin{align}
\phi''(x)+(\beta+x)\phi'(x)+(\lambda +1)\phi(x)&=0,\ x<0,\label{eqS1}\\
\phi''(x)+(\beta+\eta x)\phi'(x)+(\lambda +\eta)\phi(x)&=0,\ x>0,\label{eqS2}
\end{align}
and the interface conditions are $\phi(0^-)=\phi(0^+)$
and $\phi'(0^-)=\phi'(0^+)$. Furthermore, we can
transform~(\ref{eqS1}) and~(\ref{eqS2}) into the self-adjoint
form by setting
\begin{equation}
\label{eqS3}
\phi(x)=\begin{cases}
e^{-\beta x/2}e^{-x^2/4}\psi(x),&x<0,\\
e^{-\beta x/2}e^{-\eta x^2/4}\psi(x),&x>0,
\end{cases}
\end{equation}
which leads to the Schr\"{o}dinger equation
\begin{equation}\label{eqS4}
-\psi''(x)+ \V(x)\psi(x)=\E\psi(x),\ -\infty<x<\infty,
\end{equation}
where $\E$ and $\lambda$ are related by
\begin{equation}
\label{eqS5}
\E=\lambda +\dfrac{1}{2},
\end{equation}
and the ``potential'' function $\V(x)$ is
\begin{equation}
\label{eqS6}
\V(x)=\begin{cases}
\frac{1}{4}(\beta+x)^2,&x<0,\\
\frac{1}{4}(\beta+\eta x)^2+\frac{1-\eta}{2},&x>0.
\end{cases}
\end{equation}
We also require the eigenfunctions $\psi(x)$ to
satisfy $\psi(0^+)=\psi(0^-)$, and $\psi'(0^+)=\psi'(0^-)$.
Since the problem is defined over the
entire real line, additional conditions must be
imposed at $x=\pm\infty$, and most often
it is required that $\psi(x)\in L^2(-\infty,\infty)$, i.e.,
$\int^{\infty}_{-\infty}|\psi(x)|^2\, dx<\infty$. However, for parabolic
and piecewise parabolic potentials, such as
the one in~(\ref{eqS6}), this condition is equivalent
to simply rejecting solutions of~(\ref{eqS4}) that
have Gaussian growth as $x\to \pm \infty$.

We note that if the potential $\V(x)$
were exactly quadratic, say $\V(x)=x^2/4$, then
the problem is just the quantum harmonic
oscillator (or, for our application, the standard
Ornstein--Uhlenbeck process), and then the
eigenvalues are $\E_\N=\N+1/2$ and the corresponding
eigenfunctions are $\psi_\N(x)=c_\N e^{-x^2/4}\Hup e_\N(x)$, where
$c_\N$ is a normalizing constant and $\Hup e_\N(x)$ is the $\N^{\text{\upshape
th}}$ Hermite polynomial. Thus the spectrum
is purely discrete for quadratic
potentials.

We can also view the differential
equation in~(\ref{eqS4}) as constituting a singular Sturm--Liouville
boundary value problem. The study
of such problems dates back to the work of
Sturm in the nineteenth century, and they
are discussed in detail in the books of
Titchmarsh~\cite{T1962}, Stakgold~\cite{STAK}, Reid~\cite{REID} and
Coddington and Levinson \cite[Chapters 7--12]{CL}.
The problem in~(\ref{eqS4}) is singular since it is
defined over the infinite interval $x\in (-\infty,\infty)$.

Singular Sturm--Liouville problems are classified as either of
limit circle or limit point type. For
limit point problems the condition that
the solution be square
integrable is sufficient to determine it,
while limit circle problems require a more
explicit boundary condition at the singular
point(s) (which are at $x=\pm \infty$ for~(\ref{eqS4})).

Since Sturm--Liouville problems and
Schr\"{o}dinger equations are self-adjoint, their
spectra are confined to the real axis.
Singular problems may have both
discrete and continuous spectra. However,
there is a general result, originally due
to Weyl, with simplified proofs by
Titchmarsh appearing in~\cite{TMa,TMb}
(see also the book~\cite{T1962}), that guarantees
that~(\ref{eqS4}) will have a purely discrete
spectrum. This needs only the conditions that
\begin{enumerate}
\item[(i)]
$\V(x)$ be finite on finite intervals.
\item[(ii)]
$\V(x)\to +\infty$ as $x\to \pm\infty$.
\end{enumerate}
Our potential in~(\ref{eqS6}) clearly satisfies these
conditions and thus has a purely discrete
spectrum, for any $\eta > 0$. The fact
that $\V(x)$ has a jump discontinuity at $x=0$
does not affect the spectrum; it only
means that some jump conditions must be
specified at $x=0$. However, if $\eta=0$, then
the potential does not grow at $x=+\infty$,
and then in fact, as we discussed in~\cite{vlk},
the problem has a continuous spectrum
in the range $\lambda>\beta^2/4$ $(\E>\beta^2/4+1/2)$, and
may also have any number of discrete
eigenvalues, depending on the value of~$\beta$.
Much of the asymptotic work here assumes
that $\eta\to 0^+$, so we are looking at a
very singular limit where the discrete
spectrum begins to resemble a continuous
one, in certain ranges of~$\lambda$.

For the problem in~(\ref{eqS6}) the smallest
eigenvalue is $\E_0=1/2$ (thus $\lambda_0=0$) with
the corresponding eigenfunction being the
piecewise Gaussian
\[
\psi_0(x)=\begin{cases}
e^{-\beta x/2}e^{-x^2/4},&x<0\\
e^{-\beta x/2}e^{-\eta x^2/4},&x>0
\end{cases}
\]
and this corresponds to the steady state
distribution in our model.

Given the discrete spectrum we order
the eigenvalues $\E_\N$ as
\begin{equation}\label{eqS7}
\dfrac{1}{2}=\E_0<\E_1<\E_2<\dots<\E_\N<\dots
\end{equation}
with $\lambda_\N=\E_\N-1/2$. By general
results for Sturm--Liouville problems the
sequence $\{\E_\N\}$ satisfies $\E_\N\to \infty$ as
$\N\to \infty$. Also, for every eigenvalue there
is only one linearly independent eigenfunction,
so all eigenvalues are simple.
This can be shown directly
from~(\ref{eqS4}), for if $\psi(x)$ and $\widetilde{\psi}(x)$ corresponded
to the same eigenvalue~$\E$, then $\widetilde{\psi}(x)\psi''(x)
-\widetilde{\psi}''(x)\psi(x)=\frac{d}{dx}\left[\widetilde{\psi}(x)\psi'(x)-\widetilde{\psi}'(x)\psi(x)\right]=0$
so that the Wronskian $\widetilde{\psi}\psi'-\widetilde{\psi}'\psi$ is a
constant, say~$c$.
But by letting $x\to \pm \infty$ we conclude
that $c=0$ so that
\begin{equation}\label{eqS8}
\widetilde{\psi}^2(x)\dfrac{d}{dx}\left[\dfrac{\psi(x)}{\widetilde{\psi}(x)}\right]=0
\end{equation}
and thus $\widetilde{\psi}$ must be a multiple of~$\psi$. Since
$V(x)$ in~(\ref{eqS6}) is discontinuous at $x=0$, we can
use~(\ref{eqS8}) to conclude that $\psi(x)=k\widetilde{\psi}(x)$ for
$x>0$ and $\psi(x)=k'\widetilde{\psi}(x)$ for $x<0$, but
then continuity of $\psi,\widetilde{\psi}$ at $x=0$ forces $k=k'$.

There are two other singular Sturm--Liouville
problems that are relevant to
the analysis here. First consider
\begin{equation}\label{eqS17}
-\Psi''(X)+\left[\dfrac{1}{4}(X+\gamma)^2-\dfrac{1}{2}\right]\Psi(X)=\widetilde{\E}\Psi(X),\
0<X<\infty
\end{equation}
with the boundary condition
\begin{equation}\label{eqS18}
\Psi'(0)+\dfrac{\gamma}{2}\Psi(0)=0.
\end{equation}
This problem has a regular point at $X=0$,
where a standard boundary condition is
applied, and a singular point at $X=\infty$,
where we require that $\Psi(X)\in L^2(0,\infty)$. This
is a singular problem of limit point type
at $X=\infty$ which may be explicitly solved
in terms of parabolic cylinder functions,
with $\Psi(X)=D_{\widetilde{\E}}(X+\gamma)$. Then~(\ref{eqS18}) leads
to the eigenvalue condition
\begin{equation}\label{eqS19}
D_{\widetilde{\E}}'(\gamma)+\dfrac{\gamma}{2}D_{\widetilde{\E}}(\gamma)=0.
\end{equation}

The results of Weyl and Titchmarsh
again guarantee that the problem has a
purely discrete spectrum and thus an
infinite sequence of eigenvalues $\{\widetilde{\E}_{\N}\}$.
Also, $\widetilde{\E}_0=0$ is the lowest eigenvalue
with $\Psi_0(X)=e^{-(X+\gamma)^2/4}$. Note that~(\ref{eqS19})
is essentially the same as
equation~\eqref{eqrt} in Proposition~\ref{lem4}, since $D'_{\widetilde{\E}}(\gamma)
=\frac{1}{2}\gamma D_{\widetilde{\E}}(\gamma)-D_{\widetilde{\E}+1}(\gamma)$.
Thus the existence
of infinitely many real solutions to \eqref{eqrt} follows
from Sturm--Liouville ODE theory, though
in the next section we shall also show that
it follows from the oscillatory nature of
the parabolic cylinder functions, as functions
of their index~$\widetilde{\E}$.

Another singular Sturm-Liouville problem
is
\begin{equation}\label{eqS20}
-\psi''(x)+x\psi(x)=\E\psi(x),\ 0<x<\infty
\end{equation}
with the boundary conditions
\begin{equation}\label{eqS21}
\psi'(0)+\omega\, \psi(0)=0
\end{equation}
and $\psi(x)\in L^2(0,\infty)$. Here $\omega$ is a real
parameter. Again, since $V(x)=x$ grows linearly
as $x\to \infty$ and $x=0$ is
a regular point, we have a purely discrete
spectrum. But~(\ref{eqS20}) is related to the Airy
equation, with solutions proportional to
$\Ai(x-\E)$, and the eigenvalues are determined
by~(\ref{eqS21}), hence
\begin{equation}\label{eqS22}
\Ai'(-\E)+\omega\, \Ai(-\E)=0.
\end{equation}
But (\ref{eqS22}) is equivalent to equation~\eqref{199} in Proposition~\ref{lem5},
so again ODE theory can be used to establish
the existence of infinitely many solutions.
Note also that if $\omega=0$ the eigenvalues
are the roots of $\Ai'(\cdot)$, while if $\omega=\infty$
the eigenvalues are the roots of~$\Ai(\cdot)$.

To summarize we have given some
basic background on Sturm-Liouville theory,
and on Schr\"{o}dinger equations and their
spectral properties, that are useful in
the present investigations. In particular
this theory guarantees infinitely many
discrete solutions to the equations that
arise in Propositions~6 and~7, and in Theorem~\ref{thmspec}.

\section{Parabolic cylinder functions and their properties}\label{S5}

The parabolic cylinder equation is
the second order ordinary differential equation
\begin{equation}\label{eqP1}
y''(z)+\Big(p+\dfrac{1}{2}-\dfrac{z^2}{4}\Big)y(z)=0,
\end{equation}
where $z$ is a complex variable and $p$ is
a parameter. Since (\ref{eqP1}) has no singular
points (except at $z=\infty$) its solutions
are entire functions of~$z$ (see \cite{CL}).

One solution of (\ref{eqP1}) is denoted by
$D_p(z)$, which is called a parabolic cylinder
function of order~$p$, and it is defined
by the integral representation
\begin{equation}\label{eqP2}
D_{p}(z)=\dfrac{1}{i\sqrt{2\pi}}e^{z^2/4}\int_{\Br}u^p
e^{-zu}e^{u^2/2}\, du.
\end{equation}
Here $\Br$ is a vertical Bromwich contour
on which ${\rm Re}(u)>0$, and the branch
of~$u^p$ will be defined by $u^p=|u|^pe^{i\arg(u)}$
where $-\pi<\arg(u)\leq \pi$. Then the integrand
in the right-hand side of (\ref{eqP2}) is analytic exterior
to the branch cut where ${\rm Im}(u)=0$ and ${\rm Re}(u)\leq 0$.
The function $D_p(-z)$ provides a second
linearly independent solution to~(\ref{eqP1}), so that
\begin{equation}\label{eqP3}
c_1D_p(z)+c_2D_p(-z)
\end{equation}
is the general solution, with $c_1$ and $c_2$ being
complex constants. When $p=0,1,2,\dots$ is a
non-negative integer we can obtain $D_p(z)$ in
a closed form as
\begin{equation}\label{eqP4}
D_p(z)=D_\N(z)=e^{-z^2/4}\He_\N(z);\ N=0,1,2,\dots
\end{equation}
where $\He_\N(\cdot)$ is the $\N^{\text{\upshape th}}$ Hermite polynomial.
Note that if $p=\N$ the integrand in~(\ref{eqP2}) becomes
an entire function of~$u$. Here we use
the notation $\He(\cdot)$ for the Hermite polynomials,
so that $\He_0(z)=1$, $\He_1(z)=z$ and in
general $\He_\N(z)\sim z^\N$ as $z\to \infty$. When
$p=\N$, $D_\N(-z)=(-1)^\N D_\N(z)$ and another
linearly independent solution must be used
in~(\ref{eqP3}), but we shall not need it in the present
analysis. The function $D_p(z)$ is real valued
when $z$ and $p$ are real. As discussed in
\cite{TEMME}, $D_p(z)$ is an entire function
of both $p$ and $z$, and indeed one can
easily compute derivatives of all orders
from the integral representation in~(\ref{eqP2}). For example, we have
\begin{equation}\label{eqP5}
D_p'(z)=\dfrac{\partial}{\partial
z}D_p(z)=\dfrac{1}{i\sqrt{2\pi}}e^{z^2/4}\int_{\Br}\left(\dfrac{z}{2}-u\right)u^pe^{-zu}e^{u^2/2}\,
du
\end{equation}
and
\begin{equation}\label{eqP6}
\dfrac{\partial}{\partial
p}D_p(z)=\dfrac{1}{i\sqrt{2\pi}}e^{z^2/4}\int_{\Br}u^p(\log
u)e^{-zu}e^{u^2/2}\, du.
\end{equation}
In (\ref{eqP6}) $\log u$ is real for $u$ real and positive,
and analytic exterior to the cut $\text{\upshape Im}(u)=0$,
${\rm Re}(u)\leq 0$.

If $z=0$ the integrals in (\ref{eqP2}) and (\ref{eqP5})
may be expressed in terms of the Gamma
function, with
\begin{equation}\label{eqP7}
D_p(0)=\dfrac{\sqrt{\pi}2^{p/2}}{\Gamma\left(\tfrac{1-p}{2}\right)},\quad
D'_p(0)=\dfrac{-\sqrt{\pi}2^{(p+1)/2}}{\Gamma\left(\tfrac{-p}{2}\right)}.
\end{equation}
Since $\Gamma(z)$ has simple poles at $z=0,-1,-2,\dots$
it follows that $D_p(0)$ has simple zeros at
$p=1,3,5,\dots$, while $D'_p(0)$ has simple zeros
at $p=0,2,4,\dots$. This also shows that
the functions in (\ref{eqP7}), as functions of
$p$, tend to oscillate for $p>0$, but have
one sign for $p<0$ (actually, for all $p<1$ for
$D_p(0)$). Also, in view of the growth
of $\Gamma(z)$ as $z\to +\infty$, the functions in~(\ref{eqP7})
decay very rapidly for $p\to-\infty$.

Using (\ref{eqP2}) and (\ref{eqP5}) we can easily
derive the recurrence relations
\begin{align}\label{eqP8}
D'_p(z)+\dfrac{1}{2}zD_p(z)-pD_{p-1}(z)&=0,
\\
D'_p(z)-\dfrac{1}{2}zD_p(z)+D_{p+1}(z)&=0,
\label{eqP9}
\end{align}
which we shall use in the present analysis.

The Wronskian of $D_p(z)$ and $D_p(-z)$
is defined as $D_p(z)D'_p(-z)+D'_p(z)D_p(-z)$ and
it has a very simple form, with
\begin{equation}
\label{eqP9A}
D_p(z)D'_p(-z)+D'_p(z)D_p(-z)=-\dfrac{\sqrt{2\pi}}{\Gamma(-p)},
\end{equation}
which vanishes if $p=0,1,2,\dots$. We also note
that $D_p(z)$ and $D'_p(z)$ cannot both vanish
simultaneously. For if $D_p(z_*)=D'_p(z_*)=0$
for some $p$ and $z_*$, then (\ref{eqP1}) shows that
$D''_p(z_*)=0$. Then repeated differentiation of~(\ref{eqP1}) would
show that all derivatives of $D_p(z)$ vanish at $z=z_*$.
Then we could expand $D_p(z)$ in Taylor series about
$z=z_*$ to conclude that $D_p(z)=0$ in some
neighborhood of $z=z_*$. But since $D_p(z)$ is an entire
function of $z$ this would imply that $D_p(z)$
is identically zero, which is clearly not the case.

To better understand the behavior of
these special functions, many asymptotic
formulas have been derived for $z$ and/or
$p$ large. We summarize some of these below,
since they are used to establish our
main results. First, for $z$ large and
positive, we have (see for example \cite[p.~1093]{gradshteyn})
\begin{equation}\label{eqP10}
D_p(z)=z^pe^{-z^2/4}\left[1-\dfrac{p(p-1)}{2z^2}+\Oup(z^{-4})\right],\ z\to
\infty, \quad |\arg(z)|<3\pi/4.
\end{equation}
A more general result, which allows $z$ to
be negative, is the following (see for example \cite[p.~1094]{gradshteyn}):
\begin{align}\label{eqP11}
D_p(z)={}&e^{-z^2/4}z^p\left[1-\dfrac{p(p-1)}{2z^2}+\Oup(z^{-4})\right]\\
&-\dfrac{\sqrt{2\pi}}{\Gamma(-p)}e^{p\pi
i}e^{z^2/4}z^{-p-1}\left[1+\Oup(z^{-2})\right],\quad z\to \infty,\
\dfrac{\pi}{4}<\arg (z)<\dfrac{5\pi}{4}.\notag
\end{align}
Here we let $z=|z|e^{i\,\arg(z)}$ where $|\cdot|$ denotes
the complex modulus.

The limit of $z$ large and negative
corresponds to setting $\arg(z)=\pi$ in~(\ref{eqP11}),
and then the leading term becomes
\begin{equation}\label{eqP12}
D_p(z)\sim \dfrac{\sqrt{2\pi}}{\Gamma(-p)}e^{z^2/4}(-z)^{-p-1},\ z\to
-\infty,
\end{equation}
which holds as long as
we are away from the zeros of $1/\Gamma(-p)$,
which occur at $p=0,1,2,\dots$. If $p=\N=0,1,2,\dots$ we have $D_p(z)=e^{-z^2/4}$
$\He_N(z)$ and then~(\ref{eqP10}) holds for all values
of $\arg(z)$.

In some of the analysis that follows
we will need to consider cases where $p$
is not exactly a non-negative integer, but
is very close to one. For $p\neq \N$ the first
series in~(\ref{eqP11}) has Gaussian decay as $z\to-\infty$
$(\Oup(e^{-z^2/4}))$, while the second series in~(\ref{eqP11})
(cf.\ also~(\ref{eqP12})) has Gaussian growth as
$z\to-\infty$. But if $p$ is very close to an
integer these two terms may be of comparable
magnitude. For example if $p=\varepsilon$ is small
$\Gamma(-p)\sim -1/\varepsilon$ and if $\varepsilon\to 0$ and $z\to
-\infty$
in such a way that $e^{z^2/2}\varepsilon$ is $\Oup(1)$, then
the two parts of~(\ref{eqP11}) are roughly comparable.

The asymptotic results in (\ref{eqP10})--(\ref{eqP12})
follow easily by expanding~(\ref{eqP2}),
using techniques for the asymptotic
evaluation of integrals, such as the saddle
point method and singularity analysis.
General references for such methods are
the books of Bleistein and Handelsman
\cite{BH}, Wong \cite{WONG}, Szpankowski \cite{SPA},
and Flajolet and Sedgewick~\cite{FS}. The
integrand in~(\ref{eqP2}) has a saddle point at
$u=z$ and a branch point at $u=0$, and
one of these (or both) determines the asymptotic
behavior of~$D_p(z)$ as $z\to\infty$, for any
direction $\arg(z)$ in the complex plane.
The saddle leads to~(\ref{eqP10}) and the first
part of~(\ref{eqP11}), while the branch point leads
to~(\ref{eqP12}) and the second part of~(\ref{eqP11}).

We next consider a fixed (real) $z$ and
expand $D_p(z)$ in the limits of $p\to\pm\infty$.
Then
\begin{equation}\label{eqP13}
D_p(z)\sim\dfrac{1}{\sqrt{2}}\exp\left[\dfrac{p}{2}\log(-p)-\dfrac{p}{2}+z\sqrt{-p}\right], \quad
p\to-\infty,
\end{equation}
\begin{equation}\label{eqP14}
D_p(z)=\sqrt{2}\exp\left[\dfrac{p}{2}\log
p-\dfrac{p}{2}\right]\left[\cos\left(p\dfrac{\pi}{2}-\sqrt{p}z\right)+\Oup\left(p^{-1}\right)\right],\quad
p\to +\infty.
\end{equation}
From (\ref{eqP14}) we see faster than exponential
growth with~$p$, coupled with oscillations,
in view of the trigonometric factor. Near
a zero of the cosine the $\Oup(p^{-1})$ error
term may become important, and it may
be also explicitly obtained. Thus $D_p(z)$ has an infinite number of
zeros as $p$ increases toward $+\infty$, for any
fixed real~$z$. This is in sharp contrast
to viewing $D_p(z)$ for a fixed~$p$ as a function
of~$z$, in which case it has at most finitely many
zeros. It is
known \cite[p.~696]{a&s} that $D_p(z)$ has
no zeros in the range $p+1/2<z^2/4$.
Also, (\ref{eqP14}) shows that the large zeros can be estimated by
\begin{equation}\label{eqP15}
p=2\M+1+2\dfrac{z}{\pi}\sqrt{2\M+1}+O(1),\ \M\to \infty.
\end{equation}
The
results in~(\ref{eqP13}) and~(\ref{eqP14}) may be obtained, for
example, by expanding the integral in~(\ref{eqP2}). When $p\to-\infty$
the asymptotics are governed by a single
saddle point at $u=\sqrt{-p}$, while as $p\to +\infty$
two saddle points, at $u=\pm i\sqrt{p}$, contribute.

We next consider asymptotic limits
where $p$ and $z$ are simultaneously large,
restricting ourselves to real $p$ and~$z$.
The results quoted below are taken out of
Abramowitz and Stegun~\cite{a&s}, where we note
that in~\cite{a&s} the results are given for
the function $\U(a,z)$, which is related to $D_p(z)$
by $\U(a,z)=D_{-a-1/2}(z)$, or $D_p(z)=\U(-p-1/2,z)$.
A complete summary of the
asymptotics of~$\U$ is also given in Temme~\cite{TEMME}.

When $p+1/2>0$ and $z^2-4p$ is large and
positive, the so-called Darwin's expansions
apply, where \cite[pp.~689--690]{a&s}
\begin{align}\label{eqP16}
D_p(z)={}&\dfrac{\sqrt{\Gamma(p+1)}}{(2\pi)^{1/4}}\left(z^2-4p-2\right)^{-1/4}\\
&\times\exp\left[-\dfrac{1}{4}z\sqrt{z^2-4p-2}+\left(p+\dfrac{1}{2}\right)\log\left(\dfrac{z+\sqrt{z^2-4p-2}}{\sqrt{4p+2}}\right)\right]\notag\\
&\times\left[1+\Oup\left(\left(z^2-4p\right)^{-3/2}\right)\right],\notag
\end{align}
and
for $p+1/2<0$ and $z^2-4p$ large and positive
\begin{align}\label{eqP17}
D_p(z)={}&\dfrac{(2\pi)^{1/4}}{\sqrt{\Gamma(-p)}}\left(z^2-4p-2\right)^{-1/4}\\
&\times\exp\left[-\dfrac{1}{4}z\sqrt{z^2-4p-2}+\left(p+\dfrac{1}{2}\right)\log\left(\dfrac{z+\sqrt{z^2-4p-2}}{\sqrt{-4p-2}}\right)\right]\notag\\
&\times\left[1+\Oup\left(\left(z^2-4p\right)^{-3/2}\right)\right].\notag
\end{align}
The results in (\ref{eqP16}) and (\ref{eqP17}) apply for
$|z|\to \infty$ and $p\to \pm\infty$ with $z^2/(4|p|)$
fixed, and are uniform in the interval
$z^2/(4|p|)\in [1+\varepsilon,\infty)$ for any $\varepsilon>0$.
The expressions in (\ref{eqP16}) and~(\ref{eqP17}) are more uniform
than (\ref{eqP10}), (\ref{eqP12}) and~(\ref{eqP13}), and contain these as special cases.

Now consider $z$ and $p$ large with
$z^2-4p\to-\infty$. For $p>0$ the appropriate
expansion is now~\cite[p.~690]{a&s}
\begin{align}\label{eqP18}
D_p(z)\sim{}&\dfrac{2\sqrt{\Gamma(p+1)}}{(2\pi)^{1/4}}\left(4p+2-z^2\right)^{-1/4}\\
&\times\cos\left[\dfrac{\pi
p}{2}-\dfrac{z}{4}\sqrt{4p+2-z^2}-\left(p+\dfrac{1}{2}\right)\sin^{-1}\left(\dfrac{z}{\sqrt{4p+2}}\right)\right].\notag
\end{align}
Here $\sin^{-1}(\cdot)\in(-\pi/2,\pi/2)$.
The result in~(\ref{eqP18}) is uniform for $p\to\infty$
and $z\to\pm\infty$ for intervals where $z^2/(4p)
\in[0,1-\varepsilon]$ for any $\varepsilon>0$, except if
we are at or near a zero of the cosine
function, in which case the correction term(s)
to~(\ref{eqP18}) must be considered.

For a fixed large~$z$, $p$~has to
increase past approximately $z^2/4$ in order
for the zeros of~ $D_p(z)$ to become
evident. The expansions in~(\ref{eqP16})
and~(\ref{eqP18}) develop non-uniformities when
$z^2/(4p)\approx 1$ and there yet other
expansions apply. The following result
\cite[p.~689]{a&s} is more uniform and applies
for all $z^2/(4p)\in [0,\infty]$, as long as
$p$ and $|z|$ are large:
\begin{equation}\label{eqP19}
D_p(z)\sim2^{p/2}\Gamma\left(\dfrac{p+1}{2}\right)\left(4p+2\right)^{1/6}
\left(\dfrac{\tau}{\xi^2-1}\right)^{1/4}\Ai\left(\left(4p+2\right)^{2/3}\tau\right),
\end{equation}
\begin{equation}\label{eqP20}
\xi=\dfrac{z}{\sqrt{4p+2}},
\end{equation}
and
\begin{align}
\tau&=-\left(\dfrac{3}{8}\cos^{-1}(\xi)-\dfrac{3\xi}{8}\sqrt{1-\xi^2}\right)^{2/3},\ \xi\leq 1,
\label{eqP21}\\
\tau&=\left(\dfrac{3}{8}\xi\sqrt{\xi^2-1}-\dfrac{3}{8}\cosh^{-1}(\xi)\right)^{2/3},\
\xi\geq 1.
\label{eqP22}
\end{align}
Here $\Ai(\cdot)$ is the Airy function, which has the following
asymptotic behaviors as $z\to\pm\infty$ (see
\cite[p.~448]{a&s})
\begin{align}
\Ai(z)&\sim
\dfrac{1}{2\sqrt{\pi}}z^{-1/4}\exp\left(-\dfrac{2}{3}z^{3/2}\right),\ z\to
+\infty\label{eqP23}\\
\
\Ai(z)&\sim\dfrac{1}{\sqrt{\pi}}(-z)^{-1/4}\sin\left(\dfrac{2}{3}(-z)^{3/2}+\dfrac{\pi}{4}\right),\
z\to -\infty.\label{eqP24}
\end{align}
For $p\to\infty$ with a fixed $\tau>0$, we
can simplify (\ref{eqP19}) by using (\ref{eqP23}) to approximate
the Airy function, and then we obtain
(\ref{eqP16}) as a special case, when $\xi>1$.

Similarly, for $\xi<1$, we can use (\ref{eqP24}) to
approximate the Airy function in~(\ref{eqP19}) and
then (\ref{eqP19}) reduces to~(\ref{eqP18}), up to a
Stirling approximation for the Gamma factors.
When
$\xi=1$ $(\tau=0)$ or $\xi\approx 1$ we can also simplify
(\ref{eqP19}) considerably. Suppose that $p\to\infty$
and $\tau\to 0$ in such a way that $p^{2/3}\tau$ is
fixed. Setting
\begin{equation}\label{eqP27} 
p=\dfrac{z^2}{4}-\left(\dfrac{z}{2}\right)^{2/3}\delta
\end{equation}
with $\delta$ fixed as $z\to \infty$, we have
$
(4p+2)^{2/3}\tau\sim \delta
$
and then (\ref{eqP19}) simplifies to
\begin{equation}\label{eqP28} 
D_p(z)\sim
e^{-z^2/8}\left(\dfrac{z}{2}\right)^p\sqrt{2\pi}\left(\dfrac{z}{2}\right)^{1/3}\Ai(\delta),
\end{equation}
We shall use \eqref{eqP28} to establish Propositions \ref{lem3} and \ref{lem5}.

We have thus summarized the various
``uniform'' asymptotic approximations to
$D_p(z)$, where both $z$ and~$p$ become large.
Despite the seeming complexity of these
results, they are easily obtained from~(\ref{eqP2})
via the saddle point method. Indeed, setting
$u=zv$ (with $z>0$) in~(\ref{eqP2}) leads to
\begin{equation}\label{eqP29} 
D_p(z)=\dfrac{z^{p+1}e^{z^2/4}}{i\sqrt{2\pi}}\int_{\Br}e^{z^2\Phi(v;z,p)}dv, \quad \Phi(v)=\dfrac{p}{z^2}\log v+\dfrac{v^2}{2}-v.
\end{equation}
The saddle point equation is $\Phi'(v)=0$ which
is the quadratic equation
\begin{equation}
v^2-v+\dfrac{p}{z^2}=0 \quad \Rightarrow \quad v=v_{\pm}\equiv \frac{1}{2}\big[1\pm\sqrt{1-4pz^{-2}}\big].
\end{equation}
To expand the integral in (\ref{eqP29}) for $z\to\infty$,
$p\to\infty$ with $p/z^2$ fixed we find that for
$4p/z^2<1$ the real saddle at $v=v_+$
determines the asymptotic behavior, and we
ultimately obtain (\ref{eqP16}), up to a Stirling approximation
of $\Gamma(p+1)$ and the equivalence $z^2-4p-2
\sim z^2-4p$. In contrast, when $4p/z^2>1$ two
complex saddles, at $\frac{1}{2}\big[1\pm i\sqrt{4pz^{-2}-1}\big]$,
contribute to the asymptotics and then
we obtain~(\ref{eqP18}). The transition range in~(\ref{eqP28})
corresponds to $4pz^{-2}\approx 1$ and then the two
saddles coalesce to form a higher order
saddle, and such transitions invariably involve
Airy functions (see Chapter~9 in~\cite{BH}).

We have discussed here only approximations
to $D_p(z)$, but some of our main results
involve also the derivative $D'_p(z)$ (see, for example, Theorem~\ref{thmspec}). Its
asymptotics
follow from the integral in~(\ref{eqP5}), but the
same results can be obtained by formally
differentiating the results for $D_p(z)$, as in
this case term by term differentiation of
the asymptotic series is permissible. For
example, the logarithmic derivatives of
$D_p(z)$ and $\Ai(z)$ satisfy
\begin{equation}\label{eqP31} 
\begin{split}
\dfrac{D'_p(z)}{D_p(z)}&=-\dfrac{z}{2}+\dfrac{p}{z}+\dfrac{p(p-1)}{z^3}+\Oup\left(z^{-5}\right),\ z\to
+\infty,\\
\dfrac{\Ai'(z)}{\Ai(z)}&= -\sqrt{z}-\dfrac{1}{4z}+\Oup\left(z^{-5/2}\right),\
z\to +\infty,
\end{split}
\end{equation}
and later we shall make use of these results.

\section{Proofs}\label{sec:proofs}

\subsection{Proof of Proposition \ref{lemsym}}
Here we establish the symmetry relations \eqref{sym9}-\eqref{sym11}. These may be obtained without solving explicitly for $p(x,t)$. Consider the problem in \eqref{diffusioneq} with $p=p(x,t;x_0;\beta,\eta)$ and set
 \begin{equation}\label{4a1}
x=-\frac{x'}{\sqrt{\eta}},\quad t=\frac{t'}{\eta},\quad \beta=-\beta'\sqrt{\eta}, \quad x_0=-\frac{x_0'}{\sqrt{\eta}}
\end{equation}
with
\begin{equation}\label{4a2}
p(x,t;x_0;\beta,\eta)=R(x',t';x_0';\beta',\eta').
\end{equation}
Then \eqref{diffusioneq} becomes
\begin{align}\label{4a3}
\eta R_{t'}&=\eta R_{x'x'}+\eta\beta'R_{x'}+x'R_{x'}+R, \quad x'>0\nonumber,\\
\eta R_{t'}&=\eta R_{x'x'}+\eta\beta'R_{x'}+\eta[x'R_{x'}+R], \quad x'<0,
\end{align}
where $R_{x'}=\partial R/\partial x'$
and the initial condition becomes
\begin{align}\label{4a4}
R\big|_{t'=0}=\delta\Big(\frac{x_0'-x'}{\sqrt{\eta}}\Big)=\sqrt{\eta}\delta(x_0'-x')=\sqrt{\eta} \delta (x'-x_0'),
\end{align}
where
 we have used the scaling law of the delta function.

 Dividing \eqref{4a3} by $\eta$ and setting $R=\sqrt{\eta}\tilde R$ we see that $p$ and $\tilde R$ satisfy the same problem in the $x', t'$ variables, with $\eta$ replaced by $1/\eta$. Hence
 \begin{equation}\label{4a5}
\frac{1}{\sqrt{\eta}}R(x',t';x_0';\beta',\eta')=p(x',t';x_0';\beta,1/\eta')
\end{equation}
and then \eqref{4a2}, with \eqref{4a5} and \eqref{4a1}, is equivalent to \eqref{sym9}.

To establish \eqref{sym10} we again replace $\beta$ by $-\beta'\sqrt{\eta}$ and also set $\theta=\eta \theta'$. Then using \eqref{mainD} we obtain
\begin{align}\label{4a6}
\funcDD(\theta;\eta,\beta)&=-\sqrt{\eta}D_{-\eta\theta'}(\beta'\sqrt{\eta})
\Big[-\frac{d}{d\beta'}D_{-\theta'}(-\beta')\Big]+\frac{1}{\sqrt{\eta}}
\frac{d}{d\beta'}\Big[D_{-\eta\theta'}(\beta'\sqrt{\eta})\Big]D_{-\theta'}(-\beta')\nonumber\\
&=\sqrt{\eta}D_{-\eta\theta'}(\beta'\sqrt{\eta})D_{-\theta'}'(-\beta')+D_{-\eta\theta'}'(\beta'\sqrt{\eta})D_{-\theta'}(-\beta')\nonumber\\
&=-\sqrt{\eta}\funcDD(\theta';1/\eta,\beta')=-\sqrt{\eta}\funcDD(\theta/\eta;1/\eta,-\beta/\sqrt{\eta}),
\end{align}
which establishes \eqref{sym10}. Then since by definition $r=r(\beta, \eta)$ is the minimal negative root of $\funcDD=0$, the right-hand side of \eqref{4a6} has a root where $-\theta/\eta=r(-\beta/\sqrt{\eta},1/\eta)$ and then \eqref{sym11} follows immediately.

\subsection{Proof of Theorem \ref{greenpos}}\label{prooflaplace}

We let $\hat p (x;\theta)=\int_0^\infty e^{-\theta t}p(x,t)dt$ and note that $\hat p$ will be analytic in the right half-plane ${\rm Re}(\theta)>0$.
If $p$ satisfies (\ref{diffusioneq}) its Laplace transform satisfies
\begin{equation}\label{diffeqn}
\theta \hat{p}(x;\theta)-\delta(x-x_0)=-\frac{d}{d x}[a(x)\hat{p}(x;\theta)]+\frac{d^2\hat{p}(x;\theta)}{d x^2},
\end{equation}
where
\begin{equation}\label{3666}
-\frac{d}{d x}[a(x)\hat{p}(x;\theta)]=\left\{
                                        \begin{array}{ll}
                                          (x\eta+\beta) \frac{d}{d x}\hat{p}(x;\theta)+\eta\hat{p}(x;\theta), & \hbox{$x>0$,} \\
                                          (x+\beta) \frac{d}{d x}\hat{p}(x;\theta)+\hat{p}(x;\theta), & \hbox{$x<0$.}
                                        \end{array}
                                      \right.
\end{equation}
Assume that $x_0<0$ and $x>0$, so that $\delta(x-x_0)=0$.
By writing $\hat{p}=e^{-\eta x^2/4}e^{-\beta x/2}v(x;\theta)$, \eqref{diffeqn} reduces to the differential equation
\begin{equation}\label{pce}
v''(x;\theta)+\left(\tfrac{1}{2}\eta-\theta-\tfrac{1}{4}(x\eta+\beta)^2\right)v(x;\theta)=0, \quad x>0,
\end{equation}
where $v'=dv/dx$.
This is the parabolic cylinder equation (Erdelyi \cite{erdelyi}, p.~116)
and as we discussed in Section \ref{S5} two linearly independent solutions (at least for ${\rm Re}(\theta)>0$) are given by
\begin{equation}\label{3888}
D_{-\theta/\eta}(\tfrac{x\eta+\beta}{\sqrt{\eta}}), \quad D_{-\theta/\eta}(\tfrac{-x\eta-\beta}{\sqrt{\eta}}).
\end{equation}
Thus we write
\begin{equation}\label{rA}
\hat{p}(x;\theta)=\gamma_4(\theta)e^{-\eta x^2/4}e^{-\beta x/2}D_{-\theta/\eta}(\tfrac{x\eta+\beta}{\sqrt{\eta}}), \quad x>0,
\end{equation}
as the second solution in \eqref{3888} must be rejected due to its Gaussian growth (see \eqref{eqP12}).
For $x<0$ we must use the second expression in \eqref{3666}, and then solve \eqref{pce} with $\eta=1$. Now we reject solutions with Gaussian growth as $x\rightarrow-\infty$, so that for $x<x_0$ the appropriate solution to the second equation in \eqref{3666} (using \eqref{pce} with $\eta=1$) is
\begin{equation}\label{rB}
\hat{p}(x;\theta)=\gamma_1(\theta)e^{-\frac{1}{4}x^2}e^{-\frac{1}{2}\beta x}D_{-\theta}(-\beta-x), \quad x<x_0<0.
\end{equation}
But in the range $x_0<x<0$ the solution will involve both of the parabolic cylinder functions $D_{-\theta}(-\beta-x)$
and $D_{-\theta}(\beta+x)$, hence
\begin{equation}\label{rC}
\hat{p}(x;\theta)=e^{-\frac{1}{4}x^2}e^{-\frac{1}{2}\beta x}\left[\gamma_2(\theta)D_{-\theta}(-\beta-x)+\gamma_3(\theta)D_{-\theta}(\beta+x)\right].
\end{equation}
The functions $\gamma_j(\theta)$ are determined from continuity conditions at $x=0$ and $x=x_0$ (cf. \eqref{444} and \eqref{555}). Continuity of $\hat{p}$ and $\frac{d}{d x}\hat{p}$ at $x=0$ leads to
\begin{align}
\gamma_2D_{-\theta}(-\beta)+\gamma_3D_{-\theta}(\beta)&=\gamma_4 D_{-\theta/\eta}(\tfrac{\beta}{\sqrt{\eta}}),\label{rD}\\
-\gamma_2D_{-\theta}'(-\beta)+\gamma_3D_{-\theta}'(\beta)&=\gamma_4\sqrt{\eta} D_{-\theta/\eta}'(\tfrac{\beta}{\sqrt{\eta}}).\label{rE}
\end{align}
Continuity of $\hat{p}$ at $x=x_0$ yields
\begin{align}
\gamma_1D_{-\theta}(-\beta-x_0)=&\gamma_2D_{-\theta}(-\beta-x_0)+\gamma_3D_{-\theta}(\beta+x_0),\label{rF}
\end{align}
and the jump condition of $\frac{d}{d x}\hat{p}$ at $x=x_0$, i.e.,
\begin{equation}\label{eqq64}
\hat{p}(x_0^+;\theta)-\hat{p}(x_0^-;\theta)=-\int_{x_0^-}^{x_0^+}\delta(x-x_0){\rm d}x=-1,
\end{equation}
leads to
\begin{align}
-e^{\frac{1}{4}x_0^2}e^{\frac{1}{2}\beta x_0}= &-\gamma_2D_{-\theta}'(-\beta-x_0)+\gamma_3D_{-\theta}'(\beta+x_0)+\gamma_1D_{-\theta}'(-\beta-x_0)\Big. . \label{rG}
\end{align}
Here we used \eqref{rB} to compute $\hat{p}(x_0^-)$, \eqref{rC} to compute $\hat p (x_0^+)$, multiplied \eqref{eqq64} by $e^{\frac{1}{4}x_0^2}e^{\frac{1}{2}\beta x_0}$, and used also \eqref{rF}.
Equations \eqref{rD}-\eqref{rF} and \eqref{rG} give a $4\times 4$ linear system for the $\gamma_j$, whose solution leads to Theorem \ref{greenpos}.

We note that the calculations assumed that  ${\rm Re}(\theta)>0$. If ${\rm Re}(\theta)\leq 0$ the $4\times 4$ system may become singular, and in fact this occurs when $\theta=0$ and at the eigenvalues $-\theta=\lambda_N$, $N\geq 1$.
Theorem \ref{greenpos} thus gives the Laplace transform $\hat p$ for ${\rm Re}(\theta)>0$, and then the expression can be analytically continued to the left half-plane, since we know how to continue the parabolic cylinder functions, which are entire functions of $\theta$. After the continuation, locating the singularities in \eqref{green1}-\eqref{green3} in the range ${\rm Re}(\theta)\leq 0$ can be used, for example, to obtain the spectral representation in \eqref{starr}.

\subsection{General considerations for establishing Propositions \ref{lem1}-\ref{lem5}}\label{sec63}

Here we discuss some general principles about solving $\funcDD=0$ in Theorem
\ref{thmspec}, for the minimal root $r(\beta,\eta)$, in the limit of $\eta\rightarrow 0^+$. As discussed in Section \ref{S1}, general results for Schr\"{o}dinger equations and Sturm-Liouville problems show that the roots of $\funcDD=0$ are all on the real axis, and that the sequence of roots (or eigenvalues) $-\theta_N=\lambda_N$ satisfies $\lambda_N\rightarrow\infty$ as $N\rightarrow\infty$, for any fixed $\beta$ and $\eta>0$. Also, we know from Section \ref{S1} that the roots are all simple, and thus $\partial\funcDD/\partial \theta\neq 0$ when $\theta=-\lambda_N$.

Consider $r$ as a function of $\beta$ and $\eta$. Then $\funcDD(-r(\beta,\eta);\eta,\beta)=0$ and by implicit differentiation we obtain
\begin{equation}
\frac{\partial\funcDD}{\partial \theta} \Big|_{\theta=-r}\cdot \frac{\partial r}{\partial \beta}+\frac{\partial\funcDD}{\partial \beta}=0
\end{equation}
and since $\partial\funcDD/\partial \theta |_{\theta=-r}\neq 0$ we can use this relation to compute $\partial r/\partial \beta$ (for any $\beta$ and any $\eta>0$).
By taking higher order derivatives of $\funcDD=0$ with respect to $\beta$, a similar argument shows that $r$ has derivatives of all orders with respect to $\beta$. Also, by differentiating $\funcDD=0$ implicitly with respect to $\eta$, we conclude that $r$ has derivatives of all orders with respect to $\eta$, for any $\eta>0$. Thus $r(\beta,\eta)$ is infinitely smooth for all real $\beta$ and for $\eta>0$, and this is true for the higher roots also. Note that since $D_p(z)$ is an entire function of both $p$ and $z$ (see \cite{TEMME}), $\funcDD$ is an entire function of $\theta$ and $\beta$, and real analytic for $\eta>0$. However, the limit $\eta\rightarrow 0^+$ is quite singular, as we shall show.

The above discussion shows that the roots of $\funcDD=0$ vary smoothly with $\beta$ and $\eta$, and a root cannot simply appear/disappear, say at some critical value $\eta_c$. Thus for $\eta\rightarrow 0^+$ the roots have to lie in some range(s) of $\theta$. In Section \ref{S5} we gave detailed results of the different asymptotic expansions of the parabolic cylinder functions $D_p(z)$, for different ranges of $p,z$. Applying these results to the equation $\funcDD=0$, the function $\funcDD$ can be approximated by simpler functions in the limit of $\eta\rightarrow 0^+$, but these approximations are different in different ranges.

Different expansions arise for the ranges $-\theta=O(\eta)$, $0<-\theta<\beta^2/4$, $-\theta\approx \beta^2/4$ and $-\theta>\beta^2/4$. Depending on the value of $\beta$, we shall need to consider different ranges of $-\theta$ in order to locate the minimal root.


 In what follows we shall use the following principle: suppose the equation $F(u,\epsilon)=0$
 has roots $u_j=u_j(\epsilon)$ which depend on the small parameter $\epsilon$, and these roots are smooth functions of $\epsilon$. Also, suppose that $F(u,\epsilon)$ is an analytic function of both $u$ and $\epsilon$, with an expansion of the form
 \begin{align}\label{431}
F(u,\epsilon)=F_0(u)+\epsilon F_1(u)+O(\epsilon^2).
\end{align}
Then if $F_0(u)$ has a simple root at $u_*$ then $F(u,\epsilon)$ has a root close to $u_*$ for $\epsilon\rightarrow 0$. The same conclusion holds if $F$ is not analytic in $\epsilon$, but has an asymptotic expansion of the form \eqref{431}, where the expansion holds uniformly on some (finite) $u$ interval that contains $u_*$. In our case $-\theta$ plays the role of $u$ and $\epsilon$ will correspond to $\eta$ or a fractional power of $\eta$, such as $\sqrt{\eta}$ or $\eta^{1/3}$.

\subsection{Proof of Proposition \ref{lem4}}
For $\beta=O(\sqrt{\eta})$ and $\eta\rightarrow 0^+$, we shall show that $\funcDD=0$ has solutions in the range $-\theta=O(\eta)$.
We use the facts that $D_0(0)=1$ and $D_0'(x)\sim-x/2$ as $x\rightarrow 0$. We let $\eta\rightarrow\infty$ and from \eqref{mainD} we
obtain
\begin{align}\label{4777}
\funcDD(\theta;\eta,\beta)\rightarrow \frac{\beta}{2} D_{-\theta}(-\beta)-D_{-\theta}'(-\beta), \quad \eta\rightarrow\infty.
\end{align}
The error term in \eqref{4777} is uniformly $O(\eta^{-1})$ on finite $\theta$ intervals, and we note that the right-hand side of \eqref{4777} is an entire function of $\theta$, as will be the error terms.

Then the symmetry relation for $\funcDD$ in Proposition \ref{lemsym}
implies that if we scale $\theta=\eta S$ and $\beta=\gamma\sqrt{\eta}$ we obtain
\begin{align}\label{}
\funcDD(\theta;\eta,\beta)\sim -\sqrt{\eta}\left[ D_{-S}(\gamma)\frac{\gamma}{2}+D_{-S}'(\gamma)\right], \quad \eta\rightarrow 0,
\end{align}
with an error term that is uniformly $O(\eta)$ on finite intervals of $S$ and $\gamma$.
Proposition \ref{lem4} follows upon setting $S=-R$ and using the identity $D_{1-S}(\gamma)+D_{-S}'(\gamma)=\frac{1}{2}\gamma D_{-S}(\gamma)$
(see \cite{gradshteyn}, p.~1066).

Finally, we show that for $\beta=O(\sqrt{\eta})$ and $\eta\to 0^+$, there can be no roots of
$\funcDD=0$ in \eqref{mainD}, other than $\theta=0$, in the range $\theta=o(\eta)$. Consider the scales $\theta=\Theta(\eta^M)$ for $M>1$, setting $\theta=\theta_*\eta^M$. Then we have
\begin{align}\label{65-1}
\funcDD&=-\sqrt{\eta}D_{-\theta_*\eta^M}(-\sqrt{\eta}\gamma)D_{-\theta_*\eta^{M-1}}'(\gamma)-D_{-\theta_*\eta^M}'(-\sqrt{\eta}\gamma)D_{-\theta_*\eta^{M-1}}(\gamma).
\end{align}
Then both $\theta_*\eta^{M}$ and $\theta_*\eta^{M-1}$ are small and expanding the parabolic cylinder functions in \eqref{65-1} in Taylor series in both index and argument yields
\begin{align}\label{65-2}
\funcDD=&-\sqrt{\eta}\Big[D_0(-\sqrt{\eta}\gamma)+O(\theta_*\eta^{M})\Big]D_{-\theta_*\eta^{M-1}}'(\gamma)\nonumber\\
&-\Big[D_{0}'(-\sqrt{\eta}\gamma)+O(\theta_*\eta^{M})\Big]D_{-\theta_*\eta^{M}}(\gamma)\nonumber\\
&=\theta_*\eta^{M-1/2}\frac{d}{dp}\left[D_p'(\gamma)+\tfrac\gamma 2 D_p(\gamma)\right]\Big|_{p=0}[1+o(1)]\nonumber\\
&=\theta_*\eta^{M-1/2}D_{-1}(\gamma)[1+o(1)],
\end{align}
where the error term is $o(1)$ for $\eta\to 0^+$ and this is uniform in finite $\theta_*$ intervals. But then we conclude that $\theta_*=0$, contradicting our assumption that there is a root in the range $\theta=\Theta(\eta^M)$ for $M>1$. To obtain the last expression in \eqref{65-2} we also used the recurrence \eqref{eqP8}.

The calculation that led to \eqref{65-2} only used the fact that $\theta$ and $\theta/\eta$ are both small. Indeed, for any $\theta=o(\eta)$ we obtain \eqref{65-2}  with $\theta_*\eta^M$ replaced by $\theta$. Thus if there is a root in any range where $\theta=o(\eta)$, we again conclude that $\theta=0$. Hence, there can be no roots in ranges where $\theta=\Theta(\frac{\eta}{\log(1/\eta)})$, $\theta=\Theta(\frac{\eta}{\log\log(1/\eta)})$, etc.

\subsection{Proof of Proposition \ref{lem2}}
Here we establish Proposition \ref{lem2}. We represent the parabolic cylinder function (see also \eqref{eqP2} and \eqref{eqP29}) by the contour integral
\begin{align}
D_{-\theta/\eta}\Big(\frac{\beta}{\sqrt{\eta}}\Big)&=\frac{e^{\beta^2/4\eta}}{i\sqrt{2\pi}}\int_\mathcal{C} t^{-\theta/\eta}e^{-\frac{\beta}{\sqrt{\eta}}t}e^{t^2/2}dt\nonumber\\
&= \frac{e^{\beta^2/4\eta}}{i\sqrt{2\pi}}\eta^{\frac{\theta}{2\eta}-\frac{1}{2}}\int_\mathcal{C} \exp\left[
\tfrac{1}{\eta}\phi(z;\beta,\theta)\right]dz, \label{saddle123}
\end{align}
where $\phi=\frac{1}{2}z^2-\beta z-\theta\log z$ and $\mathcal{C}$ is a vertical Bromwich contour with ${\rm Re}(t)$ (or ${\rm Re}(z)$) larger than zero.
Now assume that $\theta+\beta^2/4>0$ and $\eta\rightarrow 0^+$. From the discussion in Section \ref{S5} we can use \eqref{eqP17} (or the more uniform result in \eqref{eqP19}) to approximate $D_{-\theta/\eta}(\beta/\sqrt{\eta})$ in this range. We briefly derive the asymptotic formula below.

The integrand  in \eqref{saddle123} has a saddle point where $\phi'(z)=0$, which occurs at
\begin{equation}
z_*=\frac{1}{2}\left[\beta+\sqrt{\beta^2+4\theta}\right].
\end{equation}
For $\beta^2+4\theta>0$ this saddle point lies on the real axis and the directions of steepest descent are ${\rm arg}(z-z_*)=\pm \pi/2$. Then expanding $\phi$ in Taylor series about $z=z_*$ and noting that
\begin{align}
\int_\mathcal{C} e^{\phi(z_*;\beta,\theta)/\eta} \exp\left[
\frac{1}{2\eta}\phi''(z_*;\beta,\theta)(z-z_*)^2\right]dz
=\exp\left[\frac{1}{\eta}\left(\frac{1}{2}z_*^2-\beta z_* -\theta \log z_*\right)\right]\frac{\sqrt{2\pi\eta}\ i}{\sqrt{\phi''(z_*;\beta,\theta)}},
\end{align}
we obtain
\begin{align}\label{5111}
D_{-\theta/\eta}\Big(\frac{\beta}{\sqrt{\eta}}\Big)&\sim e^{\beta^2/4\eta}\eta^{\frac{\theta}{2\eta}}
\Big(1+\frac{\theta}{z_*^2}\Big)^{-1/2}\exp\left[\tfrac{1}{\eta}(\tfrac{1}{2}z_*^2-\beta z_*-\theta\log z_*)\right].
\end{align}
In view of \eqref{eqP5}, the integral representation of $D_{-\theta/\eta}'(\beta/\sqrt{\eta})$ corresponds to multiplying the integrand in \eqref{saddle123} by $(\beta/2-z)/\sqrt{\eta}$ (after scaling $t=z/\sqrt{\eta}$).
Then again applying the saddle point method we see that the leading term in the expansion of $\sqrt{\eta}D_{-\theta/\eta}'(\beta/\sqrt{\eta})$ is essentially the same as \eqref{5111}, with the additional factor $\beta/2-z_*$, and thus
\begin{align}\label{5222}
\frac{\sqrt{\eta}D_{-\theta/\eta}'(\beta/\sqrt{\eta})}{D_{-\theta/\eta}(\beta/\sqrt{\eta})}
\sim \frac{\beta}{2}-z_*=-\frac{1}{2}\sqrt{\beta^2+4\theta}.
\end{align}
In \eqref{5222} we divided by $D_{-\theta/\eta}(\beta/\sqrt{\eta})$, since this is  nonzero in the range $\theta+\beta^2/4>0$, as discussed in Section \ref{S5}. Using \eqref{5222} we see that  $\funcDD=0$ leads to (after dividing by $D_{-\theta/\eta}(\beta/\sqrt{\eta})$)
\begin{align}\label{4e1}
D_{-\theta}(-\beta)\sqrt{\beta^2/4+\theta}-D_{-\theta}'(-\beta)+O(\eta)=0.
\end{align}
For $\eta\to 0$ we obtain the limiting equation in \eqref{defr0}. Up to now the calculation did not distinguish between $\beta<\beta_*$ or $\beta>\beta_*$, but in the former case \eqref{defr0} (or \eqref{4e1} with $\eta=0$) has no roots (other than $\theta=0$ if $\beta>0$).
By computing the first correction term in \eqref{5111}, which is $O(\eta)$ relative to the leading term, we can refine \eqref{5222} to
\begin{align}\label{4e2}
\frac{\sqrt{\eta}D_{-\theta/\eta}'(\beta/\sqrt{\eta})}{D_{-\theta/\eta}(\beta/\sqrt{\eta})}
=-\sqrt{\beta^2/4+\theta}-\frac{\eta}{\beta^2+4\theta}\Big(\frac{\beta}{2}-\sqrt{\beta^2/4+\theta}\Big)+O(\eta^2).
\end{align}
Here the error is uniform for $\eta\to 0$ with $\sqrt{\beta^2/4+\theta}>\epsilon>0$. However, the asymptotics break down as $\beta^2/4+\theta\to 0$, and a separate analysis is needed for locating the roots of $\funcDD=0$ in the range $-\theta\approx \beta^2/4$, which we carry out in the proofs of Propositions \ref{lem3} and \ref{lem5}.
With \eqref{4e2} we can improve \eqref{4e1} to
\begin{align}\label{4e3}
D_{-\theta}(-\beta)\sqrt{\beta^2/4+\theta}-D_{-\theta}'(-\beta)+
\frac{D_{-\theta}(-\beta)\eta}{\beta^2+4\theta}\Big(\frac{\beta}{2}-\sqrt{\beta^2/4+\theta}\Big)+O(\eta^2)=0.
\end{align}
Then the leading term in \eqref{13thirt} follows by letting $\eta\to 0$ in \eqref{4e3}, and the correction term follows by dividing \eqref{4e3} by $\eta$ and then letting $\eta\to 0$ and $-\theta\to r_0$ simultaneously, noting also that
\begin{align}\label{4e4}
\lim_{\eta\to 0, -\theta\to r_0 }& \frac{1}{\eta}\left[D_{-\theta}(-\beta)\sqrt{\beta^2/4+\theta}-D_{-\theta}'(-\beta)\right]\nonumber\\
&=\lim_{\eta\to 0, -\theta\to r_0 }\left[-\frac{1}{\eta}\tilde\funcDD(-\theta,\beta)\right]=
\frac{\partial \tilde\funcDD}{\partial \theta}\Big|_{-\theta=r_0}\cdot \lim_{\eta\to 0}\left[\frac{r(\beta,\eta)-r_0(\beta)}{\eta}\right],
\end{align}
where by definition in \eqref{13thirt}, $\mathcal{A}(\beta)=\lim_{\eta\to 0}[r(\beta,\eta)-r_0(\beta)]/\eta$.

\subsection{Proof of Proposition \ref{lem1}}

For $\beta<0$  we analyze the range $\theta=O(\eta)$, and we shall see that there are roots in this range. Note that for $\beta>0$ the calculation in the previous subsection (since \eqref{5222} applies for $\theta=O(\eta)$) shows that the only root in this range is at $\theta=0$. Since we are examining ranges where $\theta$ is small, we again use the fact that $D_{-\theta}(-\beta)$ is an entire function of $\theta$, and hence by Taylor series, noting that $D_0(-\beta)=e^{-\beta^2/4}$ and $D_0'(-\beta)=\frac12 \beta e^{-\beta^2/4}$, we have
\begin{align}\label{e52}
\frac{1}{\sqrt{\eta}}
\frac{\frac{d}{d\beta}D_{-\theta}(-\beta)}{D_{-\theta}(-\beta)}=\frac{1}{\sqrt{\eta}}\left[-\frac{\beta}{2}-\theta\mathcal{R}(\beta)+O(\theta^2)\right],
\end{align}
where (with $z=-\beta$)
\begin{align}
\mathcal{R}(\beta)=-\frac{d}{dp}\left.\frac{D_p'(z)}{D_p(z)}\right|_{p=0}
&=-\frac{e^{z^2/2}}{i\sqrt{2\pi}}\int_\mathcal{C}\frac{1}{t}e^{t^2/2-zt}dt\nonumber\\
&=-e^{\beta^2/2}\int_{-\infty}^\beta e^{-u^2/2}du. \label{eq133}
\end{align}
We then rewrite $\funcDD=0$ as
\begin{align}\label{4g1}
\frac{1}{\sqrt{\eta}}
\frac{\frac{d}{d\beta}D_{-\theta}(-\beta)}{D_{-\theta}(-\beta)}D_{-\theta/\eta}(\beta/\sqrt{\eta})&=D_{-\theta/\eta}'(\beta/\sqrt{\eta})\nonumber\\
&=\frac{1}{\sqrt{\eta}}\left[-\frac{\beta}{2}-\theta\mathcal{R}(\beta)+O(\theta^2)\right]D_{-\theta/\eta}(\beta/\sqrt{\eta}),
\end{align}
and rearranging terms we obtain
\begin{align}\label{5666}
\sqrt{\eta}\frac{d}{d\beta}D_{-\theta/\eta}(\tfrac{\beta}{\sqrt{\eta}})+\frac{\beta}{2\sqrt{\eta}}D_{-\theta/\eta}(\tfrac{\beta}{\sqrt{\eta}})
=\frac{-\theta}{\sqrt{\eta}}\mathcal{R}(\beta)D_{-\theta/\eta}(\tfrac{\beta}{\sqrt{\eta}})+ O(\theta^2).
\end{align}
Note that $\theta=0$
is a solution to \eqref{5666}. Now consider $\theta<0$ with $-\theta=-q\eta=O(\eta)$. For $\beta<0$ and $\eta\to 0$, \eqref{eqP12} shows that
\begin{align}
D_q(\beta/\sqrt{\eta})=D_q(-|\beta|/\sqrt{\eta})\sim\frac{\sqrt{2\pi}}{\Gamma(-q)}e^{\beta^2/(4\eta)}(-\beta/\sqrt{\eta})^{-q-1}
\end{align}
as long as $q\neq 0,1,2,\ldots$. Using a similar formula for $D'_q(\beta/\sqrt{\eta})$
we then have
\begin{align}
\frac{D'_q(\beta/\sqrt{\eta})}{D_q(\beta/\sqrt{\eta})}\sim\frac{\beta}{2\sqrt{\eta}}, \quad  q\neq 0,1,2,\ldots.
\end{align}
But then \eqref{4g1} cannot be satisfied asymptotically. We conclude that if $\funcDD=0$ has roots in the range $\theta=O(\eta)$ they must occur where $-\theta/\eta=q\approx 1,2,\ldots$ (in addition to the root at $\theta=0$). To obtain the minimal root we examine the range where $q\approx 1$, thus setting $q=1+\varepsilon$ where $\varepsilon=\varepsilon(\eta) \to 0$ (which corresponds to $-\theta=\eta+\eta\varepsilon(\eta))$. Then \eqref{5666} is equivalent to the asymptotic relation
\begin{align}\label{544}
\sqrt{\eta}\frac{d}{d\beta}\left[e^{\frac{\beta^2}{4\eta}}D_{1+\varepsilon}(\tfrac{\beta}{\sqrt{\eta}})\right]
\sim\sqrt{\eta}e^{\frac{\beta^2}{4\eta}}\mathcal{R}(\beta)D_{1+\varepsilon}(\tfrac{\beta}{\sqrt{\eta}}).
\end{align}

For a fixed $\beta<0$ we have $\beta/\sqrt{\eta}\rightarrow-\infty$ and we use  the asymptotic expansion of $D_p(z)$ as $z\rightarrow-\infty$ which applies as $z\to -\infty$
 ($\arg (z)=\pi$) in \eqref{eqP11}, even if $p$ is close to a positive integer. With $p=1+\varepsilon$ we have
\begin{align}\label{5888}
e^{z^2/4}D_{1+\varepsilon}(-|z|)= & \frac{\sqrt{2\pi}|z|^{-\varepsilon-2}}{\Gamma(-1-\varepsilon)}e^{z^2/2}\left[1+O(z^{-2})\right]\nonumber \\
&+z^{1+\varepsilon}\left[1-\frac{\varepsilon(1+\varepsilon)}{2z^2}+O(z^{-4})\right].
\end{align}
For $\varepsilon$ small we furthermore approximate $\Gamma(-1-\varepsilon)$ by $\varepsilon^{-1}+O(1)$, which follows from the Laurent expansion of $\Gamma(z)$ near $z=-1$. We thus obtain from \eqref{5888}
\begin{align}\label{566}
e^{\frac{\beta^2}{4\eta}}D_{1+\varepsilon}(\tfrac{\beta}{\sqrt{\eta}})&\sim \frac{\beta}{\sqrt{\eta}}+\frac{\eta}{\beta^2}\sqrt{2\pi}\varepsilon e^{\frac{\beta^2}{2\eta}}
\end{align}
and since term by term differentiation is permissible (see the discussion in Section \ref{S5})
\begin{align}\label{577}
\sqrt{\eta}\frac{d}{d\beta}\left[e^{\frac{\beta^2}{4\eta}}D_{1+\varepsilon}(\tfrac{\beta}{\sqrt{\eta}})\right]&\sim 1+\sqrt{2\pi}\varepsilon e^{\frac{\beta^2}{2\eta}}\frac{\sqrt{\eta}}{\beta}.
\end{align}
Using \eqref{566} and \eqref{577} in \eqref{544} we see that $\varepsilon$ must be roughly of the order $\exp(-\beta^2/2\eta)$; more precisely,
\begin{align}\label{}
\varepsilon(\eta)\sim-\frac{\beta}{\sqrt{2\pi\eta}}e^{-\frac{\beta^2}{2\eta}}\left[1-\beta\mathcal{R}(\beta)\right]
\end{align}
and then $-\theta=\eta(1+\varepsilon)$ leads to Proposition \ref{lem1}.

We have thus shown that for $\beta<0$ there is a root of $\funcDD=0$ where $\theta=O(\eta)$ with $\theta\sim-\eta$. Now we show that there can be no roots in the range $\theta=o(\eta)$, other than $\theta=0$ (which is an exact root for all $\eta$ and $\beta$). For $\beta<0$, $\beta/\sqrt{\eta}\to -\infty$ and we use \eqref{eqP12} and set $\theta=\theta_*\eta^M$. Thus
\begin{align}\label{66-1}
D_{-\theta_*\eta^{M-1}}(-|\beta|/\sqrt{\eta})
&\sim \exp\Big(\frac{\beta^2}{4\eta}\Big)\frac{\sqrt{2\pi\eta}}{-\beta \Gamma(\theta_*\eta^{M-1})}\nonumber\\
&\sim \theta_*\eta^{M-1}\frac{\sqrt{2\pi\eta}}{-\beta}\exp\Big(\frac{\beta^2}{4\eta}\Big),
\end{align}
since $1/\Gamma(p)=p+O(p^2)$ by Taylor series. Then using an analogous formula for $D'_{-\theta/\eta}(\cdot)$ and approximating $D_{-\theta_*\eta^M}(-\beta)\sim \e^{-\beta^2/4}$ and $D_{-\theta_*\eta^M}'(-\beta)\sim \frac12 \beta \e^{-\beta^2/4}$ we obtain
\begin{align}\label{66-2}
\funcDD= \e^{-\beta^2/4} \theta_*\eta^{M-1/2}\exp\Big(\frac{\beta^2}{4\eta}\Big)\sqrt{2\pi}[1+o(1)].
\end{align}
This asymptotic relation holds for $\theta_*$ fixed, $\eta\to 0$, and the $o(1)$ error term holds uniformly on finite $\theta_*$ intervals. But since the leading term is proportional to $\theta_*$ we again conclude that $\theta_*=0$, contradicting our assumption that a root has $\theta=\Theta(\eta^M)$ for $M>1$. If $M=1$ the first asymptotic relation in \eqref{66-1} still holds and then we would find that $\funcDD$ is proportional to $1/\Gamma(S)$ (with now $\theta=S\eta$), which has roots at $S=0,-1,-2,\ldots$. The first root regains $\theta=0$, the second is the spectral gap we analyzed above, and the roots $\theta\sim-N\eta$  for $N\geq 2$ correspond to the higher eigenvalues.

The expansion in \eqref{66-2} relies only on $\theta/\eta$ being small. For any $\theta=o(\eta)$ we again obtain \eqref{66-2}, with $\theta_*\eta^M$ replaced by $\theta$, and this excludes roots where $\theta/\eta=o(1)$, except for $\theta=0$.

\subsection{Proofs of Propositions \ref{lem3} and \ref{lem5}}
We analyze $\funcDD=0$ for $\eta\rightarrow 0$ and $0<\beta\leq \beta_*$. We can no longer use
\eqref{5111} to approximate the parabolic cylinder function, as we will have $\hat \theta\sim-\beta^2/4$. This corresponds to  two saddle points in \eqref{saddle123} (at $\frac{1}{2}(\beta\pm\sqrt{\beta^2+4\theta})$) coalescing, see the discussion in Section \ref{S5} below \eqref{eqP28}. Now we must approximate $D_{-\theta/\sqrt{\eta}}(\beta/\sqrt{\eta})$ and its derivative by Airy functions. We use the following proposition, which follows from \cite[p. 689]{a&s}, and was discussed in \eqref{eqP19}-\eqref{eqP28}.

\begin{proposition}\label{lempara}
If $A,B\rightarrow\infty$ with $A=-\frac{1}{4}B^2+(\frac{1}{2}B)^{2/3}\delta$ and $\delta=O(1)$,
\begin{align}
D_{-A}(B)=e^{-B^2/8}&\Big(\frac{2}{B}\Big)^A\sqrt{2\pi}\Big(\frac{B}{2}\Big)^{1/3}\nonumber\\
&\times\left[{\rm Ai}(\delta)+\frac{1}{2^{4/3}B^{2/3}}
\left(\delta^2{\rm Ai}(\delta)-2{\rm Ai}'(\delta)\right)+O(B^{-4/3})\right] \label{6222}\\
D_{-A}'(B)=e^{-B^2/8}&\Big(\frac{2}{B}\Big)^A\sqrt{2\pi}\Big(\frac{B}{2}\Big)^{2/3}\nonumber\\
&\times\left[{\rm Ai}(\delta)+\frac{1}{2^{4/3}B^{2/3}}
\left(\delta^2{\rm Ai}'(\delta)-2\delta{\rm Ai}(\delta)\right)+O(B^{-4/3})\right]\label{6333}.
\end{align}
The error terms in \eqref{6222} and \eqref{6333} are uniform on finite $\delta$ intervals.
\end{proposition}

We  let  $A=\theta/\eta$, $B=\beta/\sqrt{\eta}$, $\delta=\chi$, and note that
\begin{align}
\frac{\theta}{\eta}=-\frac{\beta^2}{4\eta}+\frac{\beta^{2/3}\chi}{2^{2/3}\eta^{1/3}} \quad \Rightarrow \quad \theta=-\frac{\beta^2}{4}+
\Big(\frac{\beta}{2}\Big)^{2/3}\eta^{2/3}\chi.
\end{align}
We rewrite $\funcDD=0$ as
\begin{align}\label{6666}
-D_{-\theta}'(-\beta)D_{-\theta/\eta}(\beta/\sqrt{\eta})=\sqrt{\eta}D_{-\theta/\eta}'(\beta/\sqrt{\eta})D_{-\theta}(-\beta)
\end{align}
and recall that, by definition, $\beta_*$ is the minimal root of $D'_{\beta^2/4}(-\beta)=0$. For $-\theta=\beta^2/4+O(\eta^{2/3})$ we use \eqref{6222} and \eqref{6333} in \eqref{6666} and cancel some common factors to obtain
\begin{align}\label{eqq32}
&\left[-D_{\beta^2/4}'(-\beta)+O(\eta^{2/3})\right]
\left[{\rm Ai}(\chi)+2^{-4/3}\beta^{-2/3}\eta^{1/3}(\chi^2{\rm Ai}(\chi)-2{\rm Ai}'(\chi))+O(\eta^{2/3})\right],\nonumber\\
&=\sqrt{\eta}\Big(\frac{\beta}{2}\Big)^{1/3}\eta^{-1/6}
\left[D_{\beta^2/4}(-\beta)+O(\eta^{2/3})\right]\nonumber\\
&\quad \times \left[{\rm Ai}'(\chi)+2^{-4/3}\beta^{-2/3}\eta^{1/3}(\chi^2{\rm Ai}'(\chi)-2\chi{\rm Ai}(\chi))+O(\eta^{2/3})\right]\nonumber\\
&=O(\eta^{1/3}).
\end{align}
The error terms are uniform on finite $\chi$ intervals.
The equation \eqref{eqq32} applies both for $0<\beta<\beta_*$ and $\beta\approx\beta_*$, but its solution is different for these two cases. For $0<\beta<\beta_*$ the first factor in the left-hand side of \eqref{eqq32} is $O(1)$, while it is $o(1)$ if $\beta=\beta_*$ (or $\beta\approx\beta_*$).

First we consider $0<\beta<\beta_*$. The right-hand side of \eqref{eqq32} is $O(\eta^{1/3})$ so that $\chi$ must be such that the left-hand side vanishes (to leading order in $\eta$), which implies that ${\rm Ai}(\chi)=0$. Thus $\chi$ must be close to a root of the Airy function, and the maximal root occurs at $a_0=-2.33810\ldots$.
To obtain a more precise estimate we let $\chi-a_0=\eta^{1/3}\chi_1(\beta,\eta)$ so that ${\rm Ai}(\chi)\sim \eta^{1/3}\chi_1(\beta,\eta){\rm Ai}'(a_0)$ as $\eta\rightarrow 0$. Then \eqref{eqq32} becomes
\begin{align}\label{6888}
&-D_{\beta^2/4}'(-\beta)\eta^{1/3}\left[\chi_1(\beta,\eta){\rm Ai}'(a_0)-2^{-1/3}\beta^{-2/3}a_0{\rm Ai}'(a_0)+O(\eta^{1/3})\right] \nonumber\\
&=-\Big(\frac{\beta}{2}\Big)^{1/3}\eta^{1/3}\left[{\rm Ai}'(a_0)+O(\eta^{1/3})\right]D_{\beta^2/4}(-\beta).
\end{align}
Dividing \eqref{6888} by $\eta^{1/3} $ and letting $\eta\to 0$ we conclude that $\chi_1(\beta,\eta)\to\chi_1(\beta)$ as $\eta\to 0$, with
\begin{align}\label{}
\chi_1(\beta)=a_0\beta^{-2/3}2^{-1/3}-\Big(\frac{\beta}{2}\Big)^{1/3}\frac{D_{\beta^2/4}(-\beta)}{D_{\beta^2/4}'(-\beta)}.
\end{align}
This leads to \eqref{eq:ert} and completes the proof of Proposition \ref{lem3}.


We next consider $\beta\approx\beta_*$. Then $D_{\beta^2/4}(-\beta)\neq 0$ and by Taylor series we have
 \begin{align}\label{7111}
-\frac{D_{\beta^2/4}'(-\beta)}{D_{\beta^2/4}(-\beta)}=-L(\beta-\beta_*)+O((\beta-\beta_*)^2),
\end{align}
where $L$ is given by \eqref{LLL}. Thus if we scale $\beta-\beta_*=\eta^{1/3}W$ the left and right-hand sides of \eqref{eqq32} are both $O(\eta^{1/3})$ and we obtain the limiting equation
\begin{equation}\label{151}
-L\cdot W\cdot{\rm Ai}(\chi)=(\beta_*/2)^{1/3}{\rm Ai}'(\chi),
\end{equation}
which defines $\chi$ implicitly in terms of $W$, and leads to Proposition \ref{lem5}.


Finally, we show that when $\beta>0$ the equation $\funcDD=0$ can have no roots in the range $\theta=o(1)$. We consider scales of the form $\theta=\theta_* \eta^M$ with $M>0$ and exclude the possibility of roots that have $\theta=\Theta(\eta^M)$. We write
\begin{align}\label{67a}
\funcDD&=D_{-\theta/\eta}(\tfrac{\beta}{\sqrt{\eta}})\Big[\frac{d}{d\beta}D_{-\theta}(-\beta)-\sqrt{\eta}D_{-\theta}(-\beta)\frac{D_{-\theta/\eta}'(\tfrac{\beta}{\sqrt{\eta}})}{D_{-\theta/\eta}(\tfrac{\beta}{\sqrt{\eta}})}\Big] \end{align}
and use \eqref{eqP31} which for $\theta,\eta=o(1)$ yields
\begin{align}\label{67b}
\frac{D_{-\theta/\eta}'(\tfrac{\beta}{\sqrt{\eta}})}{D_{-\theta/\eta}(\tfrac{\beta}{\sqrt{\eta}})}=-\frac{\beta}{2\sqrt{\eta}}- \frac{\theta}{\beta\sqrt{\eta}}+O(\theta\sqrt{\eta},\theta^2/\sqrt{\eta}).
\end{align}
Also, using \eqref{eqP8} gives
\begin{align}\label{67c}
\frac{d}{d\beta}D_{-\theta}(-\beta)+\frac\beta 2 D_{-\theta}(-\beta)=\theta D_{-\theta-1}(-\beta)\sim \theta D_{-1}(-\beta).
 \end{align}
Using \eqref{67b} and \eqref{67c} in \eqref{67a} yields, to leading order in $\eta$, for $\theta=\theta_* \eta^M$,
\begin{align}\label{67d}
\funcDD&=\theta_*\eta^{M}D_{-\theta_*\eta^{M-1}}(\tfrac{\beta}{\sqrt{\eta}})\Big[\beta^{-1}\e^{-\beta^2/4}+D_{-1}(\beta)\Big][1+o(1)].
 \end{align}
For $M>1$ and $\eta\to 0$ we can approximate $D_{-\theta_*\eta^{M-1}}(\tfrac{\beta}{\sqrt{\eta}})\sim \exp[-\beta^2/(4\eta)]$, for $M=1$ we have
$D_{-\theta_*\eta^{M-1}}(\tfrac{\beta}{\sqrt{\eta}})\sim \exp[-\beta^2/(4\eta)](\sqrt{\eta}/\beta)^{\theta_*}$, while for $0<M<1$ we must approximate the parabolic cylinder function using the result in \eqref{5111}, which applies for large index and large argument. But in all cases the approximation leads to $\funcDD$ in \eqref{67d} being proportional to $\theta_*$ with a positive multiplier. Thus we again conclude that $\theta_*=0$, contradicting the existence of the root(s) where $\theta=\Theta(\eta^M)$ for any $M>0$. Note that unlike $\beta\leq 0$, the scale $\theta=\Theta(\eta)$ does not lead to roots, but only a change in the expansion of $D_{-\theta_*\eta^{M-1}}(\tfrac{\beta}{\sqrt{\eta}})$.

In obtaining  \eqref{67d} we used  only the fact that $\theta=o(1)$. For any $\theta=o(1)$, \eqref{67d} holds, with $\theta_*\eta^M$ replaced by $\theta$. Again, the expansion of $D_{-\theta/\eta}(\beta/\sqrt{\eta})$ will be different as $\theta/\eta\to 0$,  $\theta/\eta\to \infty$, or  $\theta/\eta=\Theta(1)$, but the first multiplicative factor in \eqref{67d} ($=\theta_*\eta^M=\theta$) shows that there can be no roots in any range where $\theta=o(1)$, except for the root at $\theta=0$.



\section{Monotonicity of the spectral gap}\label{S8}
The surface sketched in Figure~\ref{fig11} suggested
certain monotonicity properties of $r(\beta,\eta)$,
and these were partially confirmed by
the various asymptotic formulas in Section \ref{sec:main}.
We now establish these analytically,
for all values of $(\beta,\eta)$. We shall obtain:
\begin{proposition}\label{propM}
Let $\sgn(z)=+1$ if $z>0$,
$\sgn(z)=-1$ if $z<0$ and $\sgn(0)=0$. Then
\begin{equation}\label{eqM1}
\sgn\left(\dfrac{\partial r}{\partial\beta}\right)=-\sgn(\eta-1).
\end{equation}
Hence, for a fixed $\eta<1$ the spectral gap~$r$
is an increasing function of~$\beta$, while it decreases
with $\beta$ for fixed $\eta>1$. If $\eta=1$, $r(\beta,1)=1$
is constant.
\end{proposition}

To establish this result it is useful to
set
\begin{equation}\label{eqM2}
V(P;\eta,\beta)=\Vcal(\Ocal;\eta,\beta),\ \Ocal=-P
\end{equation}
and then in view of Theorem~\ref{thmspec}
\begin{equation}\label{eqM3}
V(r(\beta,\eta);\eta,\beta)=0.
\end{equation}
By implicit differentiation of (\ref{eqM3}) we have
\begin{equation}\label{eqM4}
\left.\dfrac{\partial V}{\partial P}\right\vert_{P=r}\cdot \dfrac{\partial
r}{\partial \beta}+\left.\dfrac{\partial V}{\partial
\beta}\right\vert_{P=r}=0.
\end{equation}
By definition, $r$ is the minimal positive
solution of $V(P;\eta,\beta)=0$ and we also note
that $V(0;\eta,\beta)=0$, as $\Ocal =0$ is a simple
pole of $\widehat{p}(x;\Ocal)$ in Theorem~\ref{greenpos},
which corresponds to the steady state limit
in~(\ref{steadyssty}). Thus $P=0$ and $P=r$
are consecutive zeros of $V(P;\eta,\beta)=0$.
To determine the sign of $\partial r/\partial \beta$ in~(\ref{eqM4})
requires
that we know the signs of $\partial V/\partial P$ and
$\partial V/\partial\beta$
at $P=r$. For the former we can compute
$\partial V/\partial P$ using the expression in Theorem~\ref{thmspec} and
the integral respresentations in~(\ref{eqP2}), (\ref{eqP5})
and~(\ref{eqP6}).
However, an indirect argument leads immediately
to the value of $\sgn(\left.\partial V/\partial P\right\vert_{P=r})$.

The solutions of $V=0$ for $P\geq 0$ correspond to
poles of the Laplace transform $\widehat{p}(x;\Ocal)$ in
Theorem~\ref{greenpos} and these are
the eigenvalues $\lambda_\N$ for $\N\geq 0$, with $\lambda_0=0$ and
$\lambda_1=r$. From the general theory of the
one-dimensional Schr\"{o}dinger equation, the
eigenvalues are all simple (see the discussion
below (\ref{eqS7})) and the equation $V=0$
has simple zeros. Hence, $\left.\partial V/\partial P\right\vert_{P=r}\neq 0$.
We also note that if $V$ had, say, a double
zero at $P=r$, then (\ref{green1}) would imply that
the spectral expansion of $p(x,t)$ would involve
the terms $e^{-\lambda_1t}=e^{-rt}$ and also $te^{-rt}$, and
this would contradict the self-adjointness of
the Schr\"{o}dinger equation in~(\ref{eqS4}). Now, since
$P=0$ and $P=r$ are consecutive simple zeros
on the real axis of the entire function
$V$ (as a function of $P$) we must have
\begin{equation}\label{eqM5}
\sgn\left(\left.\dfrac{\partial V}{\partial
P}\right\vert_{P=r}\right)=-\sgn\left(\left.\dfrac{\partial V}{\partial
P}\right\vert_{P=0}\right).
\end{equation}
Computing the right-hand side of (\ref{eqM5}) is much
easier than computing the left-hand side, as we
show below.

We define the functions $I(z)$ and $J(z)$ by
\begin{align}\label{eqM6}
I(z)&=\left.\dfrac{\partial}{\partial
P}D_P(z)\right\vert_{P=0}=\dfrac{e^{z^2/4}}{i\sqrt{2\pi}}\int_{\Br}(\log u) e^{-zu}e^{u^2/2}\,
du,\\
J(z)&=\dfrac{e^{z^2/4}}{i\sqrt{2\pi}}\int_{\Br}u(\log u)e^{-zu}e^{u^2/2}\,
du.\notag
\end{align}
Then we expand $V$ in (\ref{eqM2}) in Taylor series about
$P=0$ to obtain
\begin{align}\label{eqM7}
V={}&-\sqrt{\eta}\left[D_0(-\beta)+PI(-\beta)+\Oup (P^2)\right]\\
&\times\left[D'_0\left(\dfrac{\beta}{\sqrt{\eta}}\right)+\left.\left(\dfrac{z}{2}I(z)-J(z)\right)\right\vert_{z=\beta/\sqrt{\eta}}\dfrac{P}{\eta}+O(P^2)\right]\notag\\
&-
\left[D'_0(-\beta)+\left.\left(\dfrac{z}{2}I(z)-J(z)\right)\right\vert_{z=-\beta}P+\Oup(P^2)\right]\notag\\
&\times\left[D_0\left(\dfrac{\beta}{\sqrt{\eta}}\right)+I\left(\dfrac{\beta}{\sqrt{\eta}}\right)\dfrac{P}{\eta}+\Oup(P^2)\right].\notag
\end{align}
Using $D_0(z)=e^{-z^2/4}$, $D'_0(z)=-\dfrac{z}{2}e^{-z^2/4}$ we obtain
from~(\ref{eqM7})
\begin{align}\label{eqM8}
\left.\dfrac{\partial V}{\partial
P}\right\vert_{P=0}={}&\exp\left(-\dfrac{\beta^2}{4\eta}\right)\left[\beta
I(-\beta)+J(-\beta)\right]\\
&-\dfrac{1}{\sqrt{\eta}}e^{-\beta^2/4}\left[\dfrac{\beta}{\sqrt{\eta}}I\left(\dfrac{\beta}{\sqrt{\eta}}\right)-J\left(\dfrac{\beta}{\sqrt{\eta}}\right)\right].\notag
\end{align}
But an integration by parts shows that
\begin{align}\label{eqM9}
zI(z)-J(z)
&=\dfrac{e^{z^2/4}}{i\sqrt{2\pi}}\int_{\Br}
(z-u)\log(u)e^{-zu}e^{u^2/2}\, du\\
&=\dfrac{e^{z^2/4}}{i\sqrt{2\pi}}\int_{\Br}\log(u)\,
d\left(-e^{-zu}e^{u^2/2}\right)\notag\\
&=\dfrac{e^{z^2/4}}{i\sqrt{2\pi}}\int_{\Br}\dfrac{e^{-zu}}{u}e^{u^2/2}\,
du\notag\\
&=e^{z^2/4}\int_z^{\infty} e^{-\xi^2/2}\, d\xi.\notag
\end{align}
Here we also used the fact that a parabolic
cylinder function of order $P=-1$ can be expressed
in terms of the standard error function, or
probability integral. Using~(\ref{eqM9}) in~(\ref{eqM8}) we have
\begin{multline}\label{eqM10}
\left.\dfrac{\partial V}{\partial P}\right\vert_{P=0}=
-\Bigg[\exp\left(\dfrac{\beta^2}{4}-\dfrac{\beta^2}{4\eta}\right)\int^{\infty}_{-\beta}e^{-\xi^2/2}\,
d\xi
+\exp\left(-\dfrac{\beta^2}{4}+\dfrac{\beta^2}{4\eta}\right)\int^{\infty}_{\beta/\sqrt{\eta}}e^{-\xi^2/2}\,
d\xi\Bigg]
\end{multline}
\begin{setlength}{\multlinegap}{77pt}
\begin{multline*}
=-\exp\left(-\dfrac{\beta^2}{4\eta}-\dfrac{\beta^2}{4}\right)\Bigg[\int^{\infty}_0e^{-\beta\xi}e^{-\eta\xi^2/2}\,
d\xi
+\int^0_{-\infty}e^{-\beta\xi}e^{-\xi^2/2}\, d\xi\Bigg]
\end{multline*}
\end{setlength}
so that $\left.\partial V/\partial P\right\vert_{P=0}<0$ and hence, in
view of~(\ref{eqM4}) and~(\ref{eqM5}),
\begin{equation}\label{eqM11}
\sgn\left(\dfrac{\partial r}{\partial \beta}\right)=-\sgn\left(\left.\dfrac{\partial V}{\partial \beta}\right\vert_{P=r}\right).
\end{equation}
Now,
\begin{align}\label{eqM12}
\dfrac{\partial V(P;\eta,\beta)}{\partial\beta}={}&\sqrt{\eta}D'_P(-\beta)D'_{P/\eta}\left(\dfrac{\beta}{\sqrt{\eta}}\right)\\
&-D_P(-\beta)D''_{P/\eta}\left(\dfrac{\beta}{\sqrt{\eta}}\right)+D''_P(-\beta)D_{P/\eta}\left(\dfrac{\beta}{\sqrt{\eta}}\right) \notag\\
&-\dfrac{1}{\sqrt{\eta}}D'_P(-\beta)D'_{P/\eta}\left(\dfrac{\beta}{\sqrt{\eta}}\right).\notag
\end{align}
Using the parabolic cylinder equation $D''_P(z)=
\big(\frac{1}{4}z^2-P-\frac{1}{2}\big)D_P(z)$ we can simplify~(\ref{eqM12})
to
\begin{multline}\label{eqM13}
\dfrac{\partial
V}{\partial\beta}=\dfrac{\eta-1}{\eta}\Bigg[\sqrt{\eta}D'_P(-\beta)D'_{P/\eta}\left(\dfrac{\beta}{\sqrt{\eta}}\right)
+\left(\dfrac{\beta^2}{4}-P\right)D_P(-\beta)D_{P/\eta}\left(\dfrac{\beta}{\sqrt{\eta}}\right)\Bigg].
\end{multline}
When $P=r$ we can further use the fact
that $V(r;\eta,\beta)=0$ to simplify the right
side of~(\ref{eqM13}). We will need to separately
consider the two cases $D_r(-\beta)=0$ (a
degenerate case that occurs rarely) and $D_r(-\beta)\neq 0$
(which is typical).

In the degenerate case we have
\begin{equation}\label{eqM14}
\sgn\left(\left.\dfrac{\partial V}{\partial\beta}\right\vert_{P=r}\right)=\sgn(\eta-1)\cdot\sgn\left(D'_r(-\beta)\right)
\cdot\sgn\left(D'_{r/\eta}\left(\dfrac{\beta}{\sqrt{\eta}}\right)\right),
\end{equation}
and thus
\begin{equation}\label{eqM15}
\sgn\left(\dfrac{\partial r}{\partial
\beta}\right)=-\sgn(\eta-1)\sgn\left(D'_r(-\beta)\right)\sgn\left(D'_{r/\eta}\left(\dfrac{\beta}{\sqrt{\eta}}\right)\right).
\end{equation}
We note that if $D_r(-\beta)=0$ then certainly
$D'_r(-\beta)\neq 0$, as discussed in Section~\ref{S5} below
\eqref{eqP9A}. But if both $V=0$ and $D_r(-\beta)=0$ then
certainly $D_{r/\eta}(\beta/\sqrt{\eta})=0$. In the non-degenerate
case we can rewrite $V=0$ as
\begin{equation}\label{eqM16}
\dfrac{\sqrt{\eta}D'_{r/\eta}\left(\dfrac{\beta}{\sqrt{\eta}}\right)}{D_{r/\eta}\left(\dfrac{\beta}{\sqrt{\eta}}\right)}=-\dfrac{D'_r(-\beta)}{D_r(-\beta)}.
\end{equation}
Using (\ref{eqM16}) in (\ref{eqM13}) and (\ref{eqM11}) we conclude that
\begin{align}\label{eqM17}
\sgn\left(\dfrac{\partial r}{\partial
\beta}\right)={}&\sgn(\eta-1)\cdot\sgn\left(D_r(-\beta)\right)\\
&\cdot\sgn\left(D_{r/\eta}\left(\dfrac{\beta}{\sqrt{\eta}}\right)\right)\cdot\sgn\left\{\left(r-\dfrac{\beta^2}{4}\right)\left(D_r(-\beta)\right)^2+\left(D'_r(-\beta)\right)^2\right\}.\notag
\end{align}
We proceed to determine the signs of the
various terms in~(\ref{eqM15}) and~(\ref{eqM17}).

It proves useful to understand the
behaviors of $r(\beta,\eta)$ as $\beta\to \pm\infty$ for
a fixed~$\eta$. By a calculation completely
analogous to that used to establish
Proposition~\ref{lem1}, we find that
\begin{equation}\label{eqM18}
r(\beta,\eta)-\eta\sim\dfrac{(1-\eta)\sqrt{\eta}}{(-\beta)\sqrt{2\pi}}\exp\left(-\dfrac{\beta^2}{2\eta}\right);\
\eta\neq 1,\ \beta\to-\infty.
\end{equation}
Whereas Proposition~\ref{lem1} applies for $\eta\to 0$
with fixed $\beta<0$, (\ref{eqM18}) applies for fixed $\eta$ as
$\beta\to-\infty$. By expanding (\ref{eq:ertlem1}) for $\beta\to-\infty$
and~(\ref{eqM18}) as $\eta\to 0$ we see that the two
agree in this intermediate limit. Thus for
$\beta$ large and negative, $r$ is exponentially
close to~$\eta$, as could be expected since
then almost all of the probability mass in the
model migrates to the range $x>0$.
A similar analysis as $\beta\to+\infty$ shows
that
\begin{equation}\label{eqM19}
r(\beta,\eta)-1\sim\dfrac{\eta-1}{\beta\sqrt{2\pi}}e^{-\beta^2/2};\
\eta\neq 1,\ \beta\to+\infty,
\end{equation}
which can also be obtained simply by using
(\ref{eqM18}) and the symmetry relation in~\eqref{sym11}.
For $\beta$ large and positive the probability mass
migrates to the region $x<0$. Note that
(\ref{eqM18}) and~(\ref{eqM19}) suggest that $\partial r/\partial \beta$
has the
oppositive sign as $\eta-1$, at least for
$|\beta|$ sufficiently large, and this we proceed
to establish for any~$\beta$.

Returning to (\ref{eqM17}) we proceed to show
that $D_r(-\beta)$ and $D_{r/\eta}(\beta/\sqrt{\eta})$ always
have opposite signs. In general, suppose
that we have two real analytic functions $F(x)$
and~$G(x)$. The ratio $F(x)/G(x)$ can only change sign at a point $x_*$ where
$F(x_*)=0$ with $G(x_*)\neq 0$, or $G(x_*)=0$ with
$F(x_*)\neq 0$, or, possibly, where $F$ and $G$ both vanish
but have zeros of different orders.
Thus if $F(x)$ and $G(x)$ are non-zero, or
if their only zero(s) coincide and they
are of the same order, then $F(x)/G(x)$ cannot
change sign. Then determining $\sgn(F(x)/G(x))$
requires only that we evaluate the ratio at
a particular point, which could be $x=\pm\infty$.
But, we showed below (\ref{eqM15}) that $D_r(-\beta)$
and $D_{r/\eta}(\beta/\sqrt{\eta})$, as functions of~$\beta$ for a fixed
$\eta>0$, can only vanish simultaneously. Furthermore,
if $D_r(-\beta)$ vanishes at some~$\beta_c$, then for
$\beta$ near $\beta_c$ but $\beta\neq \beta_c$ we can rewrite the
equation $\V(r;\eta,\beta)=0$ as
\begin{equation}\label{eqM20}
\dfrac{D_r(-\beta)}{D_{r/\eta}(\beta/\sqrt{\eta})}=-
\dfrac{D'_r(-\beta)}{\sqrt{\eta}D'_{r/\eta}(\beta/\sqrt{\eta})}.
\end{equation}
Then letting $\beta\to \beta_c$ leads to
\begin{equation}\label{eqM21}
\lim_{\beta\to\beta_c}\left[\dfrac{D_r(-\beta)}{D_{r/\eta}(-\beta/\sqrt{\eta})}\right]=-\dfrac{D'_{r_c}(-\beta_c)}{\sqrt{\eta}D'_{r_c/\eta}(\beta_c/\sqrt{\eta})}
\end{equation}
where $r_c=r(\beta_c,\eta)$. But (\ref{eqM21}) shows that $D_r(-\beta)$
and $D_{r/\eta}(\beta/\sqrt{\eta})$ must vanish to the same order
at $\beta=\beta_c$ (in fact they must have simple zeros
there). Thus we conclude that
$D_r(-\beta)/D_{r/\eta}(\beta/\sqrt{\eta})$ cannot change sign. To
determine this constant sign we can let
either $\beta\to+\infty$ or $\beta\to-\infty$ as then we
have asymptotic formulas for $r$. Using
(\ref{eqM19}), (\ref{eqP12}) and the fact that
$D_1(-\beta)=-\beta e^{-\beta^2/4}$
we find that
\begin{equation}\label{eqM22}
\dfrac{D_{r/\eta}(\beta/\sqrt{\eta})}{D_r(-\beta)}\sim
-\beta^{\frac{1}{\eta}-1}\eta^{-\frac{1}{2\eta}}\exp\left[\dfrac{1}{4}\left(1-\dfrac{1}{\eta}\right)\beta^2\right],\
\beta\to+\infty
\end{equation}
and thus $D_{r/\eta}(\beta/\sqrt{\eta})$ and $D_r(-\beta)$ have
$\underline{\text{opposite}}$ signs as $\beta\to+\infty$, and thus this
is true for all~$\beta$. Note also that for $\beta\to-\infty$,
(\ref{eqM18}) and (\ref{eqP10}) lead to
\begin{equation}\label{eqM23}
\dfrac{D_{r/\eta}(\beta/\sqrt{\eta})}{D_r(-\beta)}\sim-(-\beta)^{1-\eta}\dfrac{1}{\sqrt{\eta}}\exp\left[\dfrac{1}{4}\left(1-\dfrac{1}{\eta}\right)\beta^2\right],\
\beta\to-\infty
\end{equation}
and this verifies the conclusion about opposite
signs. We have thus simplified (\ref{eqM17}) to
\begin{equation}\label{eqM24}
\sgn\left(\dfrac{\partial r}{\partial
\beta}\right)=-\sgn(\eta-1)\sgn\left\{\left(r-\dfrac{\beta^2}{4}\right)\left(D_r(-\beta)\right)^2+\left(D'_r(-\beta)\right)^2\right\}
\end{equation}
in the non-degenerate case. In the degenerate
case we conclude from (\ref{eqM21}) that $D'_{r_c}(-\beta_c)$ and
$D'_{r_c/\eta}(\beta_c/\sqrt{\eta})$ have the same sign, and then
(\ref{eqM15}) shows that $\sgn(\partial r/\partial\beta)=-\sgn(\eta-1)$,
which
establishes Proposition~\ref{propM}.

It remains to show that the last
factor in (\ref{eqM24}) has always positive sign.
%
Let us define
\begin{equation}\label{eqM25}
\Hcal(P,z)=\left(P-\dfrac{z^2}{4}\right)D_P^2(z)+\left[D'_P(z)\right]^2,
\end{equation}
and we consider $\Hcal$ as a function of both
$P$ and $z$. We clearly have $\Hcal(P,0)=
PD_P^2(0)+\left[D'_P(0)\right]^2>0$ for $P>0$, with
$\Hcal(0,0)=0$. Also, $\Hcal(0,z)=0$ for all $z$,
in view of (\ref{eqP4}). We consider $P>0$ and $z>0$.

We shall show that $\Hcal(P,z)>0$
for all $z\geq 0$ when $P>0$. For $z\to \infty$
the estimate in (\ref{eqP31}) leads to
\begin{equation}\label{eqM26}
\left(\dfrac{D'_P(z)}{D_P(z)}\right)^2=\dfrac{z^2}{4}-P+\dfrac{P}{z^2}+\Oup\left(z^{-4}\right)
\end{equation}
and then from (\ref{eqM25}) and (\ref{eqP10})
\begin{equation}\label{eqM27}
\Hcal(P,z)\sim D^2_P(z)\dfrac{P}{z^2}\sim Pz^{2P-2}e^{-z^2/2},\ z\to+\infty
\end{equation}
so that $\Hcal$ is positive for $z$ sufficiently large.
By differentiating (\ref{eqM25}) with respect to~$z$
we obtain
\begin{align}\label{eqM28}
\dfrac{\partial \Hcal}{\partial
z}={}&2D_P(z)D'_P(z)\left(P-\dfrac{z^2}{4}\right)-\dfrac{z}{2}D^2_P(z)+2D''_P(z)D'_P(z)\\
={}&-\dfrac{z}{2}D^2_P(z)-D_P(z)D'_P(z)\notag\\
={}&-PD_P(z)D_{P-1}(z).\notag
\end{align}
Here we also used (\ref{eqP1}) and the recurrence (\ref{eqP8}).
From (\ref{eqM28}) we conclude that $\Hcal$ has maximum or
minimum values at roots of $D_P(z)$ and
$D_{P-1}(z)$, as functions of~$z$. As discussed in
Section~\ref{S5}, there are at most finitely many of
these. But if $D_P(z_*)=0$ for some $z_*$ then
$D'_P(z_*)\neq0$ and $\Hcal(P,z_*)=\left[D'_P(z_*)\right]^2>0$.
If $D_{P-1}(\widetilde{z})=0$ for some $\widetilde{z}$ then (\ref{eqM25})
and
(\ref{eqP8}) show that
\begin{align}\label{eqM29}
\Hcal(P,\widetilde{z})&=\left(P-\dfrac{z^2}{4}\right)D_P^2(\widetilde{z})+\left[D'_P(\widetilde{z})\right]^2\\
&=PD^2_P(\widetilde{z}),\notag
\end{align}
which is again positive for $P>0$. Note that
we cannot have simultaneously
$D_{P-1}(\widetilde{z})=0=D_P(\widetilde{z})$,
for then (\ref{eqP8}) would imply that $D'_P(\widetilde{z})=0$ also.
We have thus shown that $\Hcal(P,z)$ is (for
$P>0$) positive at $z=0$ and as $z\to+\infty$,
and also $\Hcal>0$ at any maximum/minimum
value of~$\Hcal$. We then conclude that
\begin{equation}\label{eqM30}
\Hcal(P,z)>0\ \text{for}\ P>0\ \text{and}\ z\geq 0.
\end{equation}
Note that if $\Hcal$ becomes negative at
some $z=z'$ then $\Hcal$ would need to reach a
minimum value at a point $z''$ where $\Hcal<0$, since
for sufficiently large~$z$ we again have $\Hcal>0$.

Now we let $P=r(\beta,\eta)>0$ and $z=-\beta$ and
use (\ref{eqM25}) and (\ref{eqM30}) in (\ref{eqM24}) to conclude
that
\[
\sgn\left(\dfrac{\partial r}{\partial\beta}\right)=-\sgn(\eta-1),\ \beta\leq 0
\]
and we have thus
established Proposition~\ref{propM} for $\beta\leq 0$ and
all $\eta>0$. To show the result holds also for
$\beta>0$, we need only use the symmetry
relation in (\ref{sym11}), which shows that if $r$
increases with $\beta$ for $\beta<0$ and $0<\eta<1$, (resp.\ $\eta>1$) then
$r$ will decrease with $\beta$ for $\beta>0$ and $\eta>1$ (resp.\
$0<\eta<1$).

Alternately, we can use the
relation (\ref{eqM20}) (in the non-degenerate case)
and (\ref{eqM24}) to conclude that
\begin{align}\label{eqM31}
\sgn\left(\dfrac{\partial r}{\partial
\beta}\right)= -\sgn(\eta-1)\sgn\Bigg\{\left(r-\dfrac{\beta^2}{4}\right)\left[D_{r/\eta}\left(\dfrac{\beta}{\sqrt{\eta}}\right)\right]^2
+\eta\left[D'_{r/\eta}\left(\dfrac{\beta}{\sqrt{\eta}}\right)\right]^2\Bigg\}
\end{align}
and apply (\ref{eqM30}) with $P=r/\eta>0$ and $z=\beta/\sqrt{\eta}>0$.
This concludes the proof of Proposition~\ref{propM}, which
was suggested by our numerical and asymptotic
results.

\section*{Acknowledgments}
The work of Charles Knessl was supported partially by NSF grant DMS-05-03745 and NSA grants H 98230-08-1-0102 and H 98230-11-1-0184. The work of Johan van Leeuwaarden was supported by an NWO (The Netherlands Organization for Scientific Research) Veni grant and an ERC (European Research Council) starting grant. We thank an anonymous referee for many constructive suggestions.

\appendix

\section{Appendix A}\label{PROP6}

We now prove Proposition \ref{propAA}.
First consider $\gamma\to-\infty$.
Since (\ref{eqrt}) is equivalent to $D_{R-1}(\gamma)=0$
for $R\neq0$, the asymptotic formula in
(\ref{eqP12}) shows that $D_{R-1}(\gamma)$ is positive for
$\gamma\to-\infty$, unless $\Gamma(1-R)$ is singular, and
the minimal singularity occurs at $R=1$.
For $R$ close to $1$ we must use the asymptotic
formula in (\ref{eqP11}), with $\arg(z)=\arg(\gamma)=\pi$.
Then $D_{R-1}(\gamma)=0$ implies that
\begin{equation}\label{eqAP61}
e^{-\gamma^2/4}\gamma^{R-1}\sim\sqrt{2\pi}(-\gamma)^{-R}e^{\gamma^2/4}\dfrac{-1}{\Gamma(1-R)},\
\gamma\to -\infty.
\end{equation}
Then using $\Gamma(1-R)=(1-R)^{-1}+\Oup(1)$ as $R\to 1$,
(\ref{eqAP61}) leads to
\begin{equation}\label{eqAP62}
R-1\sim \dfrac{-\gamma e^{-\gamma^2/2}}{\sqrt{2\pi}},\ \gamma\to-\infty
\end{equation}
which is the result in (\ref{p61}).

Now consider the limit $\gamma\to+\infty$.
For a fixed $R$ the asymptotic formula in
(\ref{eqP10}) shows that $D_{R-1}(\gamma)\sim \gamma^{R-1}e^{-\gamma^2/4},\
\gamma\to +\infty$,
which is strictly positive. Thus to capture the
zeros of $D_{R-1}(\gamma)$ in this limit $R$ must be
itself large, so that we enter the oscillatory
range of the special function. As discussed
in Section~\ref{S5} the transition to oscillatory
behavior occurs when $R\approx\gamma^2/4$ and
then we can approximate $D_{R-1}(\gamma)$ by
Airy functions, with the leading term given
in (\ref{eqP28}). Thus with $R=\gamma^2/4-(\gamma/2)^{2/3}\delta$ and
$\gamma\to \infty$ the minimal root corresponds to
the maximal root of $\Ai(\delta)=0$, which
occurs at $\delta=a_0$, leading to (\ref{p62}).

\section{Appendix B}\label{T8}

We discuss the singularities of (\ref{green1})--(\ref{green3})
in the complex $\theta$-plane and thus establish Proposition \ref{propBB}. As discussed in
Section~\ref{S5}, $D_{-\theta}(\beta)$ is an entire function of~$\theta$,
so that the only singularities of~(\ref{green1})
are the zeros of $\Vcal(\theta;\eta,\beta)$.
The existence of an infinite sequence of
zeros and the fact that they lie on the
real axis $(\text{\upshape Im}(\theta)=0)$ follows from standard
ODE theory, which was discussed in Section~\ref{S1}.

Now consider (\ref{green2}) and (\ref{green3}). The factor
$\Gamma(\theta)$ has simple poles at $\theta=0,-1,-2,-3,\dots$.
If $\theta=-M$, $M\geq 0$ we can simplify $\Vcal$
by using (\ref{eqP4}), so that
\begin{align}\label{eqA81}
\Vcal (-M;\eta,\beta) ={}&-\sqrt{\eta}\He_M(-\beta)e^{-\beta^2/4}D'_{M/\eta}\left(\dfrac{\beta}{\sqrt{\eta}}\right)\\
&+\dfrac{d}{d\beta}\left[\He_M(-\beta)e^{-\beta^2/4}\right] D_{M/\eta}\left(\dfrac{\beta}{\sqrt{\eta}}\right).\notag
\end{align}
Similary, (\ref{defMM}) and (\ref{eqP4}) lead to
\begin{align}\label{eqA82}
\Mcal(-M;\eta,\beta)={}&\sqrt{\eta}\He_M(\beta)e^{-\beta^2/4}D'_{M/\eta}\left(\dfrac{\beta}{\sqrt{\eta}}\right)\\
&-\dfrac{d}{d\beta}\left[\He_M(-\beta)e^{-\beta^2/4}\right]D_{M/\eta}\left(\dfrac{\beta}{\sqrt{\eta}}\right).\notag
\end{align}
Then we use $\He_M(-\beta)=(-1)^M\He_M(\beta)$, as the
Hermite polynomials are odd/even functions
according as~$M$ is odd/even, and comparing
(\ref{eqA81}) to (\ref{eqA82}) we find that
\begin{equation}\label{eqA83}
\Vcal(-M;\eta,\beta)=(-1)^{M+1}\Mcal(-M;\eta,\beta);\ M=0,1,2,\dots.
\end{equation}

From (\ref{eqA83}) we conclude that either
$\Vcal$ and $\Mcal$ are both zero at $\theta=-M$, or
neither is zero. If $\theta=-M$ and $\Vcal(-M;\eta,\beta)
\neq 0$, then $\Mcal/\Vcal=(-1)^{M+1}$. But then
\begin{multline*}
D_M(x_0+\beta)+D_M(-x_0-\beta)\dfrac{\Mcal(-M;\eta,\beta)}{\Vcal(-M;\eta,\beta)}\\
=e^{-(\beta+x_0)^2/4}\left[\He_M(x_0+\beta)+(-1)^{M+1}\He_M(-x_0-\beta)\right]
=0
\end{multline*}
so that the last factor in (\ref{green2}) vanishes,
and thus $\theta=-M$ is \emph{not} a pole of (\ref{green2}).
Similarly, $\theta=-M$ will not be a pole of~(\ref{green3}).
This shows that $\theta=-M$ can only be a
pole of~(\ref{green2}) and~(\ref{green3}) if $\theta=-M$ and
$\Vcal$ and $\Mcal$ simultaneously vanish.

Conversely, suppose that $\Vcal$ and $\Mcal$
both vanish, say at some $\theta=\theta_*$. First
we assume that none of the three
$D_{-\theta_*}(-\beta)$, $D_{-\theta_*}(\beta)$ and
$D_{-\theta_*/\eta}(\beta/\sqrt{\eta})$
are zero. Then the equations $\Mcal=\Vcal=0$
may be rearranged to give
\begin{equation}\label{eqA84}
-\dfrac{D'_{-\theta_*}(-\beta)}{D_{-\theta_*}(-\beta)}=\dfrac{D'_{-\theta_*}(\beta)}{D_{-\theta_*}(\beta)}=\sqrt{\eta}\dfrac{D'_{-\theta_*/\eta}(\beta/\sqrt{\eta})}{D_{-\theta_*/\eta}(\beta/\sqrt{\eta})}.
\end{equation}
But the first equality in (\ref{eqA84}), along
with the Wronkskian identity in (\ref{eqP9}), leads to
\begin{equation}\label{eqA85}
0=-D'_{-\theta_*}(-\beta)D_{-\theta_*}(\beta)-D_{-\theta_*}(-\beta)D'_{-\theta_*}(\beta)=\dfrac{\sqrt{2\pi}}{\Gamma(\theta_*)}.
\end{equation}
But then $\Gamma(\theta_*)$ must be infinite, which
leads us back to the case $-\theta=M=0,1,2\dots$
which we already discussed.

Finally suppose that $D_{-\theta_*}(-\beta)=0$. Then
certainly $D'_{-\theta_*}(-\beta)\neq 0$, and $\Vcal=0$
implies that $D_{-\theta_*/\eta}(\beta/\sqrt{\eta})=0$, and then
$\Mcal=0$ leads to $D_{-\theta_*}(\beta)=0$. But
then (\ref{eqA85}) again leads to the conclusion
that $\theta_*=-M=0,-1,-2,\dots$.
Starting with the assumption that $D_{-\theta_*}(\beta)=0$
or $D_{-\theta_*/\eta}(\beta/\sqrt{\eta})=0$ leads ultimately
to the conclusion that all three denominators
in (\ref{eqA84}) must vanish, and then again (\ref{eqA85}) leads
to $\theta_*=-M$.

We have thus shown that simultaneous
zeros of $\Vcal$ and $\Mcal$ can occur only if
$\theta=0,-1,-2,\dots$. Indeed this does occur precisely
when $\theta=0$, which corresponds to the steady state
limit $p(x,\infty)$. In all cases we showed that
a singularity of (\ref{green2}) or (\ref{green3}) necessarily
has $\Vcal =0$. Thus the equation in Theorem~\ref{thmspec}
captures all of the singularities of $\widehat{p}(x;\theta)$.

\section{Appendix C}\label{3.1}

We consider the expressions in (\ref{green1})--(\ref{green3}),
in the limit $\eta\to 0^+$, and thus establish Proposition \ref{propCC}. In this limit
we can simplify $D_{-\theta/\eta}\big((\eta x+\beta)/\sqrt{\eta}\big)$
using the
asymptotic formula in (\ref{eqP17}). We can
also use (\ref{eqP5}) to obtain an analogous formula
for $D'_{-\theta/\eta}(\beta/\sqrt{\eta})$. Let us assume first that
$\theta$ is positive and real. The expansion
(\ref{eqP17}) follows from a saddle point approximation
to~(\ref{eqP2}), as discussed in~(\ref{eqP17}) and below. The
expansion of $D'_p(z)$ for $z,p\to\infty$ with
$z^2/p>4$ is the same as that of $D_p(z)$,
except that the factor $z/2-u$ in the integrand
in~(\ref{eqP5}) becomes frozen at the saddle
$u=zv_+=\frac{1}{2}\big[z+\sqrt{z^2-4p}\big]$. It follows that
\begin{equation}\label{eq311}
\dfrac{D'_p(z)}{D_p(z)}\sim-\dfrac{1}{2}\sqrt{z^2-4p};\ z,p\to \infty,\
\dfrac{z^2}{p}> 4
\end{equation}
and hence (setting $z=\beta/\sqrt{\eta}$ and $p=-\theta/\eta$ with $\eta\to
0^+$)
\begin{equation}\label{eq312}
\dfrac{
\sqrt{\eta}
D'_{-\theta/\eta}
\left(
\dfrac{\beta}{\sqrt{\eta}}\right)}{D_{-\theta/\eta}\left(\dfrac{\beta}{\sqrt{\eta}}\right)}
\to
-\dfrac{1}{2}\sqrt{\beta^2+4\theta}\ \text{as}\ \eta\to 0^+.
\end{equation}
A similar argument shows that
\begin{equation}\label{eq313}
\dfrac{D_{-\theta/\eta}\left(\dfrac{\beta+\eta
x}{\sqrt{\eta}}\right)}{D_{-\theta/\eta}\left(\dfrac{\beta}{\sqrt{\eta}}\right)}\to
e^{-x\sqrt{\theta+\beta^2/4}}\ \text{as}\ \eta\to 0^+.
\end{equation}
Then writing $\Vcal$ as
\begin{align}\label{eq314}
\Vcal&=-D_{-\theta/\eta}\left(\dfrac{\beta}{\sqrt{\eta}}\right)\left[D'_{-\theta}(-\beta)+\sqrt{\eta}\dfrac{D'_{-\theta/\eta}\left(\beta/\sqrt{\eta}\right)}{D_{-\theta/\eta}\left(\beta/\sqrt{\eta}\right)}D_{-\theta}(-\beta)\right]\\
&\sim D_{-\theta/\eta}\left(\dfrac{\beta}{\sqrt{\eta}}\right)\left[-D'_{-\theta}(-\beta)+\dfrac{1}{2}\sqrt{\beta^2+4\theta}
D_{-\theta}(-\beta)\right]\notag
\end{align}
we see that as $\eta\to 0$ the expression in (\ref{green1})
becomes that in~(\ref{277}). We also have
\begin{equation}\label{eq315}
\dfrac{\Mcal}{\Vcal}\to\dfrac{D'_{-\theta}(\beta)+\sqrt{\theta+\beta^2/4}D_{-\theta}(\beta)}{-D'_{-\theta}(-\beta)+\sqrt{\theta+\beta^2/4}D_{-\theta}(-\beta)},\
\eta\to 0^+
\end{equation}
which can be used to obtain the limits of
(\ref{green2}) and~(\ref{green3}), and this agrees with the
results we obtained in~\cite{vlk}. Throughout this
calculation we divided several times by
$D_{-\theta/\eta}\big(\beta/\sqrt{\eta}\big)$,
which was permissible since, for $\theta>0$, we are
outside of the oscillatory range of the special
function, as we discussed in Section~\ref{S5}.

\section{Appendix D}\label{DISCRETE}

Here we discuss the discrete $M/M/m+M$
model. We shall obtain an explicit, albeit
complicated, expression for the Laplace
transform (over time) of $p_n(t)=\text{\upshape Prob}[N(t)=n\mid
N(0)=n_0]$, where $N(t)$ is the number of customers
in the system. Then we will give an
alternate derivation of Theorem~\ref{thmspec}, by evaluating
the discrete model in the limit $m\to\infty$ with $\rho=\lambda/\mu=m
+ \Oup(\sqrt{m})$. The analysis here closely parallels
the proof of Theorem~\ref{greenpos}, so we just give
the main points.

We solve the following infinite system
of ODEs (we assume time is scaled to
make the service rate $\mu=1$, so that $\rho=\lambda$):
\begin{equation}\label{eqA1}
p'_0(t)=p_1(t)-\rho p_0(t)
\end{equation}
\begin{equation}\label{eqA2}
p'_n(t)=\rho\left[p_{n-1}(t)-p_n(t)\right]+(n+1)p_{n+1}(t)-np_n(t),\
1\leq n\leq m-1,
\end{equation}
\begin{equation}\label{eqA3}
p'_m(t)=\rho\left[p_{m-1}(t)-p_m(t)\right]+(m+\eta)p_{m+1}(t)-mp_m(t),
\end{equation}
\begin{align}\label{eqA4}
p_n'(t)={}&\rho\left[p_{n-1}(t)-p_n(t)\right]+\left[m+(n-m+1)\eta\right]p_{n+1}(t)\\
&-\left[m+(n-m)\eta\right]p_n(t),\
n\geq m+1\notag
\end{align}
with the initial condition $p_n(0)=\delta(n,n_0)$. We need
to consider the cases $n_0<m$ and $n_0>m$
separately, as the discrete model has no
analog of the symmetry relation in Proposition~\ref{lemsym}.

\begin{figure}
\psfrag{d}{\footnotesize ${\rm Re}(z)$}
\psfrag{e}{\footnotesize ${\rm Im}(z)$}
\psfrag{f}{\footnotesize $C_1$}
\psfrag{a}{\footnotesize ${\rm Re}(z)$}
\psfrag{b}{\footnotesize ${\rm Im}(z)$}
\psfrag{c}{\footnotesize $C_2$}
\begin{center}
 \includegraphics[width= .5 \linewidth]{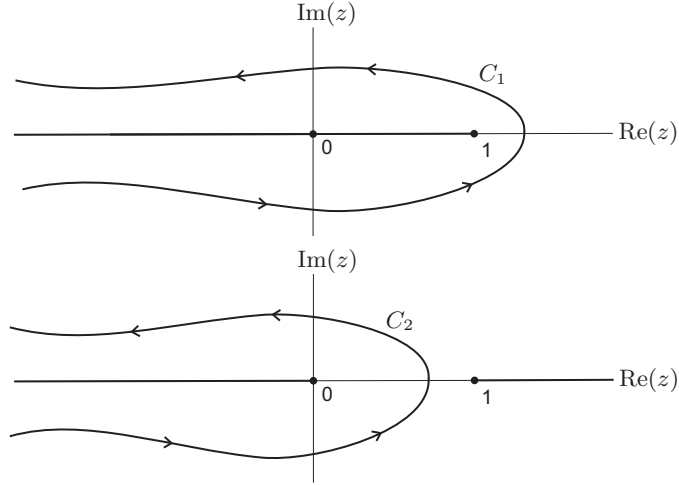}
 \end{center}
  \caption{A sketch of the branch cuts and the contours $C_1$ and $C_2$.}
  \label{figAD}
\end{figure}

Introducing the Laplace transform $\widehat{p}_n(\theta)=
\int^{\infty}_0e^{-\theta t}p_n(t)\, dt$, we first consider the case
$0<n_0<m$ and then (\ref{eqA1})--(\ref{eqA4}) become
\begin{equation}\label{eqA5}
0=\widehat{p}_1 -(\rho+\theta)\widehat{p}_0
\end{equation}
\begin{equation}\label{eqA6}
-\delta(n,n_0)=\rho
\widehat{p}_{n-1}-\left(\rho+\theta+n\right)\widehat{p}_n+(n+1)\widehat{p}_{n+1},\
1\leq n\leq m-1
\end{equation}
\begin{equation}\label{eqA7}
0=\rho\widehat{p}_{n-1}-\left[\rho+\theta+m+(n-m)\eta\right]\widehat{p}_n+\left[m+(n-m+1)\eta\right]
\widehat{p}_{n+1},\ n\geq m.
\end{equation}
We solve (\ref{eqA5})--(\ref{eqA7}) using a discrete Green's
function approach. We begin by introducing
the functions $F_n$, $G_n$, $H_n$, $I_n$; these are
defined by the contour integrals
\begin{align}\label{eqA8}
F_n(\theta)&=\dfrac{1}{2\pi i}\int_{C_0} \dfrac{e^{\rho
z}}{z^{n+1}(1-z)^{\theta}}\, dz\\
&=\sum^n_{\ell=0}\dfrac{\rho^{n-\ell}}{(n-\ell)!}\dfrac{(\theta+\ell-1)(\theta+\ell-2)\dots
(\theta+1)\theta}{\ell!},\notag
\end{align}
\begin{equation}\label{eqA9}
G_n(\theta)=\dfrac{1}{2\pi i}\int_{C_1}\dfrac{e^{\rho
z}}{z^{n+1}(z-1)^{\theta}}\, dz,
\end{equation}
\begin{equation}\label{eqA10}
H_n(\theta)=\dfrac{1}{2\pi i}\int_{C_1}\dfrac{e^{\rho
z/\eta}}{(z-1)^{\theta/\eta}z^{n+1-m}}z^{-m/\eta}\, dz,
\end{equation}
\begin{equation}\label{eqA11}
I_n(\theta)=\dfrac{1}{2\pi i}\int_{C_2}\dfrac{e^{\rho
z/\eta}}{(1-z)^{\theta/\eta}z^{n+1-m}}z^{-m/\eta}\, dz.
\end{equation}
Here $C_0$ is a small loop about $z=0$, $C_1$
goes from $-\infty-i\varepsilon$ to $-\infty+i\varepsilon$, encircling
$z=1$, while $C_2$ goes
from $-\infty-i\varepsilon$ to $-\infty+i\varepsilon$ encircling $z=0$.
The contours $C_1$ and $C_2$ are sketched in
Figure~\ref{figAD}, and in (\ref{eqA9}) we use the branch
$(z-1)^{-\theta}=|z-1|^{-\theta}\exp[-i\theta\arg(z-1)]$ with
$-\pi<\arg(z-1)\leq \pi$; in (\ref{eqA10}) $(z-1)^{-\theta/\eta}=
|z-1|^{-\theta/\eta}\exp\big[-i\theta\eta^{-1}\arg(z-1)\big]$ and
$z^{-m/\eta}
=|z|^{-m/\eta}\exp\big[-i m\eta^{-1}\arg z\big]$ with $-\pi<\arg z\leq
\pi$;
and in (\ref{eqA11})
$(1-z)^{-\theta/\eta}=|1-z|^{-\theta/\eta}\exp\big\{-i\theta\eta^{-1}[\arg(z-1)
-\pi]\big\}$ (so that $1^{-\theta/\eta}=1$ for $\theta>0$). Note
that if $\eta=1$ then $G_n=H_n$ and $F_n=I_n$
(then the contour $C_2$ in (\ref{eqA11}) may be
deformed to $C_0$).

We can easily verify that $F_n$ and $G_n$
satisfy the homogeneous form of (\ref{eqA6}) (with
$\delta(n,n_0)$ replaced by zero) and thus give two
linearly independent solutions of this
difference equation. Similarly, $H_n$ and $I_n$
give two solutions to (\ref{eqA7}). We use the
functions in (\ref{eqA8})--(\ref{eqA11}) to construct
$\widehat{p}_n(\theta)$,
making use of (\ref{eqA3}) (or (\ref{eqA7}) with $n=m$) to relate the
ranges $n<m$ and $n>m$.
The analysis is completely analogous to
the proof of Theorem~\ref{greenpos} in Section~\ref{sec:proofs}, so
we give below only the final result:
\begin{equation}\label{eqA12}
\widehat{p}_n(\theta)=\dfrac{n_0!}{m!}\rho^{m-n_0-1}\dfrac{F_{n_0}H_n}{F_mH_{m-1}-H_mF_{m-1}},\ n\geq m;
\end{equation}
\begin{equation}\label{eqA13}
\widehat{p}_n(\theta)=\dfrac{n_0!\Gamma(\theta)e^{-\rho}}{\rho^{n_0+\theta}}F_{n_0}\left[G_n+\dfrac{H_mG_{m-1}-G_mH_{m-1}}{F_mH_{m-1}-H_mF_{m-1}}F_n\right]\
n_0\leq n\leq m;
\end{equation}
\begin{equation}\label{eqA14}
\widehat{p}_n(\theta)=\dfrac{n_0!\Gamma(\theta)e^{-\rho}}{\rho^{n_0+\theta}}F_n\left[G_{n_0}+\dfrac{H_mG_{m-1}-G_mH_{m-1}}{F_mH_{m-1}-H_mF_{m-1}}F_{n_0}\right],\
0\leq n\leq n_0.
\end{equation}
Here we suppressed the dependence of $F_n$, $G_n$, $H_n$ on $\theta$.

The above holds for all $0\leq n_0\leq m$, and
if $n_0=m$ (starting with all servers occupied
but an empty queue), (\ref{eqA13}) is not needed,
and then $\widehat{p}_n(\theta)$ somewhat simplifies, to
\begin{equation}\label{eqA15}
\widehat{p}_n(\theta)=\dfrac{\rho^{-1}}{F_mH_{m-1}-H_mF_{m-1}}\begin{cases}
F_mH_n,&n\geq m\\
H_mF_n,&0\leq n\leq m.
\end{cases}
\end{equation}
Note that (\ref{eqA12}), (\ref{eqA13}) and (\ref{eqA14}) are similar in
form
to (\ref{green1}), (\ref{green3}) and (\ref{green2}), respectively,
with $H_mG_{m-1}-G_mH_{m-1}$ and $F_mH_{m-1}-
H_mF_{m-1}$ playing the roles of $\Mcal$ and $\Vcal$.
The functions in (\ref{eqA8})--(\ref{eqA11}) are all entire
functions of~$\theta$, and the singularities
of $\widehat{p}_n(\theta)$ are determined by the equation
\begin{equation}\label{eqA16}
F_m(\theta) H_{m-1}(\theta)-H_m(\theta)F_{m-1}(\theta)=0.
\end{equation}

For initial conditions $n_0>m$ we
need to solve (\ref{eqA5})--(\ref{eqA7}), but now with
the left-hand side of (\ref{eqA6}) replaced by zero
and the left-hand side of (\ref{eqA7}) replaced by
$-\delta(n,n_0)$. For $n_0\geq m$ the final result is
\begin{align}\label{eqA17}
\widehat{p}_n(\theta)={}&\dfrac{1}{\rho}
e^{-\rho/\eta}\left(\dfrac{\eta}{\rho}\right)^{n_0-1-m+(\theta+m)/\eta}\Gamma\left(\dfrac{\theta}{\eta}\right)\Gamma\left(n_0+1-m+\dfrac{m}{\eta}\right)\\
&\times\left[I_{n_0}+\dfrac{I_mF_{m-1}-I_{m-1}F_m}{F_mH_{m-1}-H_mF_{m-1}}H_{n_0}\right]H_n,\
n\geq n_0;\notag
\end{align}
\begin{align}\label{eqA18}
\widehat{p}_n(\theta)={}&\dfrac{1}{\rho}e^{-\rho/\eta}\left(\dfrac{\eta}{\rho}\right)^{n_0-1-m+(\theta+m)/\eta}\Gamma\left(\dfrac{\theta}{\eta}\right)\Gamma\left(n_0+1-m+\dfrac{m}{\eta}\right)\\
&\times\left[I_n+\dfrac{I_mF_{m-1}-I_{m-1}F_m}{F_mH_{m-1}-H_mF_{m-1}}H_n\right]H_{n_0},\
m\leq n\leq n_0;\notag
\end{align}
\begin{gather}\label{eqA19}
\widehat{p}_n(\theta)=\dfrac{1}{\rho}\left(\dfrac{\rho}{\eta}\right)^{m-n_0}\dfrac{\Gamma(n_0+1-m+m/\eta)}{\Gamma(1+m/\eta)}\dfrac{H_{n_0}F_n}{F_mH_{m-1}-H_mF_{m-1}},\\
0\leq n\leq m.\notag
\end{gather}
When $n_0=m$, (\ref{eqA18}) is not needed and then
(\ref{eqA17}) and (\ref{eqA19}) reduce to (\ref{eqA15}). Again the
singularities are determined by (\ref{eqA16}).

Now we evaluate these results
in the limit of $m\to\infty$ with the scaling
\begin{equation}\label{eqA20}
n=m+\sqrt{m}x,\ n_0=m+\sqrt{m}x_0,\ \rho=m-\sqrt{m}\beta,
\end{equation}
where $x$, $x_0$, $\beta$ are $\Oup(1)$. We shall thus
give an alternate derivation of Theorems~1
and~8 of Section~\ref{sec:main}. Let us also scale
$z=1-\xi/\sqrt{m}$ in the integrands in (\ref{eqA8})--(\ref{eqA10}).
Noting that
\begin{align*}
\rho z-n\log
z&=(m-\sqrt{m}\beta)\left(1-\dfrac{\xi}{\sqrt{m}}\right)-(m+\sqrt{m}x)\log\left(1-\dfrac{\xi}{\sqrt{m}}\right)\\
&=\rho+(x+\beta)\xi+\dfrac{1}{2}\xi^2+o(1)
\end{align*}
we obtain a limiting form of~(\ref{eqA8}):
\begin{align}\label{eqA21}
F_n(\theta)&\sim\dfrac{m^{\theta/2}e^{\rho}}{\sqrt{2\pi
m}}\dfrac{1}{\sqrt{2\pi}i}\int_{\Br}
\xi^{-\theta}e^{(x+\beta)\xi}e^{\xi^2/2}\, d\xi\\
&=\dfrac{m^{\theta/2}e^{\rho}}{\sqrt{2\pi
m}}e^{-(x+\beta)^2/4}D_{-\theta}(-x-\beta).\notag
\end{align}
Here we used the integral representation
in (\ref{eqP2}) for the parabolic cylinder
function~$D$. In (\ref{eqA21}) $\Br$ is a vertical
contour with ${\rm Re}(\xi)>0$, which can be used
to approximate $C_0$ in (\ref{eqA8}) with this scaling of~$z$.
A completely analogous expansion of~(\ref{eqA9})
leads to
\begin{equation}\label{eqA22}
G_n(\theta)\sim\dfrac{m^{\theta/2}e^{\rho}}{\sqrt{2\pi
m}}e^{-(x+\beta)^2/4}D_{-\theta}(x+\beta),
\end{equation}
and, after some calculation, we obtain from (\ref{eqA10})
\begin{equation}\label{eqA23}
H_n(\theta)\sim\sqrt{\dfrac{\eta}{2\pi
m}}e^{\rho/\eta}\left(\dfrac{m}{\eta}\right)^{\!\!\frac{\theta}{2\eta}}\exp\left[-\dfrac{(\eta
x+\beta)^2}{4\eta}\right]D_{-\theta/\eta}\left(\dfrac{\eta
x+\beta}{\sqrt{\eta}}\right).
\end{equation}
Using (\ref{eqA10}) and Stirling's formula we also have
\begin{equation}\label{eqA24}
\dfrac{n_0!}{m!}\rho^{m-n_0-1}\sim\dfrac{1}{m}e^{\beta x_0}e^{x^2_0/2}.
\end{equation}
Next we consider the limiting form of (\ref{eqA16}).
Noting that $F_m-F_{m-1}$ can be computed
from (\ref{eqA8}) by multiplying the integrand by
$1-z$ and setting $n=m$, we have
\begin{align}\label{eqA25}
F_mH_{m-1}&-F_{m-1}H_m\\
={}&
H_m(F_m-F_{m-1})-F_m(H_m-H_{m-1})\notag\\
\sim{}&H_m\dfrac{e^{\rho}m^{\theta/2}}{\sqrt{2\pi}m}e^{-\beta^2/4}D_{1-\theta}(-\beta)\notag\\
&+F_m\dfrac{e^{\rho/\eta}}{\sqrt{2\pi
m}}\sqrt{\eta}\left(\dfrac{m}{\eta}\right)^{\frac{\theta-\eta}{2\eta}}e^{-\beta^2/(4\eta)}
D_{1-\theta/\eta}\left(\dfrac{\beta}{\sqrt{\eta}}\right)\notag\\
\sim{}&\dfrac{1}{2\pi}m^{\theta/2}\left(\dfrac{m}{\eta}\right)^{\frac{\theta}{2\eta}}e^{\rho(1+1/\eta)}e^{-\beta^2/4}e^{-\beta^2/(4\eta)}\dfrac{\sqrt{\eta}}{m^{3/2}}\Vcal(\theta;\eta,\beta)\notag
\end{align}
where $\Vcal$ is as in Theorem~\ref{thmspec}. Here we
also used (\ref{eqA21}), (\ref{eqA23}) and (\ref{eqP9}).

An analogous approximation to $H_mG_{m-1}-G_mH_{m-1}$
can be used to show that
\begin{equation}\label{eqA26}
\dfrac{H_mG_{m-1}-G_mH_{m-1}}{F_mH_{m-1}-H_mF_{m-1}}\to\dfrac{\Mcal(\theta;\eta,\beta)}{\Vcal(\theta;\eta,\beta)}.
\end{equation}
Then using (\ref{eqA21})--(\ref{eqA24}) and (\ref{eqA26}), we see
that the expressions in (\ref{eqA12})--(\ref{eqA14}) reduce
to those in (\ref{green1})--(\ref{green3}), up to a factor of $1/\sqrt{m}$ in the
former, which
arises due to the fact that $p_n(t)$ is
normalized by a sum over $n$
while $p(x,t)$ is normalized by an integral
over~$x$. We have thus given an
alternate derivation of Theorem~\ref{greenpos}.
Note that all of the asymptotic calculations do
not involve scaling time~$t$ or the transform
variable~$\theta$.

Finally, we discuss the uniformity of the approximation in \eqref{eqA25} for small values of $\theta$. This will show that $\Delta(\theta)=\Delta(\theta;m,\beta,\eta)\equiv F_m H_{m-1}-F_{m-1}H_m$ can have no roots for $m\to\infty$
in ranges where $\theta=o(1)$. We set $\rho=m-\beta\sqrt{m}$ and consider finite intervals of $\beta$ and $\eta$, with $\eta>0$. First we note that, using \eqref{eqA8} and \eqref{eqA10} (with $n=m$ and then $n=m-1$)
 \begin{equation}\label{ad1}
F_m(0)=\frac{\rho^m}{m!}, \quad H_m(0)=\Big(\frac{\rho}{\eta}\Big)^{m/\eta}\frac{1}{\Gamma(m/\eta+1)}, \quad
H_{m-1}(0)=\Big(\frac{\rho}{\eta}\Big)^{m/\eta-1}\frac{1}{\Gamma(m/\eta)}.
\end{equation}
It follows that $F_m(0) H_{m-1}(0)=F_{m-1}(0)H_m(0)$ and thus $\Delta(0;m,\beta,\eta)=0$ for all values of $m,\beta,\eta$. The pole at $\theta=0$ corresponds to the steady state limit, which exists for all $m$ and $\beta$, for $\eta>0$. We expand $\Delta$ for fixed finite $\beta$ and fixed $\eta>0$, as $m\to\infty$, which will refine the leading order result in \eqref{eqA25} and show that the higher order terms remain smaller than the leading term, for $\theta=o(1)$ as $m\to\infty$. Setting $z=1-\xi/\sqrt{m}$ in the integral in \eqref{eqA8} leads to
\begin{align}\label{ad2}
m^{-\theta/2}F_m(\theta)&=\frac{\e^\rho}{\sqrt{m}}\int_{C'} \e^{\xi^2/2}\e^{\beta\xi}\xi^{-\theta}\mathcal{F}(\xi;m)d\xi,
\end{align}
where $C'$ is the image of the contour $C_0$ and
\begin{align}\label{ad3}
\mathcal{F}(\xi;m)&=\Big(1-\frac{\xi}{\sqrt{m}}\Big)^{-1}\exp\Big[-\sqrt{m}\xi-\frac{\xi^2}{2}-m\log \Big(1-\frac{\xi}{\sqrt{m}}\Big)\Big]\nonumber\\
&=\Big(1-\frac{\xi}{\sqrt{m}}\Big)^{-1}\exp\Big[\sum_{l=3}^\infty \frac{\xi^l}{l(\sqrt{m})^{l-2}}\Big]\nonumber\\
&=1+\sum_{j=1}^\infty P_j(\xi)m^{-j/2},
\end{align}
where $P_j(\xi)$ is a polynomial in $\xi$ of degree $3j$. We have $P_1(\xi)=\xi+\frac13 \xi^3$, $P_2(\xi)=\xi^2+\frac{7}{12} \xi^4+\frac{1}{18}\xi^6$, etc. Using \eqref{ad3} and some contour deformation (as $|z|<1$ in \eqref{eqA8} implies asymptotically that ${\rm Re}(\xi)>0$), we obtain the asymptotic series
\begin{align}\label{ad4}
m^{-\theta/2}F_m(\theta)&\sim\frac{\e^\rho}{\sqrt{m}}\Big[\frac{1}{\sqrt{2\pi}}\e^{-\beta^2/4}D_{-\theta}(-\beta)+\sum_{j=1}^\infty m^{-j/2} f_j(\theta,\beta)\Big],
\end{align}
where
\begin{align}\label{ad5} f_j(\theta,\beta)=\frac{1}{2\pi i}\int_{{\rm Br}_+} P_j(\xi)\e^{\xi^2/2}\e^{\beta\xi}\xi^{-\theta}d\xi
\end{align}
and ${\rm Re}(\xi)>0$ on the vertical contour ${{\rm Br}_+}$. Since the $P_j$ are polynomials, the integral in \eqref{ad5} is a finite sum of parabolic cylinder functions of different orders (for example $f_1$ involves $D_{1-\theta}(\cdot)$ and $D_{3-\theta}(\cdot)$), and thus each $f_j$ is an entire function of $\theta$. A completely analogous calculation shows that
 \begin{align}\label{ad6}
m^{-\theta/2}[F_m(\theta)-F_{m-1}(\theta)]&\sim\frac{\e^\rho}{m}\Big[\frac{1}{\sqrt{2\pi}}\e^{-\beta^2/4}D_{1-\theta}(-\beta)+\sum_{j=1}^\infty m^{-j/2} f_j(\theta-1,\beta)\Big].
\end{align}

Next consider $H_m(\theta)$ (setting $n=m$ in \eqref{eqA10}) in the same asymptotic limit. Now we scale $z=1+\xi\sqrt{\eta/m}$ and obtain
 \begin{align}\label{ad7}
\Big(\frac{m}{\eta}\Big)^{-\theta/(2\eta)}H_m(\theta)&=\sqrt{\frac{\eta}{m}}\e^{\rho/\eta}\int_{C_1'} \e^{\xi^2/2}\e^{\beta\xi}\xi^{-\theta}\mathcal{H}(\xi;m)d\xi,
\end{align}
with
\begin{align}\label{ad8}
\mathcal{H}(\xi;m)&=\Big(1+\sqrt{\frac{\eta}{m}}\xi\Big)^{-1}\exp\Big[\sqrt{\frac{m}{\eta}}\xi-\frac{\xi^2}{2}-m\log \Big(1+\sqrt{\frac{\eta}{m}}\xi\Big)\Big]\nonumber\\
&=1+\sum_{j=1}^\infty \tilde{P}_j(\xi;\eta)m^{-j/2},
\end{align}
where $\tilde{P}_j$ are again polynomials in $\xi$. Then we obtain the asymptotic expansion of \eqref{ad7} as
 \begin{align}\label{ad9}
\Big(\frac{m}{\eta}\Big)^{-\theta/(2\eta)}H_m(\theta)&\sim\sqrt{\frac{\eta}{m}}\e^{\rho/\eta}
\Big[\frac{1}{\sqrt{2\pi}}\e^{-\beta^2/(4\eta)}D_{-\theta/\eta}(\tfrac{\beta}{\sqrt{\eta}})+\sum_{j=1}^\infty m^{-j/2} h_j(\theta,\beta,\eta)\Big].
\end{align}
where
\begin{align}\label{ad10} h_j(\theta,\beta,\eta)=\frac{1}{2\pi i}\int_{{\rm Br}_+} \tilde{P}_j(\xi,\eta)\e^{\xi^2/2}\e^{-\beta\xi/\sqrt{\eta}}\xi^{-\theta/\eta}d\xi
\end{align}
so again $h_j$ is a finite sum of parabolic cylinder functions of different indices, all with argument $-\beta/\sqrt{\eta}$.
A completely  analogous calculation shows that
 \begin{align}\label{ad11}
\Big(\frac{m}{\eta}\Big)^{-\theta/(2\eta)}[H_{m-1}(\theta)-H_m(\theta)]&\sim\frac{\eta}{m}\e^{\rho/\eta}
\Big[\frac{1}{\sqrt{2\pi}}\e^{-\beta^2/(4\eta)}D_{1-\theta/\eta}(-\tfrac{\beta}{\sqrt{\eta}})+\sum_{j=1}^\infty m^{-j/2} h_j(\theta-\eta,\beta,\eta)\Big].
\end{align}
Combining \eqref{ad4}-\eqref{ad6} with \eqref{ad9}-\eqref{ad11} leads to a refinement of \eqref{eqA25} into a full asymptotic expansion of
$\Delta(\theta;m,\beta,\eta)$ in powers of $m^{-1/2}$. Then we can replace $\funcDD$ in \eqref{eqA25} by
\begin{align}\label{ad12}
\funcDD(\theta;\eta,\beta)+\sum_{j=1}^\infty m^{-j/2} \funcDD_j(\theta;\eta,\beta)
\end{align}
where
\begin{align}\label{ad13}
\funcDD_j(\theta;\eta,\beta)=\sum_{k=0}^j[h_k(\theta,\beta,\eta)f_{j-k}(\theta-1,\beta)+h_k(\theta-\eta,\beta,\eta)f_{j-k}(\theta,\beta)].
\end{align}
Here $\funcDD_0=\funcDD$ and $f_0$ and $h_0$ can be identified from \eqref{ad4} and \eqref{ad9}. For \eqref{ad13} we see that each $\funcDD_j$ is a finite sum of products of two parabolic cylinder functions, possibly of different indices and arguments. Thus each $\funcDD_j$ is an entire function of $\theta$, and the expansion in \eqref{ad12} is uniform in finite $\beta, \theta,\eta$ intervals, with $\eta>0$. Setting $\theta=0$ and using the fact that $\Delta(\theta)$, which is itself and entire function of $\theta$, vanishes for all values of $m,\beta,\eta$, we conclude that $\funcDD_j(0;\eta,\beta)=0$ for each $j$. From \eqref{eqM10}, $\funcDD$ has a simple zero at $\theta=0$ and then $\funcDD_j$ for $j\geq 1$ have zeros of orders $\geq 1$. Thus for any $\theta=o(1)$ as $m\to\infty$ the correction terms in \eqref{ad12} remain smaller than the leading term, and thus $\Delta(\theta)$ cannot have a zero for $\theta=o(1)$, except the one at $\theta=0$.
We have thus shown that for fixed $\beta$, fixed $\eta>0$, and $m\to\infty$ the zeros of $\Delta$ in \eqref{eqA25} in the range $\theta=\Theta(1)$ must approach the zeros of $\funcDD$, and that $\Delta$ has no zeros where $\theta=o(1)$, except for $\theta=0$.

\section{Appendix E}\label{appendixE}

Here we establish the positivity of the correction term in \eqref{13thirt}. We thus take $\beta>\beta_*$ and show that
$\mathcal A (\beta)>0$.

We begin by noting that $\tilde{\mathcal{V}} (0,\beta)=0$ and $\tilde{\mathcal{V}} (r_0(\beta),\beta)=0$, and since $r_0(\beta)$ is the minimal positive solution of
 $\tilde{\mathcal{V}} (p,\beta)=0$ it follows that $0$ and $r_0(\beta)$ are consecutive zeros of  $\tilde{\mathcal{V}} (p,\beta)$. These zeros are necessarily simple, due to the general results for the one-dimensional Schr\"{o}dinger equations that we discussed in Section \ref{S1} (see the discussion surrounding \eqref{eqS8}). It follows that
 \begin{equation}\label{E1}
 {\rm sgn}\Big(\frac{\partial \tilde{\mathcal{V}}}{\partial p}\Big|_{p=0}\Big)= -{\rm sgn}\Big(\frac{\partial \tilde{\mathcal{V}}}{\partial p}\Big|_{p=r_0(\beta)}\Big).
 \end{equation}
The right-hand side of \eqref{E1} appears in $\mathcal A (\beta)$ in Proposition \ref{lem2}, but the left-hand side is much easier to determine. Using \eqref{eq133} and the definition of $\tilde{\mathcal{V}}$ in Proposition \ref{lem2}, we have
 \begin{align}\label{E2}
\frac{\partial \tilde{\mathcal{V}}}{\partial p}\Big|_{p=0}&=\frac{\partial}{\partial p}[D_p'(-\beta)-D_p(-\beta)\sqrt{\beta^2/4-p}]\Big|_{p=0}\\
&=D_0(-\beta)\frac{\partial}{\partial p}[-\sqrt{\beta^2/4-p}+D_p'(-\beta)/D_p(-\beta)]\Big|_{p=0}\nonumber\\
&=\Big[\frac{1}{\beta}+\e^{\beta^2/2}\int_{-\infty}^\beta \e^{-u^2/2}du\Big]\e^{-\beta^2/4}>0.\nonumber
 \end{align}
Here we also used the facts that $\beta>0$ and $\tilde{\mathcal{V}}(0,\beta)=0$. From \eqref{E1} and \eqref{E2} we conclude that
 \begin{equation}\label{E3}
{\rm sgn}\Big(\frac{\partial \tilde{\mathcal{V}}}{\partial p}\Big|_{p=r_0(\beta)}\Big)=-1.
 \end{equation}
If we can show that $D_{r_0}(-\beta)<0$ for all $\beta>\beta_*$ then \eqref{AAA} shows that $\mathcal A (\beta)>0$.

Consider $\mathcal F (\beta)=D_{r_0(\beta)}(-\beta)$ as a function of $\beta$. As discussed in Section \ref{sec63}, $\mathcal F (\beta)$ is an infinitely smooth function of $\beta$. This function cannot change sign, for if $\mathcal F (\tilde \beta_c)=0$ for some $\tilde \beta_c$ then $D_p(-\beta)=0$ for $p=r_0(\tilde \beta_c)$ and also, since $\tilde{\mathcal{V}} (r_0(\beta), \beta)=0$, $D_p'(-\beta)=0$. But $D_p$ and $D_p'$ cannot simultaneously
vanish, as discussed in Section \ref{S5} below equation \eqref{eqP9A}. To determine the constant sign of $\mathcal F (\beta)$ we need only evaluate this at one particular point. For example, when $\beta=2$ then $r_0\approx .9323$ and $D_{r_0(\beta)}(-\beta)\approx -.8275$. We can also let $\beta\to \infty$ and use the facts that $r_0(\beta)\to 1$ and $D_{r_0(\beta)}(-\beta)\sim -\beta\e^{-\beta^2/4}<0$. Thus ${\rm sgn}(D_{r_0}(\beta))=-1$, and
\eqref{AAA} and \eqref{E3} show that ${\rm sgn}(\mathcal A (\beta))=+1$.

\nocite{*}
\bibliographystyle{alpha}

\begin{thebibliography}{99}



\bibitem{a&s} M. Abramowitz and I.~A. Stegun. \textit{Handbook of Mathematical Functions} (10th printing), 1972.



\bibitem{AtC}
J.~D. Atkinson and T.~K. Caughley.
\newblock Spectral density of piecewise linear first order systems excited by
  white noise.
\newblock {\em Int. J. Non-Linear Mechanics}, 3:137--156, 1968.

\bibitem{At67}
J.~D. Atkinson.
\newblock {\em Spectral density of first order piecewise linear systems excited
  by white noise}.
\newblock PhD thesis, CalTech, 1967.


\bibitem
{blancvandoorn}
J.~P.~C. Blanc and E.~A. van Doorn. Relaxation times for queueing systems. In \textit {Mathematics and Computer Science} (eds. J.W. de Bakker, M. Hazewinkel, J.K. Lenstra), North-Holland, Amsterdam, 139-162, 1984.




\bibitem{BH}
N.~Bleistein and R.~A. Handelsman.
\newblock {\em Asymptotic Expansions of Integrals}.
\newblock Dover, New York, 1986.

%



\bibitem{CL}
E.~A. Coddington and N.~Levinson.
\newblock {\em Theory of Ordinary Differential Equations}.
\newblock McGraw--Hill, New York, 1955.


\bibitem{cohen} J.~W. Cohen. \textit{The Single Server Queue}.
North Holland, Amsterdam, 1982.

\bibitem{daitezcan}
J.~G. Dai, S. He and T. Tezcan. Many-server diffusion limits for $G/Ph/n+M$ queues. Preprint, 2008.


\bibitem{vandoorn} E.~A. van Doorn. Conditions for exponential ergodicity and bounds for the decay parameter of a birth-death process. \textit{Advances in Applied Probability} 17:514-530, 1985.


\bibitem{erdelyi} A. Erdelyi. \textit{Higher Transcendental Functions}. Vol.~2. MacGraw-Hill, New York, 1953.


\bibitem{FS}
P.~Flajolet and R.~Sedgewick.
\newblock {\em Analytic Combinatorics}.
\newblock Cambridge University Press, Cambridge, 2009.

\bibitem{robert} C. Fricker, Ph. Robert and D. Tibi. On the rates of convergence of Erlang's model. \textit{Journal of Applied Probability} 36:1167-1184, 1999.



\bibitem{gans}  N. Gans, G. Koole and A. Mandelbaum. Telephone call centers:  Tutorial, review and research prospects.  \textit{Manufacturing and Service Operations Management} 5:79-141, 2003.

\bibitem{goldberggamarnik}
D.~Gamarnik and D.A.~Goldberg. On the exponential rate of convergence to stationarity in the Halfin-Whitt regime I: The spectral gap of the $M/M/n$ queue. Preprint, 2008.


\bibitem{garnett}
O. Garnett, A. Mandelbaum and M. Reiman. Designing a call center with impatient customers. \textit{Manufacturing \& Service Operations Management} 4:208-227, 2002.

\bibitem{gradshteyn} I.~S. Gradshteyn and I.~M.  Ryzhik. \textit{Table of Integrals, Series and Products}. 5th ed., Academic Press, New York, 1994.

\bibitem{halfinwhitt} S. Halfin and W. Whitt. Heavy-traffic limits for queues with many exponential servers. \textit{Operations Research} 29:567-588, 1981.

    \bibitem{IG} D. Iglehart. Limiting diffusion approximations for the many server queue and the repairman problem. \textit{Journal of Applied Probability} 2:429-441, 1965.



\bibitem{kangramanan} W. Kang and K. Ramanan.
Fluid limits of many-server queues with reneging. {\it Ann. Appl. Prob.} 20:2204-2260, 2010.

 \bibitem{km} S. Karlin and J.~L. McGregor.
Many server queueing processes with Poisson input and exponential service times.
{\it Pacific J. Math.} 8:87-118, 1958.


 \bibitem{keilson} J. Keilson.
A review of transient behavior in regular diffusion and birth-death processes.
{\it J. Applied. Prob.} 1:247-266, 1964.

\bibitem{LRH}
J.~Lehmann, P.~Reimann, and P.~H{\"a}nggi.
\newblock Surmounting oscillating barriers: Path integral approach for weak
  noise.
\newblock {\em Phys. Rev. E}, 62:6282--6303, 2000.

\bibitem{vlk} J.~S.~H. van Leeuwaarden and C. Knessl.
Transient analysis of the Halfin-Whitt diffusion.
{\it Stochastic Processes and their Applications} 121: 1524-1545, 2011.




\bibitem{linetsky} V. Linetsky.
On the transition densities for reflected diffusions.
\textit{Adv. Appl. Prob.}
37:435-460, 2005.


\bibitem{maglaraszeevi} C. Maglaras  and  A. Zeevi.
Diffusion approximations for a multiclass Markovian service
 system with ``guaranteed'' and ``best-effort'' service levels.
\textit{Math. Oper. Res.} 29:786--813, 2004.
	
\bibitem{MP}
A.~N. Malakhov and A.~L. Pankratov.
\newblock Exact solution of {K}ramers' problem for piecewise parabolic
  potentials.
\newblock {\em Physica A}, 229:109--126, 1996.






\bibitem{reed}
J. Reed. The G/GI/N queue in the Halfin-Whitt regime. {\it Ann. Appl. Prob.}
19: 2211-2269, 2009.


\bibitem{REID}
W.~T. Reid.
\newblock {\em Sturmian Theory for Ordinary Differential Equations}.
\newblock Springer-Verlag, New York--Berlin, 1980.

\bibitem{STAK}
I.~Stakgold.
\newblock {\em Boundary Value Problems of Mathematical Physics}, volume~I.
\newblock MacMillan, New York, 1967.

\bibitem{SPA}
W.~Szpankowski.
\newblock {\em Average Case Analysis of Algorithms on Sequences}.
\newblock Wiley-Interscience, New York, 2001.


\bibitem{TEMME}
N.~M. Temme.
\newblock Parabolic cylinder function.
\newblock In R.~F. Boisvert et~al., editors, {\em {NIST} Handbook of
  Mathematical Functions}. Cambridge University Press, 2010.
\newblock ISBN 978-0521192255.

\bibitem{TMa}
E.~C. Titchmarsh.
\newblock On the discreteness of the spectrum associated with certain
  differential equations.
\newblock {\em Ann. Math. Pure Applied}, 28:141--147, 1949.

\bibitem{TMb}
E.~C. Titchmarsh.
\newblock On the discreteness of the spectrum of a differential equation.
\newblock {\em Acta Sci. Math. Szeged}, 12:16--18, 1950.

\bibitem{T1962}
E.~C. Titchmarsh.
\newblock {\em Eigenfunction Expansions Associated with Second-order
  Differential Equations, {Part I}}.
\newblock Clarendon Press, Oxford, second edition, 1962.

\bibitem{wardglynn}
A. Ward and P. Glynn. Properties of the reflected Ornstein-Uhlenbeck process. \textit{Queueing Systems}
44:109-123, 2003.


\bibitem{whitt04}
W. Whitt. Efficieny-deriven heavy-traffic approximations for many-server queues with abandonments. \textit{Management Science}
50:1449-1461, 2004

\bibitem{whitt05}
W. Whitt. Heavy-traffic limits for the $G/H_2^*/n/m$ queue. \textit{Mathematics of Operations Research}
30:1-27, 2006.

\bibitem{whitt06}
W. Whitt. Fluid limits for many-server queues with abandonments. \textit{Operations Research}
54:363-372, 2006.


\bibitem{WONG}
R.~Wong.
\newblock {\em Asymptotic Approximation of Integrals}.
\newblock Academic Press, Inc., Boston, 1989.

\bibitem{XieKnessl1993}
S.~Xie and C.~Knessl.
\newblock On the transient behavior of the {E}rlang loss model: Heavy usage
  asymptotics.
\newblock {\em SIAM J. Appl. Math.}, 53:555--599, 1993.


\bibitem{zeltyn04}
 S. Zeltyn and A. Mandelbaum. The impact of customers' patience on delay and abandonment: some empirically-driven experiments with the $M/M/n +G$ queue. {\it OR Spectrum} 26:377-411, 2004.

\bibitem{zeltyn05}
 S. Zeltyn and A. Mandelbaum. Call centers with impatient customers: many-server asymptotics of the $M/M/n+G$ queue. {\it Queueing Systems} 51:361-402, 2005.


\end{thebibliography}
\begin{footnotesize}

\end{footnotesize}
\end{document}